\documentclass[a4paper,12pt]{amsart}
\usepackage{amsmath,amsthm,amssymb,bm,mathrsfs}

\numberwithin{equation}{section}

\newtheorem{mainthm}{Theorem}
\newtheorem{thm}{Theorem}[section]
\newtheorem{cor}[thm]{Corollary}
\newtheorem{lem}[thm]{Lemma}
\newtheorem{prop}[thm]{Proposition}

\theoremstyle{definition}

\newtheorem{rem}[thm]{Remark}
\newtheorem{ex}[thm]{Example}
\newtheorem{conj}{Conjecture}

\newcommand{\trace}{\operatorname{Tr}}
\newcommand{\Gal}{\operatorname{Gal}}
\newcommand{\Aut}{\operatorname{Aut}}
\newcommand{\End}{\operatorname{End}}
\newcommand{\Zero}{\operatorname{Zero}}
\newcommand{\sgn}{\operatorname{sgn}}

\newcommand{\diag}{\operatorname{diag}}
\newcommand{\GL}{\operatorname{GL}}
\newcommand{\SL}{\operatorname{SL}}
\newcommand{\Cl}{\operatorname{Cl}}
\newcommand{\Tor}{\operatorname{Tor}}
\newcommand{\Disc}{\operatorname{Disc}}

\newcommand{\mmod}[1]{\operatorname{\:\:(mod^* \, #1)}}
\newcommand{\proj}{\operatorname{pr}}
\newcommand{\CInv}{R^{\mathrm{inv}}}
\newcommand{\SSV}{V^{\mathrm{ss}}}
\newcommand{\HATSSV}{\hat{V}^{\mathrm{ss}}}
\newcommand{\ct}{\operatorname{ct}}
\newcommand{\legend}[2]{\left(\frac{#1}{#2}\right)}

\newcommand{\Proof}{{\it Proof.}\quad}

\title[A conjecture on the zeta functions]{A conjecture on the zeta functions of pairs of ternary quadratic forms}
\author{J. Nakagawa}
\address{Department of Mathematics, Joetsu University of Education, 
Joetsu 943-8512, Japan}
\email{jin@juen.ac.jp}
\subjclass[2010]{11E12, 11R16}
\keywords{order, quartic field, quadratic form,  prehomogeneous vector space.}

\begin{document}
\maketitle

\section{Introduction}
For any ternary quadratic form  
\[ x(v)=x_{11}v_1^2+x_{12}v_1v_2+x_{13}x_1v_3+x_{22}v_2^2+x_{23}v_2v_3+x_{33}v_3^2 \]
in three variables $v=(v_1,v_2,v_3)$, we set
\[ P(x)=4x_{11}x_{22}x_{33}+x_{12}x_{13}x_{23}-x_{11}x_{23}^2-x_{22}x_{13}^2-x_{33}x_{12}^2. \]
The action of $g_1 \in \GL(3)$ on $x$ is defined by $(g_1x)(v)=x(vg_1)$.
Then we have $P(g_1x)=(\det g_1)^2P(x)$.
We denote by $V$ the vector space of pairs of ternary quadratic forms
and put $G=\SL(3) \times \GL(2)$. Then the action of $g=(g_1,g_2) \in G$
on $x=(x_1,x_2) \in V$ is defined by
\[ g x=(p(g_1x_1)+q(g_1x_2), r(g_1x_1)+s(g_1x_2)), \]
where $g_2=\left(
\begin{array}{cc}
p & q \\
r & s
\end{array}
\right)$.  For any pair $x=(x_1,x_2) \in V$, 
we define a binary cubic form $F_x(u)$ in two variables $u=(u_1,u_2)$ by
\[ F_x(u)=P(u_1x_1-u_2x_2). \]
Further we denote by $\Disc(x)$ the discriminant $\Disc(F_x)$ of the binary cubic form $F_x(u)$.
The action of $g_2 \in \GL(2)$ on a binary cubic form $F(u)$ is defined by 
$(g_2 F)(u)=(\det g_2)^{-1}F(ug_2)$. 
Then we have 
\[ \Disc(g x)=(\det g_1)^8(\det g_2)^6  \Disc(x). \]
It is known that $(G,V)$ is a regular prehomogeneous vector space and
$\Disc(x)$ is its fundamental relative invariant. We set
\[ \SSV=\{ x \in V  \,|\, \Disc(x) \neq 0 \}. \]
Wright and Yukie proved in \cite{WY} that for any field $k$, 
there exists one to one correspondence between the set of rational orbits $G_k \backslash \SSV_k$ 
and the set of isomorphism classes of \'{e}tale quartic algebras over $k$.  
In this paper, we restrict ourselves to the case of $k=\mathbb{Q}$. 
We denote by $V_\mathbb{Z}$ the set of pairs of integral ternary quadratic forms. 
We study the orbits $\Gamma x$ of $x \in \SSV_\mathbb{Z}=\SSV\cap V_\mathbb{Z}$, 
where we set $\Gamma=\Gamma_1 \times \Gamma_2$, $\Gamma_1=\SL_3(\mathbb{Z})$, 
$\Gamma_2=\GL_2(\mathbb{Z})$.
For any $x=(x_1,x_2) \in \SSV_\mathbb{Z}$, 
the set of zeros of $x$ is defined by
\[ \Zero(x)=\{\xi=(\xi_1:\xi_2:\xi_3) \in \mathbb{P}^2_{\bar{\mathbb{Q}}} \,|\,
   x_1(\xi)=x_2(\xi)=0 \}. 
\]
Then $\Zero(x)$ is a finite set consisting of four points.
Let $y=(y_1,y_2)$ be a pair of symmetric matrices of degree three 
with coefficients in $\mathbb{Z}$.
Then we obtain a pair of ternary quadratic forms $(y_1(v),y_2(v))$ by setting
$y_i(v)=vy_i { }^t\!v$. 
We denote by $\hat{V}_\mathbb{Z}$ the subset of 
$V_\mathbb{Z}$ consisting of pairs of ternary quadratic forms 
which are obtained from pairs of symmetric matrices of degree three 
with coefficients in $\mathbb{Z}$. 
We put $\HATSSV_\mathbb{Z}=\SSV\cap \hat{V}_\mathbb{Z}$. 
By a ring of rank $n$ we mean a commutative ring with unit that is free of rank $n$ 
as a $\mathbb{Z}$-module. 
M. Bhargava proved in \cite{BH3} that quartic rings are parametrized by the 
set $\Gamma\backslash \SSV_\mathbb{Z}$. 
He also proved in \cite{BH2} that the 2-torsion subgroups of the ideal class groups
of cubic rings are parametrized by the set $\Gamma\backslash \HATSSV_\mathbb{Z}$. 
Further he obtained the density theorems of discriminants of quartic rings and fields in \cite{BH4}. 


We overview the historical background of these subjects. 
There is a discriminant preserving bijection between the set of $\GL_2(\mathbb{Z})$-equivalence classes of 
integral binary cubic forms and the set of isomorphism classes of cubic rings. 
This correspondence is called Delone-Faddeev correspondence, 
which were already essentially contained in the work of Levi \cite{LE}. 
H. Davenport obtained in \cite{DA} asymptotic formulae
for the sums of the class numbers of binary cubic forms
of positive and negative discriminants.  
Using this result, Davenport and Heilbronn obtained in \cite{DH}
the density theorems of the discriminants of cubic fields 
of positive and negative discriminants.  
T. Shintani  made a remarkable contribution 
to the study of class numbers of binary cubic forms 
by applying the theory of prehomogeneous vector spaces 
which was founded by M. Sato in 1960's (cf. \cite{SASH}, \cite{SH}).
Shintani defined in \cite{SH} the zeta functions associated 
with the prehomogeneous vector space of binary cubic forms. 
He introduced four Dirichlet series whose coefficients are 
class numbers of integral binary cubic forms.
Using the theory of prehomogeneous vector spaces, 
he proved that the four Dirichlet series are analytically continued to
meromorphic functions on the whole complex plane and 
satisfy certain functional equations. 
Y. Ohno calculated the first two hundred coefficients 
of all of the four series and presented in \cite{OH} 
a conjecture which states that the two of the four 
Dirichlet series are essentially the same 
as the remaining two series up to some elementary factors. 
The conjecture was proved by the author in \cite{NAK98}.

Taking account of the results in \cite{OH} and \cite{NAK98}, 
it is an interesting problem to find relations between the number of 
$\Gamma$-orbits of $x \in \SSV_\mathbb{Z}$ with $\Disc(x)=n$ 
and that of $y \in \HATSSV_\mathbb{Z}$ with $\Disc(y)=cn$
for any given integer $n \neq 0$, where $c \in \mathbb{Z}$ is some constant,
possibly a power of $2$.

There are three real orbits of $\SSV_\mathbb{R}$ corresponding to
the number of points in $\Zero(x)\cap \mathbb{P}^2_\mathbb{R}$. 
We denote by $V_1$, $V_2$ and $V_3$, the set of $x \in \SSV_\mathbb{R}$ such that 
$|\Zero(x)\cap \mathbb{P}^2_\mathbb{R}|$ equals $4$, $2$ and $0$, respectively. 
For any $x \in \SSV$, we denote by $\Gamma_x$ the isotropy group of $x$ in $\Gamma$. 
Then $\Gamma_x$ is a finite group. We put
\[ \mu(x)=1/|\Gamma_x|, \quad L=V_\mathbb{Z}, \quad 
\hat{L}=\hat{V}_\mathbb{Z}.
\]

For any integral binary cubic form $F(u)$, we denote by $R(F)$ the cubic ring associated with $F(u)$ 
by Delone-Faddeev correspondence. 
For any $x \in L$, $F_x(u)$ is an integral binary cubic form and 
for any $y \in \hat{L}$, $\hat{F}_y(u)=(1/4)F_y(u)$ is an integral binary cubic form. 
Hence $\Disc(y)=2^8 \Disc(\hat{F}_y)$. So we put $\Disc^*(y)=\Disc(\hat{F}_y)=2^{-8}\Disc(y)$.

For $i=1,2,3$ and $n \in \mathbb{Z}$, $n\neq 0$, we set
\[a_i(n)=\sum_{\substack{x \in \Gamma\backslash (L\cap V_i) \\ \Disc(x)=n}} \mu(x),  \qquad 
\hat{a}_i(n)=\sum_{\substack{y \in \Gamma\backslash (\hat{L}\cap V_i) \\ \Disc^*(y)=n}} \mu(y).
\]
Then the zeta functions of the prehomogeneous vector space $(G,V)$ for the lattices $L$ and $\hat{L}$ are defined by
\begin{align*}
\xi_i(L,s) 
 &= \sum_{x \in \Gamma\backslash (L \cap V_i)} \mu(x)|\Disc(x)|^{-s}
    =\sum_{n=1}^\infty \frac{a_i((-1)^{i-1}n)}{n^s}, \\
\xi_i(\hat{L},s) 
 &= \sum_{y \in \Gamma\backslash (\hat{L} \cap V_i)} \mu(y)|\Disc^*(y)|^{-s}
    =\sum_{n=1}^\infty \frac{\hat{a}_i((-1)^{i-1}n)}{n^s} 
\end{align*}
for $i=1,2,3$. 

For any cubic ring $\mathcal{O}$, we put
\[ L(\mathcal{O})=\{x \in L\,|\, R(F_x) \cong \mathcal{O}\},\quad 
   \hat{L}(\mathcal{O})=\{y \in \hat{L}\,|\, R(\hat{F}_y) \cong \mathcal{O}\}.
\]
Further we put $L_i(\mathcal{O})=L(\mathcal{O}) \cap V_i$ and 
$\hat{L}_i(\mathcal{O})=\hat{L}(\mathcal{O}) \cap V_i$ for $i=1,2,3$. 
For any number field $k$, we denote by $\mathcal{O}_k$ the maximal order of $k$. 
In this paper, we shall prove the following relations.
\begin{mainthm}\label{thm:MAIN-S4A4}
Let $k$ be a cubic field and $\mathcal{O}$ be an order of $k$ 
such that the index $(\mathcal{O}_k:\mathcal{O})$ is square free. 
If $\Disc(k)>0$, then we have 
\begin{align*}
\sum_{y \in \Gamma\backslash \hat{L}_1(\mathcal{O})} \mu(y)
 &=\sum_{x \in \Gamma\backslash L_1(\mathcal{O})} \mu(x) 
+\sum_{x \in \Gamma\backslash L_3(\mathcal{O})} \mu(x), \\
\sum_{y \in \Gamma\backslash \hat{L}_3(\mathcal{O})} \mu(y)
 &=3\sum_{x \in \Gamma\backslash L_1(\mathcal{O})} \mu(x)-\sum_{x \in \Gamma\backslash L_3(\mathcal{O})} \mu(x).
\end{align*}
If $\Disc(k)<0$, then we have
\[ \sum_{y \in \Gamma\backslash \hat{L}_2(\mathcal{O})} \mu(y)=2\sum_{x \in \Gamma\backslash L_2(\mathcal{O})} \mu(x).\]
\end{mainthm}
%

\begin{mainthm}\label{thm:MAIN}
If $n$ is a discriminant of a quadratic field, then we have
\begin{align*}
\hat{a}_1(n) &=a_1(n)+a_3(n) \quad (n>0), \\
\hat{a}_2(n) &=2 a_2(n) \quad (n<0), \\
\hat{a}_3(n) &=3a_1(n)-a_3(n) \quad (n>0).
\end{align*}
\end{mainthm}

H. Cohen and F. Thorne studied in \cite{CT} the problem of enumerating quartic fields with fixed cubic resolvent field. 
We prove the theorems by enumerating quartic rings with fixed cubic resolvent ring. 
Based on the theorems above and numerical examples, 
we present the following conjecture:

\begin{conj}\label{conj:ON}
\begin{align*}
\xi_1(\hat{L},s) &=\xi_1(L,s)+\xi_3(L,s), \\
\xi_2(\hat{L},s) &=2\xi_2(L,s), \\
\xi_3(\hat{L},s) &=3\xi_1(L,s)-\xi_3(L,s).
\end{align*}
or equivalently,
\begin{align*}
3\xi_1(\hat{L},s)+\xi_3(\hat{L},s) &= 2 \left(3\xi_1(L,s)+\xi_3(L,s)\right), \\
\xi_2(\hat{L},s) &= 2\xi_2(L,s), \\
\xi_1(\hat{L},s)-\xi_3(\hat{L},s) &= -2 \left(\xi_1(L,s)-\xi_3(L,s)\right).
\end{align*}
\end{conj}

The convergence of the zeta functions was proved by Yukie \cite{YU} in this case and 
by H. Saito \cite{SAI} for general cases. 
By Sato-Shintani \cite[Theorem 2]{SASH},  we have the following functional equation:
\begin{align*}
\lefteqn{\left(
\begin{array}{c}
\xi_1(L,1-s)  \\
\xi_2(L,1-s)  \\
\xi_3(L,1-s) 
\end{array}
\right)}\\
 &=\Gamma(s)^4\Gamma\left(s-\frac{1}{6}\right)^2
\Gamma\left(s+\frac{1}{6} \right)^2
\Gamma\left(s-\frac{1}{4} \right)^2\Gamma\left(s+\frac{1}{4} \right)^2\\
 & \times 2^{8s}3^{6s}\pi^{-12s}
\left(
\begin{array}{ccc}
u^*_{11}(s) & u^*_{21}(s) & u^*_{31}(s)  \\
u^*_{12}(s) & u^*_{22}(s) & u^*_{32}(s)  \\
u^*_{13}(s) & u^*_{23}(s) & u^*_{33}(s)  \\
\end{array}
\right)
\left(
\begin{array}{c}
\xi_1(\hat{L},s)  \\
\xi_2(\hat{L},s)  \\
\xi_3(\hat{L},s) 
\end{array}
\right),
\end{align*}
where $u^*_{ji}(s)$'s are polynomials of 
$q=\exp(\pi\sqrt{-1}s)$ and $q^{-1}$ of degree at most $6$.
However, as far as we know, the explicit determination 
of the matrix $(u^*_{ji}(s))$ is not done yet.
If we know the matrix explicitly and the conjecture is true,
then the functional equation above would be a simple form.

The conjecture was formulated by discussions with Professor Y. Ohno. 
The author would like to express grateful thanks to Professors F. Thorne and T. Taniguchi 
for their helpful comments and warm encouragement. 

The organization of this paper is as follows: 
We study the conductors of cubic rings in Section 2. 
We summarize some basic facts on ideal class groups of cubic rings in Section 3. 
In Section 4 we prove formulae for the left hand side of the equations in Theorem \ref{thm:MAIN-S4A4}
by using Bhargava's parametrization of two-torsion subgroups of ideal class groups of cubic rings. 
We summarize some basic facts on Bhargava's parametrization of quartic rings in Section 5. 
We prove formulae for the right hand side of the equations in Theorem \ref{thm:MAIN-S4A4}
in Section 6 when $k$ is a non-Galois cubic field, in Section 7 when $k$ is a Galois cubic field.
Finally, in Section 8 we prove Theorem \ref{thm:MAIN}. 

\section{Cubic rings}
For a ring $R$ of rank $n$, the discriminant $\Disc(R)$ of $R$ is defined
as the determinant $\det(\trace(\alpha_i\alpha_j))\in \mathbb{Z}$, where $\{\alpha_i\}$ is 
any $\mathbb{Z}$-basis of $R$. We call $R$ \textit{nondegenerate} if $\Disc(R)\neq 0$. 

If $T$ is a ring of rank $n$, then it has a subring $T_m = \mathbb{Z} + mT$ 
for any positive integer $m$. 
Conversely, any non-degenerate ring $R$ can be written as 
$R=T_m$ for a unique maximal $m$ which we call the \textit{content}, 
and for a unique ring $T$, which is then called \textit{primitive}. 
The content of primitive ring is $1$. We denote by $\ct(R)$ the content of $R$.

If $R$ is a nondegenerate ring of rank $n$, then 
$k=R \otimes_\mathbb{Z} \mathbb{Q}$ is an \'{e}tale algebra of degree $n$ over $\mathbb{Q}$, i. e. 
$k=k_1 \oplus \cdots \oplus k_s$, $k_i$'s are number fields and $\dim_\mathbb{Q} k=n$. 
We call $k$ a \textit{cubic algebra} or a \textit{quartic algebra} if $n=3$ or $4$. 
We denote by $\mathcal{O}_k$ the maximal order of $k$. 
Then $\mathcal{O}_k=\mathcal{O}_{k_1} \oplus \cdots \oplus \mathcal{O}_{k_s}$. 
We write $\Disc(k)=\Disc(\mathcal{O}_k)$. 
Any $\mathcal{O}_k$-module $\mathfrak{a}$ of $k$ can be written as 
$\mathfrak{a}=\mathfrak{a}_1 \oplus \cdots \oplus \mathfrak{a}_s$, 
where $\mathfrak{a}_i$ is an $\mathcal{O}_{k_i}$-module of $k_i$. 
Hence $\mathfrak{a}$ is an invertible $\mathcal{O}_k$ ideal if and only if 
each $\mathfrak{a}_i$ is a nonzero fractional ideal of $k_i$. 
We define the \textit{norm} of $\mathfrak{a}$ by $N(\mathfrak{a})=\prod_{i=1}^s N(\mathfrak{a}_i)$. 
In particular, $\mathfrak{a}\subset \mathcal{O}_k$ is an $\mathcal{O}_k$ ideal and 
the index $(\mathcal{O}_k:\mathfrak{a})$ is finite, 
then $\mathfrak{a}$ is an invertible $\mathcal{O}_k$ ideal and 
$N(\mathfrak{a})=(\mathcal{O}_k:\mathfrak{a})$. 
We consider only invertible $\mathcal{O}_k$-ideals. 
If $\mathfrak{a}$ and $\mathfrak{b}$ are invertible $\mathcal{O}_k$-ideals with 
$\mathfrak{a} \subset \mathfrak{b}$, then $\mathfrak{a}=\mathfrak{b}\mathfrak{c}$ for some 
integral invertible $\mathcal{O}_k$-ideal $\mathfrak{c}$. 
For an integral invertible $\mathcal{O}_k$-ideal $\mathfrak{a}$ and a prime number $p$, 
we denote by $\mathfrak{a}_{i,p}$ the $p$-part of $\mathfrak{a}_i$ and 
put $\mathfrak{a}_p=\mathfrak{a}_{1,p} \oplus \cdots \oplus \mathfrak{a}_{s,p}$. 
We call $\mathfrak{a}_p$ the $p$-part of $\mathfrak{a}$. 
It is easy to see that $\mathfrak{a}_p=\mathfrak{a}+p^a \mathcal{O}_k$ for a sufficiently large positive integer $a$. 

Let $\mathcal{O} \subset \mathcal{O}_k$ be an order. 
Since the index $f=(\mathcal{O}_k:\mathcal{O})$ is finite, 
$f\mathcal{O}_k$ is an $\mathcal{O}_k$-ideal contained in $\mathcal{O}$. 
The largest $\mathcal{O}_k$-ideal $\mathfrak{f}$ contained in $\mathcal{O}$ is called 
the \textit{conductor} of $\mathcal{O}$. 
Since $f\mathcal{O}_k \subset \mathfrak{f}$, 
$\mathfrak{f}$ is an invertible $\mathcal{O}_k$-ideal. 
\begin{lem}\label{lem:CONDUCTOR1}
Let $\mathcal{O}' \subset \mathcal{O} \subset \mathcal{O}_k$ be orders 
and denote by $\mathfrak{f}'$ and $\mathfrak{f}$ the conductors of 
$\mathcal{O}'$ and $\mathcal{O}$, respectively. Then $\mathfrak{f}' \subset \mathfrak{f}$. 
If $p$ is a prime number such that $p|N(\mathfrak{f})$ and $p\nmid (\mathcal{O}:\mathcal{O}')$, 
then $\mathfrak{f}'_p=\mathfrak{f}_p$. 
\end{lem}
\Proof
The first statement is clear from the definition of conductor. 
Hence we have $\mathfrak{f}'_p \subset\mathfrak{f}_p$. 
We put $m=(\mathcal{O}:\mathcal{O}')$. 
It follows from $\mathfrak{f} \subset \mathcal{O}$ and $m \mathcal{O} \subset \mathcal{O}'$ that 
$m\mathfrak{f} \subset m\mathcal{O}\subset \mathcal{O}'$, 
hence $m\mathfrak{f} \subset \mathfrak{f}'$. 
Take a sufficiently large positive integer $a$ so that $\mathfrak{f}'_p=\mathfrak{f}'+p^a \mathcal{O}_k$ and 
$\mathfrak{f}_p=\mathfrak{f}+p^a \mathcal{O}_k$. 
Since $p \nmid m$, we have $mb+p^ac=1$ for some $b,c \in \mathbb{Z}$, hence 
$\mathfrak{f}\subset m\mathfrak{f}+p^a\mathcal{O}_k$. 
Thus we have 
$\mathfrak{f}_p=\mathfrak{f}+p^a\mathcal{O}_k\subset m\mathfrak{f}+p^a\mathcal{O}_k \subset \mathfrak{f}'+p^a\mathcal{O}_k=\mathfrak{f}'_p$. 
This proves $\mathfrak{f}'_p=\mathfrak{f}_p$. 
\qed

\begin{lem}\label{lem:CONDUCTOR2}
Let $\mathcal{O}' \subset \mathcal{O} \subset \mathcal{O}_k$ be orders 
and $p$ a prime number such that $p \nmid (\mathcal{O}_k:\mathcal{O})$. 
Let $\mathfrak{g}$ be the largest $\mathcal{O}$-ideal contained in $\mathcal{O}'$ 
and assume that $(\mathcal{O}:\mathfrak{g})$ is a power of $p$. 
Then $\mathfrak{g}\mathcal{O}_k$ is the $p$-part of the conductor $\mathfrak{f}'$ of $\mathcal{O}'$. 
Further we have $\mathfrak{g}\mathcal{O}_k \cap \mathcal{O}=\mathfrak{g}$ and 
$\mathcal{O}/\mathfrak{g} \cong \mathcal{O}_k/\mathfrak{g}\mathcal{O}_k$. 
\end{lem}
\Proof
We put $(\mathcal{O}_k:\mathcal{O})=m$, $(\mathcal{O}:\mathfrak{g})=p^r$. 
Since $m\mathcal{O}_k \subset \mathcal{O}$, we have 
$m\mathfrak{g}\mathcal{O}_k \subset \mathfrak{g}\mathcal{O}=\mathfrak{g}  \subset \mathcal{O}'$. 
Since $\mathfrak{f}'$ is the conductor of $\mathcal{O}'$, 
$m\mathfrak{g}\mathcal{O}_k \subset \mathfrak{f}'$. 
Take a sufficiently large positive integer $a\geq r$ 
so that $\mathfrak{f}'_p=\mathfrak{f}'+p^a \mathcal{O}_k$. 
Since $p \nmid m$, we have $mb+p^ac=1$ for some $b,c \in \mathbb{Z}$, hence 
$\mathfrak{g}\mathcal{O}_k \subset m \mathfrak{g}\mathcal{O}_k+p^a\mathcal{O}_k \subset \mathfrak{f}'+p^a\mathcal{O}_k=\mathfrak{f}'_p$. 
On the other hand, $p^r\mathcal{O}\subset \mathfrak{g}$, hence 
$p^r\mathcal{O}_k \subset \mathfrak{g}\mathcal{O}_k$. 
By the inclusion $\mathfrak{f}' \subset \mathcal{O}'$ and the fact that 
$\mathfrak{f}'$ is an $\mathcal{O}$-ideal, we have 
$\mathfrak{f}' \subset \mathfrak{g} \subset \mathfrak{g}\mathcal{O}_k$. 
Hence 
\[ \mathfrak{f}'_p=\mathfrak{f}'+p^a\mathcal{O}_k \subset \mathfrak{g}\mathcal{O}_k+p^a\mathcal{O}_k \subset \mathfrak{g}\mathcal{O}_k+p^r \mathcal{O}_k=\mathfrak{g}\mathcal{O}_k. \]
So we have proved $\mathfrak{f}'_p=\mathfrak{g}\mathcal{O}_k$. 
In particular, $(\mathcal{O}_k:\mathfrak{g}\mathcal{O}_k)$ is a power of $p$. 
Since $p \nmid m$, $\mathcal{O}_k=\mathfrak{g}\mathcal{O}_k+m\mathcal{O}_k \subset \mathfrak{g}\mathcal{O}_k+\mathcal{O} \subset \mathcal{O}_k$. 
So the natural homomorphism $\mathcal{O}\rightarrow \mathcal{O}_k/\mathfrak{g}\mathcal{O}_k$ 
is surjective and has kernel $\mathfrak{g}\mathcal{O}_k \cap \mathcal{O}$. 
Hence $\mathcal{O}/(\mathfrak{g}\mathcal{O}_k \cap \mathcal{O}) \cong \mathcal{O}_k/\mathfrak{g}\mathcal{O}_k$. 
It remains to show $\mathfrak{g}\mathcal{O}_k \cap \mathcal{O}=\mathfrak{g}$. 
The inclusion $\mathfrak{g}\subset \mathfrak{g}\mathcal{O}_k \cap \mathcal{O}$ is obvious. 
Take any element $\alpha \in \mathfrak{g}\mathcal{O}_k \cap \mathcal{O}$. 
Then $m\alpha \in m\mathfrak{g}\mathcal{O}_k\subset \mathfrak{g}$
and $p^a\alpha \in p^r \mathcal{O} \subset \mathfrak{g}$. 
Hence $mb+p^ac=1$ implies $\alpha \in \mathfrak{g}$, 
so $\mathfrak{g}\mathcal{O}_k \cap \mathcal{O} \subset \mathfrak{g}$. 
This proves $\mathfrak{g}\mathcal{O}_k \cap \mathcal{O}=\mathfrak{g}$. 
\qed

Delone and Faddeev showed in \cite{DF} that there exists a  bijection 
between the set of isomorphism classes of cubic rings and the set of 
$\GL_2(\mathbb{Z})$-equivalence classes of integral binary cubic forms, as follows. 
Given a binary cubic form $F(u) = au_1^3 +bu_1^2u_2 +cu_1u_2^2 +du_2^3$ with $a, b, c, d \in \mathbb{Z}$, 
we associates to $F$ the ring $R(F)$ having $\mathbb{Z}$-basis $\{1, \omega, \theta\}$
and multiplication table
\begin{equation}\label{eq:CUBIC-RING-STR}
\begin{split}
\omega^2 &=-ac+b\omega-a\theta, \\
\theta^2 &=-bd+d\omega-c\theta, \\
\omega\theta &=-ad.
\end{split}
\end{equation}
The discriminant of the binary cubic form $F$ is given by
\[ \Disc(F)=18abcd+b^2c^2-4ac^3-4b^3d-27a^2d^2.\]
Then we have $\Disc(R(F))=\Disc(F)$, so the correspondence above is discriminant preserving. 
The ring $R(F)$ is primitive if and only if the binary cubic form $F$ is primitive, i.e. 
$\gcd(a,b,c,d)=1$. 
\begin{prop}\label{prop:CUBIC-SUBRING}
Let $F$, $F'$ be integral binary cubic forms such that $R(F')\subset R(F)$. 
Then there exists a matrix $\delta \in M_2(\mathbb{Z})$ such that 
$|\det \delta|=(R(F):R(F'))$ and $F'=\delta \cdot F$. 
\end{prop}
\Proof
We put $m=(R(F):R(F'))$ and denote by $a', b', c', d'$ the coefficients of $F'$. 
We denote by $\{1,\omega,\theta\}$ and $\{1,\xi,\eta\}$ 
the $\mathbb{Z}$-basis of $R(F)$ and that of $R(F')$, respectively. 
So the ring structure of $R(F)$ is given by \eqref{eq:CUBIC-RING-STR} and 
that of $R(F')$ is similarly. In particular, $\omega\theta=-ad$, $\xi\eta=-a'd'$. 
Since $R(F') \subset R(F)$, there exists a matrix $\delta \in M_2(\mathbb{Z})$ 
such that $|\det \delta|=m$ and 
${}^t(\xi,\eta)\equiv \delta\, {}^t(\omega,\theta) \pmod{\mathbb{Z}}$. 
By a theorem of elementary divisors, there exist positive integers $f_1$, $f_2$ and 
matrices $\gamma_1, \gamma_2 \in \GL_2(\mathbb{Z})$ such that $f_1|f_2$ and 
$\delta=\gamma_1 \left(
\begin{array}{cc}
f_1 & 0 \\
0 & f_2
\end{array}
\right) \gamma_2$. We put ${}^t(\omega',\theta')=\gamma_2{}^t(\omega,\theta)$, 
${}^t(\xi',\eta')=\gamma_1^{-1}{}^t(\xi,\eta)$. 
Then $\xi' \equiv f_1\omega' \pmod{\mathbb{Z}}$, $\eta'\equiv f_2\theta' \pmod{\mathbb{Z}}$. 
Translating the appropriate integral vector to ${}^t(\omega',\theta')$, 
we obtain a $\mathbb{Z}$-basis $\{1,\omega'',\theta''\}$ of $R(F)$ such that 
$\omega''\theta'' \in \mathbb{Z}$. We put $\xi''=f_1\omega''$, $\eta''=f_2\theta''$, then 
$\{1,\xi'',\eta''\}$ is a $\mathbb{Z}$-basis of $R(F')$. 
The ring structure of $R(F)$ is now given by
\[ (\omega'')^2=-\tilde{a}\tilde{c}+\tilde{b}\omega''-\tilde{a}\theta'',\quad
(\theta'')^2=-\tilde{b}\tilde{d}+\tilde{d}\omega''-\tilde{c}\theta'', \quad
\omega''\theta''=-\tilde{a}\tilde{d},
\]
where $\tilde{a}, \tilde{b}, \tilde{c}, \tilde{d}$ are the coefficients of $\gamma_2 \cdot F$. 
This implies that the ring structure of $R(F')$ is given by
\begin{align*}
(\xi'')^2 &=-f_1^2\tilde{a}\tilde{c}+f_1\tilde{b}\xi''-(f_1^2/f_2)\tilde{a}\eta'', \\
(\eta'')^2 &=-f_2^2\tilde{b}\tilde{d}+(f_2^2/f_1)\tilde{d}\xi''-f_2\tilde{c}\eta'' \\
\xi''\eta'' &=-f_1f_2\tilde{a}\tilde{d}.
\end{align*}
This proves $(\gamma_1^{-1} \cdot F')(u)=(f_1f_2)^{-1} (\gamma_2F)(f_1u_1,f_2u_2)$ 
and $F'=\delta \cdot F$. 
\qed

Let $k$ be an \'{e}tale cubic algebra and $\mathcal{O}$ be a cubic rings contained in $\mathcal{O}_k$. 
We assume that the index $f=(\mathcal{O}_k:\mathcal{O})$ is square free. 
By a theorem of elementary divisors, there exists a basis $\{1,\omega,\theta\}$ 
of $\mathcal{O}_k$ such that $\{1,f\omega,\theta\}$ is a basis of $\mathcal{O}$. 
Translating $\omega$, $\theta$ by the appropriate elements of $\mathbb{Z}$, 
we may assume that $\omega\theta \in \mathbb{Z}$. 
We call such basis \textit{normalized}. 
We take binary cubic form $F_k(u)=au_1^3+bu_1^2u_2+cu_1u_2^2+du_2^3$ 
with $a,b,c,d \in \mathbb{Z}$ such that the multiplication of $\mathcal{O}_k$ 
is given by \eqref{eq:CUBIC-RING-STR}. 
Put $\xi=f\omega$, $\eta=\theta$. 
Then the multiplication of $\mathcal{O}$ is given by
\begin{align*}
\xi^2 &=-acf^2+bf\xi-af^2\eta, \\
\eta^2 &=-bd+(d/f)\xi-c\eta, \\
\xi\eta &=-adf.
\end{align*}
Since $\mathcal{O}$ is a cubic ring with $\mathbb{Z}$-basis $\{1,\xi,\eta\}$, 
we have $f|d$. The binary cubic form
\[ F(u)=f^{-1}F_k(fu_1,u_2)=af^2u_1^3+bfu_1^2u_2+cu_1u_2^2+f^{-1}du_2^3 \]
corresponds to the cubic ring $\mathcal{O}$. 
Let $\mathfrak{f}$ be the $\mathbb{Z}$-submodule of $\mathcal{O}$ with basis 
$\{f,\xi,\eta\}$ and write $\mathfrak{f}=[f,\xi,\eta]$. 
We have
\begin{align*}
\omega\mathfrak{f}
 &=[f\omega,f\omega^2,\omega\theta]
 =[\xi,-acf+b\xi-af\eta,-af(d/f)] \subset \mathfrak{f},\\
\theta\mathfrak{f}
&=[f\theta,f\omega\theta,\theta^2]
=[f\eta,-adf, -b(d/f)f+(d/f)\xi-c\eta]\subset \mathfrak{f}.
\end{align*}
Hence $\mathfrak{f}$ is an $\mathcal{O}_k$-ideal contained in $\mathcal{O}$.
We put $\mathfrak{f}'=[f,\omega,\theta+c]$.
Then we have
\begin{align*}
\omega\mathfrak{f}'
 &=[f\omega,b\omega-a(\theta+c),-a(d/f)f+c\omega]
 \subset \mathfrak{f}',\\
\theta\mathfrak{f}'
 &=[f\theta,-a(d/f)f,-b(d/f)f+d\omega]\subset f\mathcal{O}_k \subset \mathfrak{f}'.
\end{align*}
Hence $\mathfrak{f}'$ is an $\mathcal{O}_k$-ideal. Further we have 
$\mathfrak{f}\mathfrak{f}'\subset f\mathcal{O}_k$. 
Since both of $\mathfrak{f}\mathfrak{f}'$ and $f\mathcal{O}_k$ have ideal norm $f^3$, 
we have $\mathfrak{f}\mathfrak{f}'=f\mathcal{O}_k$.
Let $\mathfrak{c}$ be an $\mathcal{O}_k$-ideal contained in $\mathcal{O}$ 
with $\mathfrak{c}\supset \mathfrak{f}$.  
We take any $\alpha \in \mathfrak{c}$ and write $\alpha=x+y\xi+z\eta$ with 
$x,y,z \in \mathbb{Z}$. Then $x=\alpha-y\xi-z\eta\in \mathfrak{c}\cap \mathbb{Z}$. 
Since $f\mathbb{Z}=\mathfrak{f}\cap \mathbb{Z}\subset \mathfrak{c}\cap \mathbb{Z}$, 
we have $\mathfrak{c}\cap \mathbb{Z}=f'\mathbb{Z}$ with $f'|f$. 
Then $f' \in \mathfrak{c}$ implies $f'\omega \in \mathfrak{c}\subset \mathcal{O}$, 
hence $f|f'$, $f'=f$. Thus $x \in f\mathbb{Z}$ and $\alpha \in \mathfrak{f}$. 
So $\mathfrak{c}\subset \mathfrak{f}$, $\mathfrak{c}=\mathfrak{f}$. 
Hence we have proved the following proposition. 
\begin{prop}\label{prop:CUNDUCTOR3}
Let $k$ be an \'{e}tale cubic algebra and let $\mathcal{O}$ be a cubic ring contained in $k$. 
Assume that the index $f=(\mathcal{O}_k:\mathcal{O})$ is square free. 
Then there exists a normalized basis $\{1,\omega,\theta\}$ of $\mathcal{O}_k$ such that 
$\mathcal{O}=[1,f\omega,\theta]$ and $\mathfrak{f}=[f,f\omega,\theta]$ is 
the conductor of $\mathcal{O}$. In particular, we have $\mathcal{O}=\mathbb{Z}+\mathfrak{f}$ and $N(\mathfrak{f})=f^2$. 
\end{prop}

\begin{prop}\label{prop:INVERTIBLE-NORMALIZED-BASIS}
Let the notations and the assumptions be as in Proposition \ref{prop:CUNDUCTOR3}. 
Then we can chose $\omega$ and $\theta$ so that $\omega, \theta \in k^\times$.
\end{prop}
\Proof
Let $F_k(u)$ and $F(u)$ be as before. 
If $k$ is a cubic field, then $F_k(u)$ is irreducible, hence $ad \neq 0$.
Next we assume that $k=\mathbb{Q}\oplus k_1$ or $k=\mathbb{Q}^3$, where 
$k_1$ is a quadratic field. If $a \neq 0$ and $d=0$, 
we take an integer $m \neq 0$ such that $af^2m^2+bfm+c \neq 0$. 
We put $F'_k(u)=F_k(u_1+fmu_2,u_2)$ and $F'(u)=F(u_1+mu_2,u_2)$. 
Then we have
\begin{align*}
F'(u)
 &=af^2u_1^3
 +(3afm+b)fu_1^2u_2+(3af^2m^2+2bfm+c)u_1u_2^2 \\
 &\quad  +m(af^2m^2+bfm+c)u_2^3.
\end{align*}
Hence the coefficients of $u_1^3$ and $u_2^3$ in $F'$ are not zero.
If $a=0$ and $d\neq 0$, we take an integer $m \neq 0$ such that $b+cm+dm^2 \neq 0$. 
We put $F'_k(u)=F_k(u_1,mu_1+u_2)$ and $F'(u)=F(u_1,mfu_1+u_2)$. 
Then the coefficients of $u_1^3$ and $u_2^3$ in $F'$ are not zero.
Finally, if $a=d=0$, then $k=\mathbb{Q}^3$. Since $\Disc(F_k)=\Disc(\mathbb{Z}^3)=1$, 
we have $b^2c^2=1$. Hence $F_k(u)=\pm u_1u_2(u_1\pm u_2)$.
We put $\alpha=\left(
\begin{array}{cc}
1+4f & 2 \\
2f & 1
\end{array}
\right)$, $\beta={}^t\! \alpha$ and put $F'_k(u)=F_k(u\alpha)$, $F'(u)=F(u\beta)$.
Then we have
\begin{align*}
F'(1,0)&=F(1+4f,2f)=\pm 2f(1+4f)(1+4f\pm 2f)\neq 0, \\
F'(0,1)&=F(2,1)=\pm 2 (2\pm 1) \neq 0.
\end{align*}
Thus the coefficients of $u_1^3$ and $u_2^3$ in $F'$ are nonzero. 
By the argument above, we may always assume that $ad \neq 0$ and $a>0$. 
\qed

Let $k$ be a number field of degree $n$ over $\mathbb{Q}$. 
We say that a prime number $p$ is of \textit{type} $f_1^{e_1}\cdots f_g^{e_g}$ in $k$ 
if the prime ideal decomposition of $p\mathcal{O}_k$ is of the form $\mathfrak{p}_1^{e_1}\cdots \mathfrak{p}_g^{e_g}$ with $N(\mathfrak{p}_i)=p^{f_i}$. 
Let $k=k_1 \oplus \cdots \oplus k_r$ be an \'{e}tale algebra of degree $n$ over $\mathbb{Q}$, where 
$k_i$'s are number fields. We say that $p$ is of type $f_{11}^{e_{11}} \cdots f_{rg_r}^{e_{rg_r}}$ 
if $p$ is of type $f_{i1}^{e_{i1}}\cdots f_{ig_i}^{e_{ig_i}}$ for each $i$.
\begin{lem}\label{lem:P-PART-CUBICRING-CONDUCTOR}
Let $k$ be an \'{e}tale cubic algebra and $p$ be a prime number. 
There exist orders $R$ of $k$ such that $(\mathcal{O}_k:R)=p$ 
if and only if $p$ is not of type $3$ in $k$. 
If this is the case, for each ideal $\mathfrak{f}$ of $\mathcal{O}_k$ such that 
$p\mathcal{O}_k \subset \mathfrak{f}$ and $N(\mathfrak{f})=p^2$, 
there exists exactly one order $R$  with index $p$ in $\mathcal{O}_k$ whose conductor is $\mathfrak{f}$, 
namely $R=\mathbb{Z}+\mathfrak{f}$. 
\end{lem}
\Proof
If there exists an order $R$ of $k$ such that $(\mathcal{O}_k:R)=p$, 
then the conductor $\mathfrak{f}$ of $R$ has norm $p^2$ and $R$ is uniquely determined by the conductor $\mathfrak{f}$ 
by Proposition \ref{prop:CUNDUCTOR3}. Hence $p$ is not of type $3$ in $k$. 
Conversely, if $p$ is not of type $3$ in $k$, then there exists an ideal $\mathfrak{f}$ such that 
$p\mathcal{O}_k \subset \mathfrak{f}$ and $N(\mathfrak{f})=p^2$. 
We put $R=\mathbb{Z}+\mathfrak{f}$. 
Then it is clear that $R$ is an order of $k$. 
Since $p\mathbb{Z}\subset \mathbb{Z} \cap \mathfrak{f}\subsetneq \mathbb{Z}$, 
we have $\mathbb{Z}\cap \mathfrak{f}=p\mathbb{Z}$, hence 
$R/\mathfrak{f} \cong \mathbb{Z}/p\mathbb{Z}$. 
So the index $(\mathcal{O}_k:R)=(\mathcal{O}_k:\mathfrak{f})/(R:\mathfrak{f})=p^2/p=p$. 
\qed

\begin{prop}\label{prop:SQUARE-FREE-INDEX-ORDER}
Let $k$ be an \'{e}tale cubic algebra and $f$ be a square free positive integer 
such that each prime number dividing $f$ is not of type $3$. 
Then for each ideal $\mathfrak{f}$ of $\mathcal{O}_k$ such that 
$f\mathcal{O}_k \subset \mathfrak{f}$ and $N(\mathfrak{f})=f^2$, 
there exists exactly one order $R$ with index $f$ in $\mathcal{O}_k$ whose conductor is $\mathfrak{f}$, 
namely $R=\mathbb{Z}+\mathfrak{f}$. 
\end{prop}
\Proof
For each prime number $p$ dividing $f$, we denote by $\mathfrak{f}_p$ the $p$-part of $\mathfrak{f}$. 
Then we have $p\mathcal{O}_k \subset \mathfrak{f}_p$ and $N(\mathfrak{f}_p)=p^2$. 
By Lemma \ref{lem:P-PART-CUBICRING-CONDUCTOR}, 
$R_p=\mathbb{Z}+\mathfrak{f}_p$ is the unique order 
with index $p$ in $\mathcal{O}_k$ whose conductor is $\mathfrak{f}_p$. 
For any positive integer $g$ dividing $f$, 
we write $f=gh$, $\mathfrak{g}=\prod_{p|g} \mathfrak{f}_p$ and $\mathfrak{h}=\prod_{p|h} \mathfrak{f}_p$. 
If $S$ is an order with index $g$ whose conductor is $\mathfrak{g}$ 
and $T$ is an order with index $h$ whose conductor is $\mathfrak{h}$, 
then  $S+T=\mathcal{O}_k$ and $S/(S\cap T) \cong \mathcal{O}_k/T$ as $\mathbb{Z}$-modules 
since $g$ and $h$ are relatively prime to each other. 
Hence $(\mathcal{O}_k:S\cap T)=(\mathcal{O}_k:S)(S:S\cap T)=gh=f$. 
So $S\cap T$ is an order with index $f$. 
Since $\mathfrak{f}=\mathfrak{g}\mathfrak{h} \subset S \cap T$, 
the conductor $\mathfrak{f}'$ of $S\cap T$ contains $\mathfrak{f}$. 
By Proposition \ref{prop:CUNDUCTOR3}, $N(\mathfrak{f}')=f^2=N(\mathfrak{f})$, 
hence $\mathfrak{f}'=\mathfrak{f}$. 
So the statement of the proposition now follows by induction on the number of prime divisors of $f$. 
\qed

\section{Ideal class groups of cubic rings}
In this section, we summarise some basic facts on ideal class groups of rings of rank $n$. 
Let $k$ be an \'{e}tale algebra of degree $n$ over $\mathbb{Q}$ 
and let $\mathcal{O}$ be an order of $k$. 
We follow the argument of Sand \cite{SA} in which $k$ is assumed to be a number field. 
Let $\mathfrak{f}$ be the conductor of $\mathcal{O}$ and put $f=(\mathcal{O}_k:\mathcal{O})$.
We denote by $I_\mathcal{O}$ the group of invertible fractional ideals of $\mathcal{O}$, 
$P_\mathcal{O}$ the subgroup of $I_\mathcal{O}$ consisting 
of principal invertible fractional ideals of $\mathcal{O}$, respectively.
The ideal class group of $\mathcal{O}$ is $\Cl_\mathcal{O}=I_\mathcal{O}/P_\mathcal{O}$. 
For the maximal order $\mathcal{O}_k$, we write $I_k=I_{\mathcal{O}_k}$, $P_k=P_{\mathcal{O}_k}$ 
and $\Cl_k=\Cl_{\mathcal{O}_k}$. 
$\Cl_\mathcal{O}$ and $\Cl_k$ are finite abelian groups. 
We denote by $h_\mathcal{O}$ (resp. $h_k$) the order of $\Cl_\mathcal{O}$ (resp. $\Cl_k$) 
and call the class numbers of $\mathcal{O}$ (resp. $k$). 
The next lemma is a restatement of Theorem 3.1 and Corollary 3.2 in \cite{SA}. 
\begin{lem}[Dedekind]\label{lem:EXTEND-RESTRICT-IDEALS}
The extension map $\mathfrak{a}\mapsto \mathfrak{a}\mathcal{O}_k$ defines
a multiplicative bijection between the monoid of integral $\mathcal{O}$-ideals which are 
relatively prime to $\mathfrak{f}$ and the monoid of integral $\mathcal{O}_k$-ideals which are 
relatively prime to $\mathfrak{f}$. 
The inverse of the extension map is the contraction map $\mathfrak{A}\mapsto \mathfrak{A}\cap \mathcal{O}$. 
In particular, for any element $\gamma \in \mathcal{O}\cap k^\times$ 
with $\gamma\mathcal{O}+\mathfrak{f}=\mathcal{O}$, 
we have $(\gamma \mathcal{O}_k)\cap \mathcal{O}=\gamma \mathcal{O}$. 
\end{lem}

The following corollary is a restatement of \cite[Corollary 3.3]{SA}.
\begin{cor}\label{cor:PRIMETOCOND-INVERTIBLE}
If $\mathfrak{a}$ is an integral ideal of $\mathcal{O}$ which is relatively 
prime to $\mathfrak{f}$ and $\mathfrak{a}\mathcal{O}_k$ is an invertible $\mathcal{O}_k$-ideal, 
then $\mathfrak{a}$ is an invertible $\mathcal{O}$-ideal. 
\end{cor}

The following lemma is \cite[Proposition 2.6]{SA}. 
\begin{lem}\label{lem:TORSION-IDEAL1}
The mapping of $I_\mathcal{O}$ to $I_k$ defined by $\mathfrak{a}\mapsto \mathfrak{a}\mathcal{O}_k$ 
is a group homomorphism. The kernel is the torsion subgroup $\Tor(I_\mathcal{O})$ of $I_\mathcal{O}$ 
and is a finite group.
\end{lem}

The following proposition is a restatement of \cite[Proposition 3.6]{SA}.
\begin{prop}\label{prop:TORSION-IDEALS3}
If $\mathfrak{a} \in \Tor(I_\mathcal{O})$, then there exists an element $\alpha \in \mathcal{O}_k$
such that $\mathfrak{a}=\alpha\mathcal{O}+\mathfrak{f}$ and 
$\alpha\mathcal{O}_k+\mathfrak{f}=\mathcal{O}_k$. Conversely, 
if $\alpha \in \mathcal{O}_k$ satisfies $\alpha\mathcal{O}_k+\mathfrak{f}=\mathcal{O}_k$, 
then $\alpha\mathcal{O}+\mathfrak{f} \in \Tor(I_\mathcal{O})$.
\end{prop}

\begin{rem}\label{rem:ALPHA-MOD-F}
In Proposition \ref{prop:TORSION-IDEALS3}, $\alpha$ can be replaced by any $\alpha'$ 
which is congruent to $\alpha$ modulo $\mathfrak{f}$. 
So we may assume that $\alpha \in \mathcal{O}_k \cap k^\times$. 
Moreover we may assume that $\alpha$ is totally positive.
\end{rem}

The following lemma is \cite[Theorem 3.7]{SA}. 
\begin{lem}\label{lem:ORDER-TORSION-IO}
The order of the torsion subgroup $\Tor(I_\mathcal{O})$ is given by
\[ |\Tor(I_\mathcal{O})|=\frac{\varphi_k(\mathfrak{f})}{\varphi_\mathcal{O}(\mathfrak{f})}, \quad
\varphi_k(\mathfrak{f})=|(\mathcal{O}_k/\mathfrak{f})^\times|, \; 
\varphi_\mathcal{O}(\mathfrak{f})=|(\mathcal{O}/\mathfrak{f})^\times|. 
\]
\end{lem}

We denote by $I_\mathcal{O}(\mathfrak{f})$ the subgroup of $I_\mathcal{O}$ 
consisting of all invertible fractional ideals of $\mathcal{O}$ which are 
relatively prime to $\mathfrak{f}$. 
Hence an element of $I_\mathcal{O}(\mathfrak{f})$ is of the form $\mathfrak{a}\mathfrak{b}^{-1}$ 
where $\mathfrak{a}$ and $\mathfrak{b}$ are integral invertible $\mathcal{O}$-ideals 
relatively prime to $\mathfrak{f}$.  
Put $P_\mathcal{O}(\mathfrak{f})=P_\mathcal{O} \cap I_\mathcal{O}(\mathfrak{f})$. 
\begin{lem}\label{lem:POF}
Any element of $P_\mathcal{O}(\mathfrak{f})$ can be written of the form 
$\gamma\mathcal{O}$, $\gamma=\alpha/\beta$, where $\alpha,\beta \in \mathcal{O}$ 
are relatively prime to $\mathfrak{f}$. 
\end{lem}
\Proof
Let $\gamma\mathcal{O}$ be an element of $P_\mathcal{O}(\mathfrak{f})$.
We write $\gamma\mathcal{O}=\mathfrak{a}\mathfrak{b}^{-1}$, 
where $\mathfrak{a}$ and $\mathfrak{b}$ are integral invertible $\mathcal{O}$-ideals relatively prime to $\mathfrak{f}$. 
It is clear that $\mathfrak{a}\mathcal{O}_k$ is an invertible $\mathcal{O}_k$-ideal, 
so we have $(\mathfrak{a}\mathcal{O}_k)^{h_k}=\xi \mathcal{O}_k$ for some $\xi \in \mathcal{O}_k \cap k^\times$ which is relatively prime to $\mathfrak{f}$. 
Then $\xi^m \in 1+\mathfrak{f}\subset \mathcal{O}$ for some positive integer $m$. 
By Lemma \ref{lem:EXTEND-RESTRICT-IDEALS}, we have $\mathfrak{a}^{h_km}=\xi^m \mathcal{O}$, 
hence $\mathfrak{b}\mathfrak{a}^{h_km-1}=\xi^m\gamma^{-1}\mathcal{O}$ 
is an integral principal invertible $\mathcal{O}$-ideal relatively prime to $\mathfrak{f}$. 
Put $\beta=\xi^m\gamma^{-1}$, $\alpha=\xi^m$. Then $\gamma=\alpha/\beta$ and 
$\alpha, \beta \in \mathcal{O}$ are relatively prime to $\mathfrak{f}$.
\qed

%
The following lemma is \cite[Lemma 4.2]{SA}. 
\begin{lem}\label{lem:TORIO-IOFPO}
$\Tor(I_\mathcal{O})\subset I_\mathcal{O}(\mathfrak{f})P_\mathcal{O}$.
\end{lem}
The following proposition is \cite[Proposition 4.3]{SA}. 
\begin{prop}\label{prop:CLASSO}
$I_\mathcal{O}(\mathfrak{f})/P_\mathcal{O}(\mathfrak{f}) \cong \Cl_\mathcal{O}$.
\end{prop}

We have the following exact sequence. 
\begin{equation}\label{eq:ORDERCG}
1\longrightarrow \Tor(I_\mathcal{O})P_\mathcal{O}/P_\mathcal{O} 
\stackrel{\iota}{\longrightarrow} \Cl_\mathcal{O}
\stackrel{\epsilon}{\longrightarrow} \Cl_k \longrightarrow 1, 
\end{equation}
where $\iota$ is induced by the inclusion  in Lemma \ref{lem:TORIO-IOFPO} 
and $\epsilon$ is induce by extension of fractional ideals. \eqref{eq:ORDERCG} was given in the 
proof of \cite[Theorem 4.4]{SA}. 
The following corollary is \cite[Theorem 4.4]{SA}. 
\begin{cor}\label{cor:DEDEKIND}
The class number $h_\mathcal{O}$ is given by
\[h_\mathcal{O}=h_k
\frac{((\mathcal{O}_k/\mathfrak{f})^\times:(\mathcal{O}/\mathfrak{f})^\times)}
{(U(\mathcal{O}_k):U(\mathcal{O}))},
\]
where $U(\mathcal{O}_k)$  \rm{(}resp. $U(\mathcal{O})$\rm{)} denotes the group of units in $\mathcal{O}_k$ \rm{(}resp. $\mathcal{O}$\rm{)}. 
\end{cor}

We denote by $P_{\mathcal{O},+}$ the subgroup of $P_\mathcal{O}$ consisting of all 
principal invertible fractional ideals of $\mathcal{O}$ 
generated by totally positive elements in $k^\times$. 
We put $\Cl_{\mathcal{O},+}=I_\mathcal{O}/P_{\mathcal{O},+}$. Further we put
$P_{\mathcal{O},+}(\mathfrak{f})=P_{\mathcal{O},+}\cap I_\mathcal{O}(\mathfrak{f})$ 
and write $\Cl_{k,+}=\Cl_{\mathcal{O}_k,+}$.
\begin{lem}\label{lem:IO-PLUS}
$I_\mathcal{O}=I_\mathcal{O}(\mathfrak{f})P_{\mathcal{O},+}\Tor(I_\mathcal{O})$.
\end{lem}
\Proof
Let $\mathfrak{a} \in I_\mathcal{O}$ and 
put $\tilde{\mathfrak{a}}=\mathfrak{a}\mathcal{O}_k$. 
Take an integral invertible $\mathcal{O}_k$-ideal $\tilde{\mathfrak{b}}$ which 
is relatively prime to $\mathfrak{f}$ and belongs to the same ideal class of 
$\tilde{\mathfrak{a}}$ in $\Cl_{k,+}$. Then $\tilde{\mathfrak{a}}=\gamma\tilde{\mathfrak{b}}$ 
for some totally positive element $\gamma$ in $k^\times$. 
Put $\mathfrak{b}=\tilde{\mathfrak{b}} \cap \mathcal{O}$. 
Then $\mathfrak{b}$ is relatively prime to $\mathfrak{f}$ and $\mathfrak{b}\mathcal{O}_k=\tilde{\mathfrak{b}}$ 
by Lemma \ref{lem:EXTEND-RESTRICT-IDEALS}. 
Hence $\mathfrak{b}$ is an invertible $\mathcal{O}$-ideal 
by Corollary \ref{cor:PRIMETOCOND-INVERTIBLE}. 
This proves $\mathfrak{b} \in I_\mathcal{O}(\mathfrak{f})$.
If we put $\mathfrak{c}=\gamma^{-1}\mathfrak{a}\mathfrak{b}^{-1}$, 
then $\mathfrak{c} \in I_\mathcal{O}$ and $\mathfrak{c}\mathcal{O}_k=\mathcal{O}_k$. 
Hence $\mathfrak{c} \in \Tor(I_\mathcal{O})$ by Lemma \ref{lem:TORSION-IDEAL1}. 
So we have $\mathfrak{a}=\mathfrak{b} \gamma\mathcal{O} \mathfrak{c} \in I_\mathcal{O}(\mathfrak{f})P_{\mathcal{O},+}\Tor(I_\mathcal{O})$. 
This proves one inclusion, and the reverse is clear.
\qed

\begin{lem}\label{lem:TORIO-IOFPOPLUS}
$\Tor(I_\mathcal{O})\subset I_\mathcal{O}(\mathfrak{f})P_{\mathcal{O},+}$.
\end{lem}
\Proof
Take an element $\mathfrak{a} \in \Tor(I_\mathcal{O})$. By Proposition \ref{prop:TORSION-IDEALS3}, 
$\mathfrak{a}=\alpha\mathcal{O}+\mathfrak{f}$ and $\mathfrak{a}^{-1}=\beta\mathcal{O}+\mathfrak{f}$ 
for some $\alpha,\beta \in \mathcal{O}_k$ such that $\alpha\beta \in \mathcal{O}$ is 
relatively prime to $\mathfrak{f}$. By Remark \ref{rem:ALPHA-MOD-F}, we may assume that 
$\alpha$ is totally positive. We write 
$\mathfrak{a}=\alpha\mathcal{O}\left(\alpha\beta\mathcal{O}+\alpha\mathfrak{f} \right)^{-1}$. 
It is obvious that $\alpha\beta\mathcal{O}+\alpha\mathfrak{f}$ is an integral $\mathcal{O}$-ideal
which is relatively prime to $\mathfrak{f}$. 
\qed

By Lemma \ref{lem:IO-PLUS} and \ref{lem:TORIO-IOFPOPLUS}, we have 
$I_\mathcal{O}=I_\mathcal{O}(\mathfrak{f})P_{\mathcal{O},+}$. Hence
\begin{equation}\label{eq:CLO-PLUS}
\Cl_{\mathcal{O},+}
=I_\mathcal{O}(\mathfrak{f})P_{\mathcal{O},+}/P_{\mathcal{O},+}
\cong I_\mathcal{O}(\mathfrak{f})/P_{\mathcal{O},+}(\mathfrak{f}).
\end{equation}
We also have the following exact sequence 
\begin{equation}\label{eq:ORDERCG-PLUS}
 1 \longrightarrow \Tor(I_\mathcal{O}) P_{\mathcal{O},+}/P_{\mathcal{O},+}
     \longrightarrow \Cl_{\mathcal{O},+} \longrightarrow \Cl_{k,+} \longrightarrow 1.
\end{equation}
We denote by $U_{+}(\mathcal{O}_k)$ (resp. $U_{+}(\mathcal{O})$) 
the group of totally positive units in $\mathcal{O}_k$ (resp. $\mathcal{O}$).
\begin{cor}\label{cor:CARDORDER-PLUS}
\[ |\Cl_{\mathcal{O},+}|=|\Cl_{k,+}|\dfrac{|(\mathcal{O}_k/\mathfrak{f})^\times|/|(\mathcal{O}/\mathfrak{f})^\times|}
{(U_{+}(\mathcal{O}_k):U_{+}(\mathcal{O}))}.
\]
\end{cor}
We used a homomorphism from $\Cl_\mathcal{O}$ to $\Cl_k$ in the proof of 
Corollary \ref{cor:DEDEKIND}.  We now use a homomorphism from a ray class group of $k$ to $\Cl_\mathcal{O}$. 
We denote by $P_{k,1}(\mathfrak{f})$ the subgroup of $P_k$ consisting of all 
invertible principal ideals $\alpha\mathcal{O}_k$ with $\alpha \equiv 1 \mmod{\mathfrak{f}}$.
Here $\alpha \equiv 1 \mmod{\mathfrak{f}}$ means $\alpha=\beta/\gamma$ for some 
$\beta, \gamma \in \mathcal{O}_k\cap k^\times$ which are relatively prime to $\mathfrak{f}$ and 
$\beta\equiv \gamma \pmod{\mathfrak{f}}$.  
We denote by $P_{k,+}(\mathfrak{f})$ the subgroup of $P_{k,1}(\mathfrak{f})$ 
consisting of all invertible principal ideals $\alpha\mathcal{O}_k$ with 
totally positive $\alpha \equiv 1 \mmod{\mathfrak{f}}$.
Then the quotient groups $\Cl_k(\mathfrak{f})=I_k(\mathfrak{f})/P_{k,1}(\mathfrak{f})$ 
and $\Cl_{k,+}(\mathfrak{f})=I_k(\mathfrak{f})/P_{k,+}(\mathfrak{f})$ 
are called \textit{the ray class group} of $k$ modulo $\mathfrak{f}$  
and that in the narrow sense, respectively.

Let $\mathfrak{C}$ be an element of $I_k(\mathfrak{f})$. 
We write $\mathfrak{C}=\mathfrak{A}\mathfrak{B}^{-1}$ where $\mathfrak{A}$ and $\mathfrak{B}$ 
are integral invertible $\mathcal{O}_k$ ideals  relatively prime to $\mathfrak{f}$. 
Then $\mathfrak{a}=\mathfrak{A}\cap \mathcal{O}$ and $\mathfrak{b}=\mathfrak{B}\cap \mathcal{O}$ 
are elements of $I_\mathcal{O}(\mathfrak{f})$ by Lemma \ref{lem:EXTEND-RESTRICT-IDEALS} and Corollary \ref{cor:PRIMETOCOND-INVERTIBLE}. 
We correspond $\mathfrak{C}$ to $\mathfrak{c}=\mathfrak{a}\mathfrak{b}^{-1} \in I_\mathcal{O}(\mathfrak{f})$. 
It is easy to see that the correspondence is well defined. 
For any $\gamma\mathcal{O}_k \in P_{k,1}(\mathfrak{f})$ (resp. $\gamma\mathcal{O}_k \in P_{k,+}(\mathfrak{f})$), 
we write $\gamma=\alpha/\beta$ where $\alpha, \beta \in \mathcal{O}_k$ are relatively prime to $\mathfrak{f}$ 
and $\alpha\equiv \beta \pmod{\mathfrak{f}}$ (resp. $\alpha\equiv \beta \pmod{\mathfrak{f}}$ and 
$\alpha/\beta$ is totally positive). Hence there exists an element $\beta' \in \mathcal{O}_k\cap k^\times$ 
such that $\alpha\beta'\equiv \beta\beta' \equiv 1 \pmod{\mathfrak{f}}$. 
Then $\alpha\beta', \beta\beta' \in \mathcal{O}$ are relatively prime to $\mathfrak{f}$ 
and $\alpha\beta'\mathcal{O}_k \cap \mathcal{O}=\alpha\beta'\mathcal{O}$, 
$\beta\beta'\mathcal{O}_k \cap \mathcal{O}=\beta\beta' \mathcal{O}$. 
By the correspondence above, 
$\gamma\mathcal{O}_k$ corresponds to $\gamma \mathcal{O} \in P_\mathcal{O}(\mathfrak{f})$ 
(resp. $\gamma\mathcal{O} \in P_{\mathcal{O},+}(\mathfrak{f})$).
By Proposition \ref{prop:CLASSO} and \eqref{eq:CLO-PLUS}, 
the correspondence $\mathfrak{C}\mapsto \mathfrak{c}$ induces 
a homomorphism $\rho:\Cl_k(\mathfrak{f})\rightarrow \Cl_\mathcal{O}$ 
(resp. $\rho_+:\Cl_{k,+}(\mathfrak{f})\rightarrow \Cl_{\mathcal{O},+}$). 
By Lemma \ref{lem:EXTEND-RESTRICT-IDEALS}, Proposition \ref{prop:CLASSO} and \eqref{eq:CLO-PLUS},
$\rho$ and $\rho_+$ are surjective.
Let $\mathfrak{C} \in I_k(\mathfrak{f})$ be a representative of an element in $\ker \rho$ 
(resp. $\ker \rho_+$). We write $\mathfrak{C}=\mathfrak{A}\mathfrak{B}^{-1}$ as before.
Put $\mathfrak{a}=\mathfrak{A} \cap \mathcal{O}$, $\mathfrak{b}=\mathfrak{B}\cap \mathcal{O}$. 
Then $\mathfrak{a}, \mathfrak{b} \in I_\mathcal{O}(\mathfrak{f})$ are integral $\mathcal{O}$-ideals 
such that $\mathfrak{a}\mathfrak{b}^{-1}=\gamma \mathcal{O} \in P_\mathcal{O}(\mathfrak{f})$ 
(resp. $P_{\mathcal{O},+}(\mathfrak{f})$). 
Since $\mathfrak{a}\mathcal{O}_k=\mathfrak{A}$ and $\mathfrak{b}\mathcal{O}_k=\mathfrak{B}$, 
we have $\mathfrak{A}=\mathfrak{a}\mathcal{O}_k=\gamma \mathfrak{b}\mathcal{O}_k=\gamma \mathfrak{B}$,  
hence $\mathfrak{C}=\gamma\mathcal{O}_k \in P_k\cap I_k(\mathfrak{f})$ 
(resp. $P_{k,+}\cap I_k(\mathfrak{f})$). 
We can write $\gamma=\alpha/\beta$ where $\alpha, \beta \in \mathcal{O}\cap k^\times$ are relatively prime to $\mathfrak{f}$ 
by Lemma \ref{lem:POF}. 

Conversely, assume that $\gamma=\alpha/\beta$ (resp. $\gamma=\alpha/\beta$ is totally positive) 
and $\alpha, \beta \in \mathcal{O}\cap k^\times$ are relatively prime to $\mathfrak{f}$.
We write $\gamma\mathcal{O}_k=\mathfrak{A}\mathfrak{B}^{-1}$ 
where $\mathfrak{A}$ and $\mathfrak{B}$ are integral invertible $\mathcal{O}_k$ ideals  relatively prime to $\mathfrak{f}$ 
and $\mathfrak{A}+\mathfrak{B}=\mathcal{O}_k$. 
Then $\alpha\mathfrak{B}=\beta\mathfrak{A}$, hence 
$(\alpha\mathfrak{B})\cap \mathcal{O}=(\beta\mathfrak{A})\cap \mathcal{O}$. 
Put $\mathfrak{a}=\mathfrak{A} \cap \mathcal{O}$ and $\mathfrak{b}=\mathfrak{B}\cap \mathcal{O}$. 
By Lemma \ref{lem:EXTEND-RESTRICT-IDEALS}, we have $(\alpha\mathcal{O})\mathfrak{b}=(\beta\mathcal{O})\mathfrak{a}$, hence 
$\mathfrak{a}\mathfrak{b}^{-1}=(\alpha/\beta)\mathcal{O}=\gamma\mathcal{O}$. 
Hence $\gamma\mathcal{O}_k$ is a representative of an element in $\ker \rho$ 
(resp. $\ker \rho_+$). 
We put
\begin{equation}\label{eq:PO-POTILDE}
\begin{split}
\tilde{P}_\mathcal{O}
 &=\left\{\alpha\beta^{-1}\mathcal{O}_k\,|\, \alpha,\beta \in \mathcal{O}\cap k^\times,\; 
 \alpha\mathcal{O}+\mathfrak{f}=\beta\mathcal{O}+\mathfrak{f}=\mathcal{O}\right\}, \\
\tilde{P}_{\mathcal{O},+}
 &=\left\{\gamma\mathcal{O}_k\in \tilde{P}_\mathcal{O}\,|\, \text{$\gamma$ is totally positive}\right\}.
\end{split}
\end{equation}
Then we have the following exact sequences of abelian groups. 
\begin{align*}
&1\longrightarrow \tilde{P}_\mathcal{O}/P_{k,1}(\mathfrak{f})
\longrightarrow I_k(\mathfrak{f})/P_{k,1}(\mathfrak{f})
\stackrel{\rho}{\longrightarrow} \Cl_\mathcal{O} \longrightarrow 1, \\
&1\longrightarrow \tilde{P}_{\mathcal{O},+}/P_{k,+}(\mathfrak{f})
\longrightarrow I_k(\mathfrak{f})/P_{k,+}(\mathfrak{f})
\stackrel{\rho_+\;}{\longrightarrow} \Cl_{\mathcal{O},+} \longrightarrow 1.
\end{align*}
Thus we have
\[ \Cl_\mathcal{O}\cong I_k(\mathfrak{f})/\tilde{P}_\mathcal{O}, \quad 
\Cl_{\mathcal{O},+} \cong I_k(\mathfrak{f})/\tilde{P}_{\mathcal{O},+}.
\]
Hence 
\begin{equation}\label{eq:CLO-RAYCLASSGROUP}
 \Cl_\mathcal{O}/\Cl_\mathcal{O}^2 \cong I_k(\mathfrak{f})/\tilde{P}_\mathcal{O}I_k(\mathfrak{f})^2, \quad 
\Cl_{\mathcal{O},+}/\Cl_{\mathcal{O},+}^2 \cong I_k(\mathfrak{f})/\tilde{P}_{\mathcal{O},+}I_k(\mathfrak{f})^2.
\end{equation}

\section{Pairs of integral symmetric matrices of degree three}
We denote by $\hat{L}$ the set of pairs of symmetric matrices of degree three with coefficients in $\mathbb{Z}$. 
Let $(A,B) \in \hat{L}$ be a pair of integral symmetric matrices of degree three. 
The group $\Gamma_1=\SL_3(\mathbb{Z})$ acts on $\hat{L}$ by 
\[ T\cdot (A,B)=(TA\,{}^tT,TB\,{}^tT), \quad T \in \SL_3(\mathbb{Z}). \]
The group $\Gamma_2=\GL_2(\mathbb{Z})$ acts on $\hat{L}$ by
\[ g\cdot (A,B)=((\det g)^{-1}(pA-qB),(\det g)^{-1}(-rA+sB)),\]
where $g=\left(
\begin{array}{cc}
p & q \\
r & s
\end{array}
\right) \in \Gamma_2$. 
Thus the group $\Gamma=\Gamma_1 \times \Gamma_2$ acts on $\hat{L}$. 
For any pair $(A,B) \in \hat{L}$, 
we defined in \S1 a binary cubic form $\hat{F}_{(A,B)}(u)$ in two variables $u=(u_1,u_2)$, 
which equals to $\det(u_1A-u_2B)$. The discriminant of $(A,B)$ is defined by $\Disc^*(A,B)=\Disc(\hat{F}_{(A,B)})$. 
We say that $(A,B)$ is \textit{nondegenerate} if $\Disc^*(A,B) \neq 0$. 
We note that the set of nondegenerate pairs of $\hat{L}$ equals $\HATSSV_\mathbb{Z}$ in the introduction. 

We consider triples $(\mathcal{O},\mathfrak{a},\delta)$, where $\mathcal{O}$ is 
a nondegenerate cubic ring, $\mathfrak{a}$ is a fractional ideal of $\mathcal{O}$ and 
$\delta$ is an invertible element of $k=\mathcal{O}\otimes_\mathbb{Z}\mathbb{Q}$ such that 
$\mathfrak{a}^2 \subset \delta\mathcal{O}$ and $N_{k/\mathbb{Q}}(\delta)=N_{\mathcal{O}}(\mathfrak{a})^2$. 
Here $N_{\mathcal{O}}(\mathfrak{a})$ is the norm of $\mathfrak{a}$ as a fractional $\mathcal{O}$-ideal, 
i. e. $N_{\mathcal{O}}(\mathfrak{a})=(\mathcal{O}:\mathfrak{a})$ for $\mathfrak{a}\subset \mathcal{O}$.
Two such triples $(\mathcal{O},\mathfrak{a},\delta)$ and $(\mathcal{O}',\mathfrak{a}',\delta')$ 
are called \textit{equivalent} if there exists an isomorphism $\phi:\mathcal{O}\rightarrow \mathcal{O}'$ 
and $\kappa \in \mathcal{O}'\otimes_\mathbb{Z}\mathbb{Q}$ such that 
$\mathfrak{a}'=\kappa \phi(\mathfrak{a})$, $\delta'=\kappa^2\phi(\delta)$. 
M. Bhargava proved the following theorem (\cite[Theorem 4]{BH2}). 
\begin{thm}\label{thm:BHARGAVA1}
There is a canonical bijection between the set of nondegenerate
$\Gamma$-orbits on $\hat{L}$ and the set of equivalence classes of triples $(\mathcal{O},\mathfrak{a},\delta)$. 
Under this bijection, the discriminant of 
a pair of integral matrices of degree three equals the discriminant of the corresponding
cubic ring.
\end{thm}
We now explain the correspondence of Theorem \ref{thm:BHARGAVA1}. 
Let $\mathcal{O}$ be a nondegenerate cubic ring and take a normalized basis $\{1,\omega,\theta\}$ of $\mathcal{O}$. 
The multiplication of $\mathcal{O}$ is given by \eqref{eq:CUBIC-RING-STR} with $a,b,c,d \in \mathbb{Z}$. 
Let $\mathfrak{a}$ be a fractional $\mathcal{O}$-ideal and $\delta$ be an invertible element of $k=\mathcal{O} \otimes_\mathbb{Z} \mathbb{Q}$ 
such that $\mathfrak{a}^2 \subset \delta\mathcal{O}$ and $N_{k/\mathbb{Q}}(\delta)=N_{\mathcal{O}}(\mathfrak{a})^2$. 
We take a $\mathbb{Z}$-basis $\{\alpha_1,\alpha_2,\alpha_3\}$ of the ideal $\mathfrak{a}$ 
having the same orientation as $\{1,\omega,\theta\}$. 
Since $\mathfrak{a}^2\subset \delta\mathcal{O}$, 
there exist integers $a_{ij}$, $b_{ij}$ and $c_{ij}$ such that
\begin{equation}\label{eq:SYMMETRIC-COND1}
\alpha_i\alpha_j=\delta(c_{ij}+b_{ij}\omega+a_{ij}\theta). 
\end{equation}
We put $A=(a_{ij})$, $B=(b_{ij})$. 
Then we have
\[ \hat{F}_{(A,B)}(u)=au_1^3+bu_1^2u_2+cu_1u_2^2+du_2^3.\]
This was shown in the proof of \cite[Therem 4]{BH2}. 

We denote by $\Gamma_{(A,B)}$ the isotropy group in $\Gamma$ of a nondegenerate pair $(A,B) \in \hat{L}$. 
The following corollary is \cite[Corollary 5]{BH2}. 
\begin{cor}\label{cor:ISOTROPY2}
For any nondegenerate pair $(A,B) \in \hat{L}$, there exists a homomorphism 
$\Gamma_{(A,B)} \rightarrow \Aut(\mathcal{O})$ with kernel isomorphic to $U^+_2(\mathcal{O}_0)$. 
Here $(\mathcal{O},\mathfrak{a})$ is the pair corresponding to $(A,B)$ as 
in Theorem \ref{thm:BHARGAVA1}, $\mathcal{O}_0=\End_\mathcal{O}(\mathfrak{a})$
is the endomorphism ring of $\mathfrak{a}$, 
and $U_2^+(\mathcal{O}_0)$ denotes the group of units of $\mathcal{O}_0$ 
having order dividing $2$ and positive norm.
\end{cor}

\begin{rem}\label{rem:ISOTROPY3}
The following example shows that the statement of \cite[Corollary 5]{BH2} is not correct. 
So we have given a weaker statement than the original one. 
\end{rem}

\begin{ex}\label{ex:ISOTROPY4}
Put $f(x)=x^3-14x^2+11x+1 \in \mathbb{Q}[x]$ and let $\omega$ be a root of $f(x)$. 
Then the discriminant of the cubic field $k=\mathbb{Q}(\omega)$ is $163^2$, 
hence $k$ is a Galois cubic field. We put $\theta=-11+14\omega-\omega^2$. Then $\{1,\omega,\theta\}$ is a $\mathbb{Z}$-basis of the maximal order $\mathcal{O}_k$
such that 
\[ \omega^2=-11+14\omega-\theta,\quad \theta^2=14-\omega-11\theta,\quad \omega\theta=1.\]
We can take a generator $\sigma$ of the Galois group $\Gal(k/\mathbb{Q})$ such that 
$\sigma(\omega)=13-\omega+\theta$, $\sigma(\theta)=1-\omega$. 
The ideal class group $\Cl_k$ is isomorphic to $(\mathbb{Z}/2\mathbb{Z})^2$. 
Since $f(x)\equiv (x+1)(x+2)(x+3) \pmod{5}$, $p=5$ splits completely in $k$. 
Put $\mathfrak{p}=[5,\omega+1,\theta+1]$. Then $\mathfrak{p}$ is a prime ideal of $\mathcal{O}_k$ and 
$5\mathcal{O}_k=\mathfrak{p} \sigma(\mathfrak{p}) \sigma^2(\mathfrak{p})$. 
The ideal class group $\Cl_k$ is generated by the ideal classes of $\mathfrak{p}$ and $\sigma(\mathfrak{p})$. 
The ideal $\mathfrak{p}^2$ is a principal ideal generated by $\omega+1$ and $N_{k/\mathbb{Q}}(\omega+1)=5^2=N(\mathfrak{p})^2$. 
If we put $\alpha_1=5$, $\alpha_2=\omega+1$, $\alpha_3=\theta+1$ and $\delta=\omega+1$, then 
\begin{alignat*}{3}
\alpha_1^2\delta^{-1} &=15-\omega-\theta,  & \quad  \alpha_1\alpha_2\delta^{-1} &=5,  & \quad \alpha_1\alpha_3\delta^{-1}&=5\theta,   \\
\alpha_2^2\delta^{-1}&=1+\omega,  & \quad \alpha_2\alpha_3\delta^{-1}&=1+\theta, & \quad  \alpha_3^2\delta^{-1}&=14-\omega-10\theta. 
\end{alignat*}
We put
\[ A=\left(
\begin{array}{ccc}
-1 & 0& 5 \\
 0 & 0& 1 \\
 5 & 1& -10
\end{array}
\right),\quad 
B=\left(
\begin{array}{ccc}
 -1 & 0 & 0\\
  0 & 1 & 0 \\
  0 & 0 & -1
\end{array}
\right).
\]
Then $\hat{F}_{(A,B)}(u)=u_1^3+14u_1^2u_2+11u_1u_2^2-u_2^3$ 
and the triple $(\mathcal{O}_k,\mathfrak{p},\delta)$ corresponds to the $\Gamma$-orbit of $(A,B)$.
We also have 
\begin{equation}\label{eq:ISOTROPY-IDEAL2}
(\delta^{-1}\alpha_i\alpha_j)
=C+\omega B+ \theta A, 
\end{equation}
where we put $C=\left(
\begin{array}{ccc}
15 & 5 &  0 \\
 5 & 1 &  1 \\
 0 & 1 & 14
\end{array}
\right)$.
Take an element $\gamma \in \Gamma$ and assume that $\gamma\cdot (A,B)=(A,B)$. 
We write $\gamma=(\gamma_1,\gamma_2)$, $\gamma_1 \in \Gamma_1$, $\gamma_2=\left(
\begin{array}{cc}
p & q \\
r & s
\end{array}
\right) \in \Gamma_2$. Then 
\begin{equation}\label{eq:ISOTROPY-EX1}
 (\gamma_1A\,{}^t\!\gamma_1,\gamma_1B\,{}^t\!\gamma_1)=(s A+q B, r A+p B).
\end{equation}
If we put ${}^t\!(\beta_1,\beta_2,\beta_3)=\gamma_1 {}^t\!(\alpha_1,\alpha_2,\alpha_3)$, 
then $\{\beta_1,\beta_2,\beta_3\}$ is also a $\mathbb{Z}$-basis of $\mathfrak{p}$. 
It follows from \eqref{eq:ISOTROPY-IDEAL2} and \eqref{eq:ISOTROPY-EX1} that
\begin{align}\label{eq:ISOTROPY-EX2}
(\delta^{-1}\beta_i\beta_j)
 &= \gamma_1(\delta^{-1}\alpha_i\alpha_j){}^t\!\gamma_1
  = \gamma_1C\,{}^t\!\gamma_1+\omega \gamma_1B{}^t\!\gamma_1+\theta \gamma_1A{}^t\!\gamma_1 \\
 &= \gamma_1C\,{}^t\!\gamma_1+\omega (rA+p B)+\theta(sA+qB). \nonumber
\end{align}
Since $\gamma_2 \cdot \hat{F}_{(A,B)}=\hat{F}_{(A,B)}$ 
and the isotropy group of $\hat{F}_{(A,B)}$ in $\Gamma_2$ is a cyclic group of order three  generated by $\xi=\left(
\begin{array}{cc}
-1 & 1 \\
-1 & 0
\end{array}
\right)$, we have $\gamma_2 \in \{1_2, \xi, \xi^2\}$.  
Suppose $\gamma_2\neq 1_2$. We may assume that $\gamma_2=\xi$. Then we have 
\begin{equation}\label{eq:ISOTROPY-IDEAL1}
(\delta^{-1}\beta_i\beta_j)=\gamma_1C\,{}^t\!\gamma_1+\omega(- A - B)+\theta B.
\end{equation}
On the other hand, applying $\sigma$ to the equation \eqref{eq:ISOTROPY-IDEAL2}, we have
\begin{align}\label{eq:ISOTROPY-IDEAL3}
(\sigma(\delta)^{-1}\sigma(\alpha_i)\sigma(\alpha_j))
 &= C+\sigma(\omega) B+ \sigma(\theta) A \\
 &= A+13B+C +\omega(-A-B) + \theta B.\nonumber 
\end{align}
Since the integral matrix $\gamma_1C\,{}^t\!\gamma_1$ in \eqref{eq:ISOTROPY-IDEAL1} is determined by 
the matrices $- A - B$ and $B$ (cf. the proof of \cite[Thererem 4]{BH2}), 
$\gamma_1C\,{}^t\!\gamma_1=A+13B+C$. By \eqref{eq:ISOTROPY-IDEAL1} and \eqref{eq:ISOTROPY-IDEAL3}, we have 
$\delta^{-1}\beta_i\beta_j=\sigma(\delta)^{-1}\sigma(\alpha_i)\sigma(\alpha_j)$ for all $i$, $j$. 
In particular, putting $i=1$ and $\kappa=\frac{\delta\sigma(\alpha_1)}{\sigma(\delta)\beta_1}$, we have
$\beta_j=\kappa \sigma(\alpha_j)$ for all $j$. 
Hence $\mathfrak{p}=\kappa \sigma(\mathfrak{p})$. 
This contradicts the fact that the ideal classes of $\mathfrak{p}$ and $\sigma(\mathfrak{p})$ generate 
$\Cl_k$ which is isomorphic to $(\mathbb{Z}/2\mathbb{Z})^2$. 
Thus we have proved that $\gamma_2=1_2$. 
By \eqref{eq:ISOTROPY-EX1}, we have $\gamma_1A\,{}^t\! \gamma_1=A$ and $\gamma_1B\,{}^t\! \gamma_1=B$. 
This implies that $\gamma_1 AB^{-1}=AB^{-1} \gamma_1$. 
Since the matrix $AB^{-1}$ has three distinct eigenvalues which are the conjugates of $-\theta$ over $\mathbb{Q}$, 
we see that $\gamma_1$ is  expressed as a polynomial of $AB^{-1}$ over $\mathbb{Q}$. 
Hence $\gamma_1=g(AB^{-1})$ for some $g(x) \in \mathbb{Q}[x]$ with $\deg g(x) \leq 2$. 
Since $\gamma_1$ is integral over $\mathbb{Z}$, we have $g(x) \in \mathbb{Z}[x]$. 
Let $\rho:k \rightarrow M_3(\mathbb{Q})$ be the regular representation of $k$ over $\mathbb{Q}$ 
with respect to the basis $\{\alpha_1,\alpha_2,\alpha_3\}$, i. e. 
\[ \lambda\, {}^t\!(\alpha_1,\alpha_2,\alpha_3)=\rho(\lambda)\, {}^t\!(\alpha_1,\alpha_2,\alpha_3), \quad \lambda \in k.\]
Then we have $\rho(-\theta)=AB^{-1}$, hence $\gamma_1=\rho(\varepsilon)$ for $\varepsilon=g(-\theta) \in \mathcal{O}_k$. 
The equation $N_{k/\mathbb{Q}}(\varepsilon)=\det \gamma_1=1$ implies $\varepsilon \in U(\mathcal{O}_k)$. 
It follows from $\gamma_1=\rho(\varepsilon)$ and ${}^t\!(\beta_i)=\gamma_1 {}^t\!(\alpha_i)$ that 
$\beta_i=\varepsilon \alpha_i$ for $i=1,2,3$. 
By \eqref{eq:ISOTROPY-EX2} with $\gamma_2=1_2$, we have
\[ (\delta^{-1}\varepsilon^2 \alpha_i\alpha_j)=\gamma_1C\,{}^t\!\gamma_1+\omega B+\theta A.\]
Then \eqref{eq:ISOTROPY-IDEAL2} implies $\delta^{-1}\varepsilon^2 \alpha_i\alpha_j=\delta^{-1}\alpha_i\alpha_j$ for all $i, j$. 
Thus $\varepsilon^2=1$, $\varepsilon=\pm 1$. Since $N_{k/\mathbb{Q}}(\varepsilon)=1$, we have $\varepsilon=1$, $\gamma_1=1_3$. 
This proves that the isotropy group $\Gamma_{(A,B)}$ is trivial. 
\end{ex}

\begin{rem}\label{rem:ISOTROPY5}
Though the action of $\Gamma$ on $L$ and on $\hat{L}$ are slightly different, 
it does not effect the $\Gamma$-orbits and the isotropy groups. 
\end{rem}

For any given nondegenerate cubic ring $\mathcal{O}$, we set $k=\mathcal{O} \otimes_\mathbb{Z} \mathbb{Q}$. 
We denote by $\mathscr{I}(\mathcal{O})$ the set of pairs $(\mathfrak{a},\delta)$ such that 
the triples $(\mathcal{O},\mathfrak{a},\delta)$ are as in Theorem \ref{thm:BHARGAVA1}. 
We denote by $\mathcal{E}(\mathcal{O})$ the set of cubic rings $\mathcal{O}_0$ such that 
$\mathcal{O} \subset \mathcal{O}_0 \subset \mathcal{O}_k$. 
Two elements $\mathcal{O}_0, \mathcal{O}'_0 \in \mathcal{E}(\mathcal{O})$ are called \textit{equivalent} 
if there exists an automorphism $\phi\in \Aut(\mathcal{O}_k)$ such that $\phi(\mathcal{O})=\mathcal{O}$ and $\phi(\mathcal{O}_0)=\mathcal{O}'_0$. 
We write $\mathcal{O}_0 \sim  \mathcal{O}'_0$ in that case. 
The notation $\mathcal{E}(\mathcal{O})/\sim$ means a set of representatives of the equivalence classes in $\mathcal{E}(\mathcal{O})$. 
For any $\mathcal{O}_0 \in \mathcal{E}(\mathcal{O})$, we put 
\[ \mathscr{I}(\mathcal{O},\mathcal{O}_0) =\left\{ (\mathfrak{a},\delta) \in \mathscr{I}(\mathcal{O})\,|\, \End_\mathcal{O}(\mathfrak{a}) \sim \mathcal{O}_0 \right\}.\]
Then it is clear that
\begin{equation}\label{eq:IO-DECOMP}
 \mathscr{I}(\mathcal{O})=\bigcup_{\mathcal{O}_0 \in \mathcal{E}(\mathcal{O})/\sim} \mathscr{I}(\mathcal{O},\mathcal{O}_0) \qquad (\text{disjoint}).
\end{equation}
We say that two elements $(\mathfrak{a},\delta)$ and $(\mathfrak{a}',\delta')$ in $\mathscr{I}(\mathcal{O})$ are \textit{equivalent} 
if there exists an invertible element $\kappa$ of $k$ such that 
$\mathfrak{a}'=\kappa\mathfrak{a}$ and $\delta'=\kappa^2\delta$. 
We note that this is stronger than that $(\mathcal{O},\mathfrak{a},\delta)$ and $(\mathcal{O},\mathfrak{a}',\delta')$ are equivalent. 
We defined the subset $\hat{L}(\mathcal{O})$ of $\hat{L}$ in \S~1. 
It is clear that $\hat{L}(\mathcal{O})$ coincides with the set of pairs $(A,B) \in \hat{L}$ such that 
the $\Gamma$-orbit of $(A,B)$ corresponds to the equivalence class of a triplet $(\mathcal{O},\mathfrak{a},\delta)$ for some $(\mathfrak{a},\delta) \in \mathscr{I}(\mathcal{O})$.
We denote by $\hat{L}(\mathcal{O},\mathcal{O}_0)$ 
the set of pairs $(A,B) \in \hat{L}$ such that 
the $\Gamma$-orbit of $(A,B)$ corresponds to the equivalence class of a triplet $(\mathcal{O},\mathfrak{a},\delta)$ for some $(\mathfrak{a},\delta) \in \mathscr{I}(\mathcal{O},\mathcal{O}_0)$. 
By \eqref{eq:IO-DECOMP}, we have
\begin{equation}\label{eq:HATLO-DECOMP}
\hat{L}(\mathcal{O})=\bigcup_{\mathcal{O}_0 \in \mathcal{E}(\mathcal{O})/\sim} \hat{L}(\mathcal{O},\mathcal{O}_0) \qquad (\text{disjoint}).
\end{equation}
We recall that $\{1,\omega,\theta\}$ is a normalized basis of $\mathcal{O}$ and 
the ring structure of $\mathcal{O}$ is given by \eqref{eq:CUBIC-RING-STR}. 
We set 
\begin{align*}
F(u)&=au_1^3+bu_1^2u_2+cu_1u_2^2+du_2^3,\\
\hat{L}(F) &=\{(A,B) \in \hat{L}(\mathcal{O})\,|\, \hat{F}_{(A,B)}=F\}, \\
\hat{L}(F,\mathcal{O}_0) &=\{(A,B) \in \hat{L}(\mathcal{O},\mathcal{O}_0)\,|\, \hat{F}_{(A,B)}=F.\}
\end{align*}
Then we have
\begin{align*}
\hat{L}(\mathcal{O})
 &=\{\gamma \cdot (A,B)\,|\, (A,B) \in \hat{L}(F),\; \gamma \in \Gamma\}, \\
\hat{L}(\mathcal{O},\mathcal{O}_0)
 &=\{\gamma \cdot (A,B)\,|\, (A,B) \in \hat{L}(F,\mathcal{O}_0),\; \gamma \in \Gamma\}.
\end{align*}
\begin{lem}\label{lem:GAMMA1-ORBITS}
For any $(A,B) \in \hat{L}(\mathcal{O},\mathcal{O}_0)$, we set
\[ S(A,B)=\Gamma (A,B) \cap \hat{L}(F,\mathcal{O}_0).\]
Then the group $\Gamma_1$ acts on $S(A,B)$ and the number of $\Gamma_1$-orbits in $S(A,B)$ equals $|\Aut(\mathcal{O})|\cdot|U^+_2(\mathcal{O}_0)|/|\Gamma_{(A,B)}|$.
\end{lem}
\Proof
We may assume $\hat{F}_{(A,B)}=F$. 
Let $\proj_2:\Gamma_{(A,B)}\rightarrow \Gamma_{2,F}$ be the projection to the second component, 
i.e. $\proj_2(\gamma_1,\gamma_2)=\gamma_2$. 
Here we denote by $\Gamma_{2,F}$ the isotropy group of $F$ in $\Gamma_2$. 
It is clear that $\Gamma_1$ acts on $S(A,B)$. 
We can take a representative of the form $(1_3,\gamma_2) \cdot (A,B)$ with $\gamma_2 \in \Gamma_{2,F}$ 
for any $\Gamma_1$-orbit in $S(A,B)$. 
For any two elements $\gamma_2, \gamma'_2 \in \Gamma_{2,F}$, 
$\Gamma_1(1_3,\gamma_2) \cdot (A,B)=\Gamma_1(1_3,\gamma'_2) \cdot (A,B)$ if and only if 
$(\gamma_1,\gamma_2^{-1}\gamma'_2) \in \Gamma_{(A,B)}$ for some $\gamma_1 \in \Gamma_1$. 
This is equivalent to $\gamma_2^{-1}\gamma'_2 \in \proj_2(\Gamma_{(A,B)})$. 
Hence the number of $\Gamma_1$-orbits in $S(A,B)$ equals 
$|\Gamma_{2,F}/\proj_2(\Gamma_{(A,B)})|$. Since $\Gamma_{2,F}$ is isomorphic to $\Aut(\mathcal{O})$, 
it follows from Corollary \ref{cor:ISOTROPY2} that $\ker(\proj_2)$ is isomorphic to $U^+_2(\mathcal{O}_0)$. 
Hence we have
\[ \frac{|\Gamma_{2,F}|}{|\proj_2(\Gamma_{(A,B)})|}
=\frac{|\Gamma_{2,F}|\cdot|\ker \proj_2|}{|\Gamma_{(A,B)}|}
=\frac{|\Aut(\mathcal{O})|\cdot|U^+_2(\mathcal{O}_0)|}{|\Gamma_{(A,B)}|}.
\]
\qed

Using the surjective mapping $\Gamma_1\backslash \hat{L}(F,\mathcal{O}_0)\rightarrow \Gamma \backslash \hat{L}(\mathcal{O},\mathcal{O}_0)$, 
Lemma \ref{lem:GAMMA1-ORBITS} implies that 
\begin{equation}\label{eq:SUM-GAMMA-BS-HATLOOZERO1}
\frac{1}{|\Aut(\mathcal{O})|}\,|\Gamma_1\backslash \hat{L}(F,\mathcal{O}_0)|
=|U^+_2(\mathcal{O}_0)|\sum_{(A,B) \in \Gamma \backslash \hat{L}(\mathcal{O},\mathcal{O}_0)} \frac{1}{|\Gamma_{(A,B)}|}.
\end{equation}
We denote by $\mathscr{C}(\mathcal{O},\mathcal{O}_0)$ the set of equivalence classes 
of $\mathscr{I}(\mathcal{O},\mathcal{O}_0)$. 
The mapping $(\mathfrak{a},\delta)\mapsto \Gamma(A,B)$ induces a mapping
of $\mathscr{C}(\mathcal{O},\mathcal{O}_0)$ to $\Gamma_1\backslash \hat{L}(F,\mathcal{O}_0)$, 
where $\Gamma (A,B)$ is the $\Gamma$-orbit corresponding to the equivalence class of the triple $(\mathcal{O},\mathfrak{a},\delta)$.  
It is easy to see that the mapping 
$\mathscr{C}(\mathcal{O},\mathcal{O}_0) \rightarrow \Gamma_1\backslash \hat{L}(F,\mathcal{O}_0)$
is bijective. By \eqref{eq:SUM-GAMMA-BS-HATLOOZERO1}, we have
\begin{equation}\label{eq:SUM-GAMMA-BS-HATLOOZERO2}
\frac{1}{|\Aut(\mathcal{O})|}|\mathscr{C}(\mathcal{O},\mathcal{O}_0)|
=|U^+_2(\mathcal{O}_0)|\sum_{(A,B) \in \Gamma\backslash \hat{L}(\mathcal{O},\mathcal{O}_0)} \frac{1}{|\Gamma_{(A,B)}|}.
\end{equation}

From now on, we restrict ourselves to the case when the index 
$f=(\mathcal{O}_k:\mathcal{O})$ is square free. 
By using some results of Fr\"{o}hlich \cite{FRO}, we shall give 
formulae for $|\mathscr{C}(\mathcal{O},\mathcal{O}_0)|$ in terms of 
$|\Cl_{\mathcal{O}_0}^{(2)}|$. The key is to investigate 
the ideal $\delta^{-1}\mathfrak{a}^2$ for $(\mathfrak{a},\delta) \in \mathscr{I}(\mathcal{O},\mathcal{O}_0)$. 

Let $k$ be an \'{e}tale cubic algebra over $\mathbb{Q}$ and 
let $\mathcal{O}$ be an order of $k$ with square free index $f=(\mathcal{O}_k:\mathcal{O})$. 
As we have seen in \S 2, we can take a normalized basis $\{1,\omega,\theta\}$ of $\mathcal{O}_k$ with $\omega, \theta \in k^\times$ 
such that $\{1,f\omega,\theta\}$ is a basis of $\mathcal{O}$ 
and the ring structure of $\mathcal{O}_k$ is given by \eqref{eq:CUBIC-RING-STR}. 
We note that $f$ divides $d$. 
Let $(\mathcal{O},\mathfrak{a},\delta)$ be a triple, where $\mathfrak{a}$ is 
a fractional ideal of $\mathcal{O}$ and $\delta$ is an invertible element of $k$ such that 
$\mathfrak{a}^2\subset \delta\mathcal{O}$ and $N_{k/\mathbb{Q}}(\delta)=N_{\mathcal{O}}(\mathfrak{a})^2$. 
We put $\mathcal{O}_0=\End_\mathcal{O}(\mathfrak{a})$. 
Since $\mathcal{O} \subset \mathcal{O}_0 \subset \mathcal{O}_k$ 
and $\mathcal{O}_k/\mathcal{O} \cong \mathbb{Z}/f\mathbb{Z}$, 
$\mathcal{O}_0=[1,g\omega,\theta]$ for some positive divisor $g$ of $f$. We write $f=gh$. 
By Proposition \ref{prop:CUNDUCTOR3}, the conductor $\mathfrak{f}$ of $\mathcal{O}$ is 
given by $\mathfrak{f}=[f,f\omega,\theta]$. Similarly 
the conductor $\mathfrak{g}$ of $\mathcal{O}_0$ is given by $\mathfrak{g}=[g,g\omega,\theta]$.
We put $\mathfrak{j}=[h,f\omega,\theta]$. 
Then $\mathfrak{j}\subset \mathcal{O}$ and $\mathfrak{j}$ is an $\mathcal{O}_0$-ideal. 
In fact, 
\begin{align*}
g\omega\mathfrak{j}
 &=[gh\omega,fg\omega^2,g\omega\theta]
 =[f\omega,fg(-ac+b\omega-a\theta),-adg]\\
 &\subset [f,f\omega,\theta]\subset [h,f\omega,\theta]=\mathfrak{j}, \\
\theta\mathfrak{j}
 &=[h\theta,f\omega\theta,\theta^2]=[h\theta,-adf,-bd+d\omega-c\theta] \\
 &\subset [f,f\omega,\theta]\subset [h,f\omega,\theta]=\mathfrak{j}.
\end{align*}
Since $f=gh$ is square free, $g$ and $h$ are coprime to each other.
Hence $\mathfrak{j}+\mathfrak{g}=\mathcal{O}_0$. 
We put $\mathfrak{h}=[h,h\omega,\theta]$ and $\mathfrak{h}'=[h,\omega,\theta+c]$. 
As we have seen before Proposition \ref{prop:CUNDUCTOR3}, 
$\mathfrak{h}$ and $\mathfrak{h}'$ are integral ideals of $\mathcal{O}_k$ 
such that $\mathfrak{h}\mathfrak{h}'=h\mathcal{O}_k$. 
So $\mathfrak{h}$ is an invertible $\mathcal{O}_k$-ideal. 
Further we have $\mathfrak{j}\mathcal{O}_k=\mathfrak{h}$. 
By Corollary \ref{cor:PRIMETOCOND-INVERTIBLE}, $\mathfrak{j}$ is an invertible $\mathcal{O}_0$-ideal. 
By the same argument as just before Proposition \ref{prop:CUNDUCTOR3}, 
we see that $\mathfrak{j}$ is the largest $\mathcal{O}_0$-ideal contained in $\mathcal{O}$. 
So we write $\mathfrak{j}=\mathfrak{j}(\mathcal{O},\mathcal{O}_0)$. 

We denote by $\widehat{\mathcal{O}}_0$ the dual lattice of $\mathcal{O}_0$ 
with respect to the trace pairing. So we have 
$\widehat{\mathcal{O}}_0=\{\lambda \in k\,|\, \trace(\lambda \mathcal{O}_0)\subset \mathbb{Z}\}$. 
We put $\xi=-\omega/a$ and $\Delta=3a\xi^2+2b\xi+c \in k$.
Since $N_{k/\mathbb{Q}}(\Delta)=-a^{-1}\Disc(\mathcal{O}_k)\neq 0$, $\Delta$ is an invertible element of $k$.
The Lagrange's interpolation formula implies that 
the dual basis of $\{1,\xi,\xi^2\}$ with respect to the trace pairing is
\[ \left\{\Delta^{-1}(a\xi^2+b\xi+c),\Delta^{-1}(a\xi+b),\Delta^{-1} a \right\}, \]
which equals to 
$\{-\Delta^{-1}\theta,\Delta^{-1}(b-\omega),\Delta^{-1}a\}$ by \eqref{eq:CUBIC-RING-STR}.
Since $\mathcal{O}_0=[1,g(\omega-b),\theta]$, 
we see that $\widehat{\mathcal{O}}_0=\Delta^{-1}[1,g^{-1}\xi,a\xi^2]$.
It is easy to see that 
$[ag^2,g(\omega-b),\theta]$ is an integral $\mathcal{O}_0$-ideal. 
Moreover we have
\begin{align*}
\widehat{\mathcal{O}}_0\Delta[ag^2,g(\omega-b),\theta]
 &=[1,g^{-1}\xi,a\xi^2][ag^2,g(\omega-b),\theta] \\
 &=[ag^2,bg,c,g^{-1}d,g\omega,\theta].
\end{align*}
Since the index $(\mathcal{O}_k:\mathcal{O}_0)=g$ is square free, 
the binary cubic form $ag^2u_1^3+bgu_1^2u_2+cu_1u_2^3+(d/g)u_2^3$ 
corresponding to $\mathcal{O}_0$ is primitive. 
Hence we have $\widehat{\mathcal{O}}_0\Delta[ag^2,g(\omega-b),\theta]=[1,g\omega,\theta]=\mathcal{O}_0$.
This proves that $\widehat{\mathcal{O}}_0$ is an invertible $\mathcal{O}_0$-ideal and 
$\widehat{\mathcal{O}}_0^{-1}=\Delta[ag^2,g(\omega-b),\theta]$.

We now show that $\widehat{\mathcal{O}}_0\cap (f/\xi)\widehat{\mathcal{O}}_0
=\Delta^{-1}\mathfrak{j}$. 
We have 
\begin{align*}
\widehat{\mathcal{O}}_0\cap (f/\xi)\widehat{\mathcal{O}}_0
 &=\Delta^{-1}[1,g^{-1}\xi,a\xi^2]\cap \Delta^{-1}[(f/\xi),f/g,af\xi] \\
 &=\Delta^{-1}\left( [1,g^{-1}\xi,a\xi^2]\cap [h,af\xi,(f/d)(a\xi^2+b\xi+c)]\right).
\end{align*}
If $x+yg^{-1}\xi+za\xi^2=x'h+y'af\xi+z'(f/d)(a\xi^2+b\xi+c)$, $x,y,z\in \mathbb{Z}$, 
$x',y',z' \in \mathbb{Z}$, then 
$z'=(d/f)z$, $y=afgy'+bgz$ and $x=hx'+cz$. This proves that 
$\widehat{\mathcal{O}}_0 \cap (f/\xi)\widehat{\mathcal{O}}_0
=\Delta^{-1}[h,f\omega,\theta]=\Delta^{-1}\mathfrak{j}$.

We put $\omega_0=1$, $\omega_1=f\omega$, $\omega_2=\theta$. 
We take a basis $\{\alpha_1,\alpha_2,\alpha_3\}$ of $\mathfrak{a}$ 
having the same orientation as $\{\omega_0,\omega_1,\omega_2\}$. 
Thus there exists a matrix $T\in \GL_3(\mathbb{Q})$ with $\det T>0$ such that
${}^t (\alpha_1, \alpha_2, \alpha_3)=T\,{}^t(\omega_0,\omega_1,\omega_2)$. 
Since $\mathfrak{a}^2\subset \delta\mathcal{O}$, there exist integers $a_{ij}$, $b_{ij}$ and $c_{ij}$ such that 
\[ \alpha_i\alpha_j=\delta(c_{ij}+b_{ij}f\omega+a_{ij}\theta). \]
In Theorem \ref{thm:BHARGAVA1}, the corresponding pair $(A,B)$ is given by
$A=(a_{ij})$ and $B=(b_{ij})$.  We put 
\begin{equation}\label{eq:A0B0F}
A_0=\left(
\begin{array}{ccc}
 0 &  0   & 1 \\
 0 & -af^2&  0 \\
 1 &  0   & -c
\end{array}
\right), \quad 
B_0=\left(
\begin{array}{ccc}
0 & 1 & 0 \\
1 & bf& 0 \\
0 & 0 & d/f
\end{array}
\right).
\end{equation}
Then we have 
\[ (A,B)=(T\rho(\delta^{-1})A_0{}^tT,T\rho(\delta^{-1})B_0{}^tT), \]
where $\rho:k\rightarrow M_3(\mathbb{Q})$ is the regular representation of $k$ with 
respect to the basis $\{\omega_0,\omega_1,\omega_2\}$. 
Since the dual basis of $\{1,\xi,\xi^2\}$ with respect to the trace pairing is 
$\{-\Delta^{-1}\theta,\Delta^{-1}(b-\omega),\Delta^{-1}a\}$, we have
\[ \left(\trace(\omega_i\omega_j\Delta^{-1}) \right)_{0\leq i,j\leq 2}=-A_0,\quad
   \left(\trace(\omega_i\omega_j\Delta^{-1}(\xi/f)) \right)_{0\leq i,j\leq 2}=-B_0.
\]
We also have
\[ (\Delta^{-1}\delta^{-1}\alpha_i\alpha_j)
  =T\rho(\delta^{-1}) \left(\omega_i\omega_j\Delta^{-1} \right){}^tT.
\]
Hence 
\[ \left(\trace(\Delta^{-1}\delta^{-1}\alpha_i\alpha_j) \right)
=-T\rho(\delta^{-1})A_0{}^tT=-A.
\]
We have similarly
\[ \left(\trace(\Delta^{-1}\delta^{-1}\alpha_i\alpha_j(\xi/f)) \right)
=-T\rho(\delta^{-1})B_0{}^tT=-B.
\]
Since $\delta^{-1}\mathfrak{a}^2$ is generated by $\{\delta^{-1}\alpha_i\alpha_j\}$ over $\mathbb{Z}$,
we have $\trace(\alpha\Delta^{-1}) \in \mathbb{Z}$ and $\trace(\alpha\Delta^{-1}(\xi/f)) \in \mathbb{Z}$
for all $\alpha \in \delta^{-1}\mathfrak{a}^2$. 
Since $\delta^{-1}\mathfrak{a}^2$ is an $\mathcal{O}_0$-ideal, we have 
$\lambda\alpha \in \delta^{-1}\mathfrak{a}^2$ for any $\lambda \in \mathcal{O}_0$ and 
$\alpha \in \delta^{-1}\mathfrak{a}^2$. Hence 
\[ \trace(\lambda\alpha\Delta^{-1}) \in \mathbb{Z}, \quad \trace(\lambda\alpha\Delta^{-1}(\xi/f)) \in \mathbb{Z}. \]
This proves that $\alpha\Delta^{-1} \in \widehat{\mathcal{O}}_0$ and 
$\alpha\Delta^{-1}(\xi/f) \in \widehat{\mathcal{O}}_0$ for any $\alpha \in \delta^{-1}\mathfrak{a}^2$. 
Hence $\Delta^{-1}\delta^{-1}\mathfrak{a}^2 \subset \widehat{\mathcal{O}}_0\cap (f/\xi)\widehat{\mathcal{O}}_0=\Delta^{-1}\mathfrak{j}$. 
So we have $\delta^{-1}\mathfrak{a}^2\subset \mathfrak{j}$. 
We say that a cubic ring $R$ is \textit{weakly self dual} if every proper $R$-ideal is an invertible $R$-ideal. 
Since the index $g=(\mathcal{O}_k:\mathcal{O}_0)$ is square free, 
$\mathcal{O}_0$ is weakly self dual by Fr\"{o}hlich \cite[Therem 10, Proposition 8.1]{FRO}. 
Hence any proper $\mathcal{O}_0$-ideal is an invertible $\mathcal{O}_0$-ideal. 
In particular, $\mathfrak{a}$ is an invertible $\mathcal{O}_0$-ideal. 
We denote by $N_{\mathcal{O}_0}(\mathfrak{a})$ the norm of $\mathfrak{a}$ as 
an $\mathcal{O}_0$-ideal. Then we have 
\[ N_{\mathcal{O}_0}(\mathfrak{a})=(\mathcal{O}_0:\mathcal{O})N_\mathcal{O}(\mathfrak{a})=hN_\mathcal{O}(\mathfrak{a}).\]
It follows from \cite[Corollary 1 to Theorem 4]{FRO} that 
an $\mathcal{O}_0$-ideal $\mathfrak{b}$ is invertible if and only if 
$(\mathcal{O}_k\mathfrak{b}:\mathfrak{b})=(\mathcal{O}_k:\mathcal{O}_0)$. 
Hence for any proper integral $\mathcal{O}_0$-ideal $\mathfrak{b}$, 
we have 
\begin{align*}
N_{\mathcal{O}_0}(\mathfrak{b})
 &=(\mathcal{O}_0:\mathfrak{b}) 
 =\frac{(\mathcal{O}_k:\mathcal{O}_0)(\mathcal{O}_0:\mathfrak{b})}{(\mathcal{O}_k:\mathcal{O}_0)}
 =\frac{(\mathcal{O}_k:\mathcal{O}_k\mathfrak{b})(\mathcal{O}_k\mathfrak{b}:\mathfrak{b})}{(\mathcal{O}_k:\mathcal{O}_0)}\\
 &=(\mathcal{O}_k:\mathcal{O}_k\mathfrak{b})=N(\mathcal{O}_k\mathfrak{b}).
\end{align*}
Thus for any two integral proper $\mathcal{O}_0$-ideals $\mathfrak{b}_1$ and $\mathfrak{b}_2$, we have
\[ N_{\mathcal{O}_0}(\mathfrak{b}_1\mathfrak{b}_2)=N(\mathcal{O}_k\mathfrak{b}_1\mathcal{O}_k\mathfrak{b}_2)
=N(\mathcal{O}_k\mathfrak{b}_1)N(\mathcal{O}_k\mathfrak{b}_2)
=N_{\mathcal{O}_0}(\mathfrak{b}_1)N_{\mathcal{O}_0}(\mathfrak{b}_2).
\]
In particular, we have 
\[ N_{\mathcal{O}_0}(\delta^{-1}\mathfrak{a}^2)=N_{k/\mathbb{Q}}(\delta)^{-1}N_{\mathcal{O}_0}(\mathfrak{a})^2
=N_{k/\mathbb{Q}}(\delta)^{-1}h^2N_\mathcal{O}(\mathfrak{a})^2=h^2.
\]
Since $\mathfrak{j}=[h,f\omega,\theta]$ and $\mathcal{O}_0=[1,g\omega,\theta]$, 
we have $N_{\mathcal{O}_0}(\mathfrak{j})=h^2$. 
Then the inclusion $\delta^{-1}\mathfrak{a}^2\subset \mathfrak{j}$ implies 
$\delta^{-1}\mathfrak{a}^2=\mathfrak{j}$. 
Conversely, if $\mathfrak{a}$ is an invertible $\mathcal{O}_0$-ideal and 
$\delta$ is an invertible element of $k$ with positive norm 
such that $\delta^{-1}\mathfrak{a}^2=\mathfrak{j}$, 
then we have $\delta^{-1}\mathfrak{a}^2\subset \mathcal{O}$ and 
$N_{\mathcal{O}_0}(\mathfrak{a})^2=h^2 N_{k/\mathbb{Q}}(\delta)$, hence 
$N_\mathcal{O}(\mathfrak{a})^2=N_{k/\mathbb{Q}}(\delta)$. 
%
%

Since $\mathcal{O}_k=[1,\omega,\theta]$ and $\mathcal{O}=[1,f\omega, \theta]$, 
for any $\phi \in \Aut(\mathcal{O}_k)$ with $\phi(\mathcal{O})=\mathcal{O}$, 
we have $\phi(\mathcal{O}_0)=\mathcal{O}_0$. 
Hence $\mathcal{E}(\mathcal{O})=\mathcal{E}(\mathcal{O})/\sim$. 
We have seen that 
\[ \mathscr{I}(\mathcal{O},\mathcal{O}_0)
=\{(\mathfrak{a},\delta) \in I_{\mathcal{O}_0}\times k^\times\,|\, 
\delta^{-1}\mathfrak{a}^2=\mathfrak{j}(\mathcal{O},\mathcal{O}_0),\;N_{k/\mathbb{Q}}(\delta)>0\}.
\]
We also have
\begin{equation}\label{eq:IO-DECOMP-SF}
\mathscr{I}(\mathcal{O})=\bigcup_{\mathcal{O}\subset \mathcal{O}_0\subset \mathcal{O}_k} \mathscr{I}(\mathcal{O},\mathcal{O}_0) \qquad (\text{disjoint}).
\end{equation}
We define two subgroups of $I_{\mathcal{O}_0} \times k^\times$ by 
\begin{align*}
G(\mathcal{O}_0)
 &=\{(\mathfrak{b},\lambda) \in I_{\mathcal{O}_0}\times k^\times\,|\, 
 \mathfrak{b}^2=\lambda \mathcal{O}_0, \,N_{k/\mathbb{Q}}(\lambda)>0\}, \\
H(\mathcal{O}_0)
 &=\{(\kappa\mathcal{O}_0,\kappa^2)\,|\, \kappa \in k^\times\}.
\end{align*}
$H(\mathcal{O}_0)$ is a subgroup of $G(\mathcal{O}_0)$. 
If the set $\mathscr{I}(\mathcal{O},\mathcal{O}_0)$ is not empty, then the group $G(\mathcal{O}_0)$ acts 
transitively on it 
by $(\mathfrak{b},\lambda)\cdot (\mathfrak{a},\delta)=(\mathfrak{b}\mathfrak{a},\lambda\delta)$, 
hence the number of $H(\mathcal{O}_0)$-orbits of $\mathscr{I}(\mathcal{O},\mathcal{O}_0)$ is given by 
\begin{equation}\label{eq:HOIO}
|H(\mathcal{O}_0)\backslash \mathscr{I}(\mathcal{O},\mathcal{O}_0)|=(G(\mathcal{O}_0):H(\mathcal{O}_0)).
\end{equation}
We denote by $\Cl_{\mathcal{O}_0}^{(2)}$ the subgroup of $\Cl_{\mathcal{O}_0}$ 
consisting of all elements $c$ satisfying $c^2=1$. 
For any ideal class $c \in \Cl_{\mathcal{O}_0}^{(2)}$, we take an invertible $\mathcal{O}_0$-ideal 
$\mathfrak{b}$ in $c$. Then $\mathfrak{b}^2=\lambda \mathcal{O}_0$ for some $\lambda \in k^\times$. 
If $N_{k/\mathbb{Q}}(\lambda)<0$, then we replace $\lambda$ by $-\lambda$.  So we can take $N_{k/\mathbb{Q}}(\lambda)>0$. 
Hence $(\mathfrak{b},\lambda) \in G(\mathcal{O}_0)$. 
Thus the correspondence $(\mathfrak{b},\lambda)\mapsto \mathfrak{b}P_{\mathcal{O}_0}$ 
defines a surjective homomorphism $\varphi:G(\mathcal{O}_0)\rightarrow \Cl_{\mathcal{O}_0}^{(2)}$. 
The kernel is obviously given by 
\[ \ker \varphi=\{(\kappa\mathcal{O}_0,\varepsilon\kappa^2)\,|\, \kappa \in k^\times,\,
\varepsilon \in U^+(\mathcal{O}_0)\}, 
\]
where $U^+(\mathcal{O}_0)$ denotes the group of units of $\mathcal{O}_0$ having positive norm. 
We have $\ker \varphi/H(\mathcal{O}_0) \cong U^+(\mathcal{O}_0)/U^+(\mathcal{O}_0)^2$. 
Hence we have 
\begin{equation}\label{eq:INDEX-G0H0}
(G(\mathcal{O}_0):H(\mathcal{O}_0))=|\Cl_{\mathcal{O}_0}^{(2)}| (U^+(\mathcal{O}_0):U^+(\mathcal{O}_0)^2). 
\end{equation}
We put 
\[F(u)=af^2u_1^3+bfu_1^2u_2+cu_1u_2^2+(d/f)u_2^3. \]
Then $F$ corresponds to the cubic ring $\mathcal{O}=[1,f\omega,\theta]$. 
By definition, we have $\mathscr{C}(\mathcal{O},\mathcal{O}_0)=H(\mathcal{O}_0)\backslash \mathscr{I}(\mathcal{O},\mathcal{O}_0)$. 
By \eqref{eq:SUM-GAMMA-BS-HATLOOZERO2}, \eqref{eq:HOIO} and \eqref{eq:INDEX-G0H0}, we have
\[\frac{|\Cl_{\mathcal{O}_0}^{(2)}| (U^+(\mathcal{O}_0):U^+(\mathcal{O}_0)^2)}{|\Aut(\mathcal{O})|\,|U^+_2(\mathcal{O}_0)|}
=\sum_{(A,B) \in \Gamma\backslash \hat{L}(\mathcal{O},\mathcal{O}_0)} \frac{1}{|\Gamma_{(A,B)}|}
\]
provided that $\mathscr{I}(\mathcal{O},\mathcal{O}_0)$ is not empty.
The set $\mathscr{I}(\mathcal{O},\mathcal{O}_0)$ is not empty if and only if 
the ideal class of $\mathfrak{j}(\mathcal{O},\mathcal{O}_0)$ belongs to $\Cl_{\mathcal{O}_0}^2$. 
We denote by $X(\mathcal{O},\mathcal{O}_0)$ the subgroup of $\Cl_{\mathcal{O}_0}/\Cl_{\mathcal{O}_0}^2$ 
generated by the ideal class of $\mathfrak{j}(\mathcal{O},\mathcal{O}_0)$. 
Then $X(\mathcal{O},\mathcal{O}_0)$ has order two if $\mathscr{I}(\mathcal{O},\mathcal{O}_0)$ is empty, 
otherwise it is trivial. Hence we have
\begin{align}\label{eq:HATLO-ORBITNUMBER}
\Gamma\backslash \hat{L}(\mathcal{O})
 &=\bigcup_{\mathcal{O}\subset \mathcal{O}_0 \subset \mathcal{O}_k} \Gamma\backslash \hat{L}(\mathcal{O},\mathcal{O}_0)\qquad \text{(disjoint)},\nonumber \\
\sum_{(A,B) \in \Gamma\backslash \hat{L}(\mathcal{O},\mathcal{O}_0)} \frac{1}{|\Gamma_{(A,B)}|}
&=\frac{(U^+(\mathcal{O}_0):U^+(\mathcal{O}_0)^2)}{|\Aut(\mathcal{O})|\,|U^+_2(\mathcal{O}_0)|}\,|\Cl_{\mathcal{O}_0}^{(2)}|\, (2-|X(\mathcal{O},\mathcal{O}_0)|).
\end{align}
If $k$ is a cubic field  or $\mathcal{O}_0=\mathcal{O}_k$, then 
\begin{equation}\label{eq:UPLUS-INDEX}
(U^+(\mathcal{O}_0):U^+(\mathcal{O}_0)^2)=\left\{
\begin{array}{ll}
4 , &\quad  \Disc(k)>0, \\
2 , &\quad \Disc(k)<0,
\end{array}
\right.
\end{equation}
where $\Disc(k)$ denotes the discriminant of $\mathcal{O}_k$. 
Further we have
\begin{equation}\label{eq:UPLUS2-ORDER}
|U^+_2(\mathcal{O}_0)|=\left\{
\begin{array}{ll}
1 , &\quad \text{$k$ is a cubic field,} \\
2 , &\quad \text{$k=\mathbb{Q}\oplus k_1$, $k_1$ is a quadratic field,} \\
4 , &\quad \text{$k=\mathbb{Q}^3$}.
\end{array}
\right.
\end{equation}

We now study the case of $\Disc(k)>0$ more precisely. 
For an element $(A,B) \in \hat{L}(F,\mathcal{O}_0)$, 
we shall obtain a necessary and sufficient condition for $(A,B) \in V_1$, 
where $V_1$ is the real orbit defined in \S 1. 
Let $(\mathcal{O},\mathfrak{a},\delta)$ be the corresponding triple. 
Then we have $\delta^{-1}\mathfrak{a}^2=\mathfrak{j}(\mathcal{O},\mathcal{O}_0)$ and 
$(A,B)$ is given by
\[ A=-(\trace(\Delta^{-1}\delta^{-1}\alpha_i\alpha_j)),
\quad B=-(\trace(\Delta^{-1}\delta^{-1}\alpha_i\alpha_j(\xi/f))), 
\]
where $\{\alpha_1,\alpha_2,\alpha_3\}$ is a $\mathbb{Z}$-basis of $\mathfrak{a}$ 
having the same orientation as $\{1,f\omega,\theta\}$. 
For any $\lambda \in k$, we denote by $\lambda^{(i)}$ ($i=1,2,3$) 
the conjugates of $\lambda$. If $k$ is a cubic field, then $\lambda^{(i)}$'s are as usual. 
If $k=\mathbb{Q}\oplus k_1$, $k_1$ is a quadratic field, then 
for $\lambda=(\lambda_1,\lambda_2)$, $\lambda_1\in \mathbb{Q}$, $\lambda_2 \in k_1$, 
the conjugates of $\lambda$ are $\lambda_1$, $\lambda_2$ and $\lambda'_2$, 
where $\lambda'_2$ is the conjugate of $\lambda_2$. 
If $k=\mathbb{Q}^3$, then the conjugates of $\lambda=(\lambda_1,\lambda_2,\lambda_3)$ are 
$\lambda_1$, $\lambda_2$ and $\lambda_3$. We put $U=(\alpha^{(i)}_j)\in \GL_3(\mathbb{R})$ 
and $\beta=-\Delta^{-1}\delta^{-1}$. Then we have
\begin{align*}
A &={}^tU \diag[\beta^{(1)},\beta^{(2)},\beta^{(3)}]\,U,  \\
B &={}^tU \diag[\beta^{(1)}\xi^{(1)}/f,\beta^{(2)}\xi^{(2)}/f,\beta^{(3)}\xi^{(3)}/f]\,U,
\end{align*}
where $\diag[a_1, \ldots, a_n]$ denotes the diagonal matrix 
whose diagonal entries are $a_1, \dots, a_n$.

Let $v=(v_1:v_2:v_3)$ be an element of $\Zero(A,B)\subset \mathbb{P}^2(\mathbb{C})$ 
and put $w=(w_1:w_2:w_3)=v\,{}^tU$.  
Then we have
\begin{align*}
\beta^{(1)}w_1^2+\beta^{(2)}w_2^2+\beta^{(3)}w_3^2 &=0, \\
\xi^{(1)}\beta^{(1)}w_1^2+\xi^{(2)}\beta^{(2)}w_2^2+\xi^{(3)}\beta^{(3)}w_3^2 &=0.
\end{align*}
These equations imply that $w_1^2/\delta^{(1)}=w_2^2/\delta^{(2)}=w_3^2/\delta^{(3)}$. 
Since $(A,B)$ belongs to $V_1$ if and only if $|\Zero(A,B) \cap \mathbb{P}^2(\mathbb{R})|=4$. 
Hence $(A,B) \in V_1$ if and only if $\delta$ is totally positive. 
If $\delta$ is totally positive, then the ideal class 
$\mathfrak{j}(\mathcal{O},\mathcal{O}_0)P_{\mathcal{O}_0,+}$ belongs to $\Cl_{\mathcal{O}_0,+}^2$. 
So we have proved that $\hat{L}(F,\mathcal{O}_0) \cap V_1$ is not empty if and only if 
the ideal class $\mathfrak{j}(\mathcal{O},\mathcal{O}_0)P_{\mathcal{O}_0,+}$ 
belongs to $\Cl_{\mathcal{O}_0,+}^2$. 
We denote by $k_+^\times$ the subgroup of $k^\times$ consisting of 
totally positive elements. 
We put 
\[ \mathscr{I}_+(\mathcal{O},\mathcal{O}_0)= 
\{(\mathfrak{a},\delta) \in I_{\mathcal{O}_0}\times k_+^\times \,|\, 
\delta^{-1}\mathfrak{a}^2=\mathfrak{j}(\mathcal{O},\mathcal{O}_0) \}.
\]
We define a subgroup of $I_{\mathcal{O}_0} \times k_+^\times$ by 
\[ G_+(\mathcal{O}_0)
 =\{(\mathfrak{b},\lambda) \in I_{\mathcal{O}_0}\times k_+^\times\,|\,
 \mathfrak{b}^2=\lambda \mathcal{O}_0\}.
\] 
Then $G_+(\mathcal{O}_0)$ contains $H(\mathcal{O}_0)$. 
If the set $\mathscr{I}_+(\mathcal{O},\mathcal{O}_0)$ is not empty, then the group $G_+(\mathcal{O}_0)$ acts 
transitively on it as before, 
hence the number of $H(\mathcal{O}_0)$-orbits of $\mathscr{I}_+(\mathcal{O},\mathcal{O}_0)$ is given by 
\begin{equation}\label{eq:H0IPLUS0}
|H(\mathcal{O}_0)\backslash \mathscr{I}_+(\mathcal{O},\mathcal{O}_0)|
=(G_+(\mathcal{O}_0):H(\mathcal{O}_0)).
\end{equation}
We denote by $\mathscr{C}_+(\mathcal{O},\mathcal{O}_0)$ the set of equivalence classes of $\mathscr{I}_+(\mathcal{O},\mathcal{O}_0)$. 
Then we have $\mathscr{C}_+(\mathcal{O},\mathcal{O}_0)=H(\mathcal{O}_0)\backslash \mathscr{I}_+(\mathcal{O},\mathcal{O}_0)$ by definition. 
We put $\hat{L}_1(\mathcal{O},\mathcal{O}_0)=\hat{L}(\mathcal{O},\mathcal{O}_0) \cap V_1$. 
We have similarly as in \eqref{eq:SUM-GAMMA-BS-HATLOOZERO2}
\begin{equation}\label{eq:SUM-GAMMA-BS-TOTPOS-HATLOOZERO}
\frac{1}{|\Aut(\mathcal{O})|}\,|\mathscr{C}_+(\mathcal{O},\mathcal{O}_0)|
=|U^+_2(\mathcal{O}_0)|\sum_{(A,B) \in \Gamma\backslash \hat{L}_1(\mathcal{O},\mathcal{O}_0)} \frac{1}{|\Gamma_{(A,B)}|}.
\end{equation}
We denote by $\Cl_{\mathcal{O}_0,+}^{(2)}$ the subgroup of $\Cl_{\mathcal{O}_0,+}$ 
consisting of all elements $c$ satisfying $c^2=1$. 
For any ideal class $c \in \Cl_{\mathcal{O}_0,+}^{(2)}$, 
we take an invertible $\mathcal{O}_0$-ideal 
$\mathfrak{b}$ in $c$. Then $\mathfrak{b}^2=\lambda \mathcal{O}_0$ for 
some $\lambda \in k_+^\times$. 
Hence $(\mathfrak{b},\lambda) \in G_+(\mathcal{O}_0)$. 
Thus the correspondence $(\mathfrak{b},\lambda)\mapsto \mathfrak{b}P_{\mathcal{O}_0,+}$ 
defines a surjective homomorphism 
$\psi:G_+(\mathcal{O}_0)\rightarrow \Cl_{\mathcal{O}_0,+}^{(2)}$. 
The kernel is obviously given by 
\[ \ker \psi=\{(\kappa\mathcal{O}_0,\varepsilon\kappa^2)\,|\, \kappa \in k_+^\times,\,
\varepsilon \in U_+(\mathcal{O}_0)\}.
\]
where $U_+(\mathcal{O}_0)$ denotes the group of totally positive units in $\mathcal{O}_0$.  
We put $H_+(\mathcal{O}_0)=\{(\kappa\mathcal{O}_0,\kappa^2)\,|\,\kappa \in k_+^\times\}$. 
Then we have 
\[ \ker \psi/H_+(\mathcal{O}_0) \cong U_+(\mathcal{O}_0)/U_+(\mathcal{O}_0)^2\]
and
\[ H(\mathcal{O}_0)/H_+(\mathcal{O}_0) \cong k^\times/k_+^\times U_2(\mathcal{O}_0), \]
where $U_2(\mathcal{O}_0)$ denotes the group of units in $\mathcal{O}_0$ having order dividing 2. Hence 
\begin{align*}
(G_+(\mathcal{O}_0):H(\mathcal{O}_0))
 &=\frac{(G_+(\mathcal{O}_0):\ker \psi)(\ker \psi:H_+(\mathcal{O}_0))}{(H(\mathcal{O}_0):H_+(\mathcal{O}_0))}\\
 &=\frac{|\Cl_{\mathcal{O}_0,+}^{(2)}|(U_+(\mathcal{O}_0):U_+(\mathcal{O}_0)^2)}
 {(k^\times:k_+^\times U_2(\mathcal{O}_0))} \\
 &=\frac{|\Cl_{\mathcal{O}_0,+}^{(2)}|(U_+(\mathcal{O}_0):U_+(\mathcal{O}_0)^2)(k^\times_+U_2(\mathcal{O}_0):k^\times_+)}
 {(k^\times:k_+^\times)}.
\end{align*}
Since $\Disc(k)>0$, $U_2(\mathcal{O}_0) \cap k^\times_+=\{1\}$, hence 
$k^\times_+U_2(\mathcal{O}_0)/k^\times_+\cong  U_2(\mathcal{O}_0)$. 
We also have $(k^\times:k_+^\times)=2^3$. 
By Dirichlet's unit theorem, the index $(U_+(\mathcal{O}_0):U_+(\mathcal{O}_0)^2)$ is 
$4$, $2$ or $1$ according as $k$ is a totally real cubic field, 
$k$ is a direct sum of $\mathbb{Q}$ and a real quadratic field 
or $k=\mathbb{Q}^3$. 
Hence we have 
\begin{equation}\label{eq:INDEX-GPLUS0H0}
(G_+(\mathcal{O}_0):H(\mathcal{O}_0))=
2^{-3}|\Cl_{\mathcal{O}_0,+}^{(2)}|(U_+(\mathcal{O}_0):U_+(\mathcal{O}_0)^2)|U_2(\mathcal{O}_0)|.
\end{equation}
We denote by $X_+(\mathcal{O},\mathcal{O}_0)$ the subgroup of $\Cl_{\mathcal{O}_0,+}/\Cl_{\mathcal{O}_0,+}^2$ 
generated by the ideal class of $\mathfrak{j}(\mathcal{O},\mathcal{O}_0)$. 
Then $X_+(\mathcal{O},\mathcal{O}_0)$ has order two if $\mathscr{I}_+(\mathcal{O},\mathcal{O}_0)$ is empty, 
otherwise it is trivial. 
By \eqref{eq:H0IPLUS0}, \eqref{eq:SUM-GAMMA-BS-TOTPOS-HATLOOZERO} and \eqref{eq:INDEX-GPLUS0H0}, we have
\begin{align}\label{eq:HATLO-ORBITNUMBER-TOTPOSITIVE}
\Gamma\backslash \hat{L}_1(\mathcal{O})
 &=\bigcup_{\mathcal{O}\subset \mathcal{O}_0 \subset \mathcal{O}_k} \Gamma\backslash \hat{L}_1(\mathcal{O},\mathcal{O}_0), \nonumber\\
\sum_{(A,B) \in \Gamma\backslash \hat{L}_1(\mathcal{O},\mathcal{O}_0)} \frac{1}{|\Gamma_{(A,B)}|}
&=\frac{(U_+(\mathcal{O}_0):U_+(\mathcal{O}_0)^2)|U_2(\mathcal{O}_0)|}{2^3 |\Aut(\mathcal{O})|\, |U^+_2(\mathcal{O}_0)|}\,|\Cl_{\mathcal{O}_0,+}^{(2)}|
\\
&\quad \times (2-|X_+(\mathcal{O},\mathcal{O}_0)|). \nonumber
\end{align}
If $k$ is a totally real cubic field  or $\mathcal{O}_0=\mathcal{O}_k$, then 
\[ |U_2(\mathcal{O}_0)|=\left\{
\begin{array}{ll}
2 , &\quad \text{if $k$ is a totally real cubic field,} \\
4 , &\quad \text{if $k=\mathbb{Q} \oplus k_1$, $k_1$ is a real quadratic field,} \\
8 , &\quad \text{if $k=\mathbb{Q}^3$}.
\end{array}
\right.
\]
So in that case \eqref{eq:HATLO-ORBITNUMBER-TOTPOSITIVE} becomes
\[ \sum_{(A,B) \in \Gamma\backslash \hat{L}_1(\mathcal{O},\mathcal{O}_0)} \frac{1}{|\Gamma_{(A,B)}|}
=\frac{|\Cl_{\mathcal{O}_0,+}^{(2)}|}{|\Aut(\mathcal{O})|\,|U^+_2(\mathcal{O}_0)|}\,(2-|X_+(\mathcal{O},\mathcal{O}_0)|).
\]
We write $y=(A,B) \in \hat{L}$ and $\mu(y)=1/|\Gamma_{(A,B)}|$. 
We now set 
\[\hat{L}_i(\mathcal{O},\mathcal{O}_0)=\hat{L}(\mathcal{O},\mathcal{O}_0) \cap V_i \quad (i=1,2,3).\]
By \eqref{eq:HATLO-ORBITNUMBER} and \eqref{eq:HATLO-ORBITNUMBER-TOTPOSITIVE}, we have
\begin{prop}\label{prop:LHAT-CLN1}
Let $k$ be an an \'{e}tale cubic algebra over $\mathbb{Q}$ and let $\mathcal{O}$ and $\mathcal{O}_0$ be 
orders of $k$ such that $\mathcal{O}\subset \mathcal{O}_0 \subset \mathcal{O}_k$. 
We assume that the index $(\mathcal{O}_k:\mathcal{O})$ is square free. 
If $\Disc(k)>0$, then we have
\begin{align*}
\lefteqn{\sum_{y \in \Gamma\backslash \hat{L}_1(\mathcal{O},\mathcal{O}_0)} \mu(y)}\\
 &=\frac{(U_+(\mathcal{O}_0):U_+(\mathcal{O}_0)^2)\cdot|U_2(\mathcal{O}_0)|}{2^3|\Aut(\mathcal{O})|\cdot|U_2^+(\mathcal{O}_0)|}\,|\Cl_{\mathcal{O}_0,+}^{(2)}|\,
(2-|X_+(\mathcal{O},\mathcal{O}_0)|), \\
\lefteqn{\sum_{y \in \Gamma\backslash \hat{L}_1(\mathcal{O},\mathcal{O}_0)} \mu(y)
+\sum_{y \in \Gamma\backslash \hat{L}_3(\mathcal{O},\mathcal{O}_0)} \mu(y)}\\
 &=\frac{(U^+(\mathcal{O}_0):U^+(\mathcal{O}_0)^2)}{|\Aut(\mathcal{O})|\cdot|U_2^+(\mathcal{O}_0)|}\,|\Cl_{\mathcal{O}_0}^{(2)}|\, (2-|X(\mathcal{O},\mathcal{O}_0)|).
\end{align*}
If $\Disc(k)<0$, then we have
\[\sum_{y \in \Gamma\backslash \hat{L}_2(\mathcal{O},\mathcal{O}_0)} \mu(y)
=\frac{(U^+(\mathcal{O}_0):U^+(\mathcal{O}_0)^2)}{|\Aut(\mathcal{O})|\cdot|U_2^+(\mathcal{O}_0)|}\,|\Cl_{\mathcal{O}_0}^{(2)}| \,(2-|X(\mathcal{O},\mathcal{O}_0)|).
\]
\end{prop}
If $k$ is a cubic field, then the proposition above becomes as follows.
\begin{cor}\label{cor:LHAT-CLN2}
Let $k$ be a cubic field and let $\mathcal{O}$ and $\mathcal{O}_0$ be 
orders of $k$ such that $\mathcal{O}\subset \mathcal{O}_0 \subset \mathcal{O}_k$ and 
the index $(\mathcal{O}_k:\mathcal{O})$ is square free. 
If $\Disc(k)>0$, then we have
\begin{align*}
\sum_{y \in \Gamma\backslash \hat{L}_1(\mathcal{O},\mathcal{O}_0)} \mu(y)
 &=\frac{|\Cl_{\mathcal{O}_0,+}^{(2)}|(2-|X_+(\mathcal{O},\mathcal{O}_0)|)}{|\Aut(\mathcal{O})|}, \\
\sum_{y \in \Gamma\backslash \hat{L}_1(\mathcal{O},\mathcal{O}_0)} \mu(y)
+\sum_{y \in \Gamma\backslash \hat{L}_3(\mathcal{O},\mathcal{O}_0)} \mu(y)
 &=\frac{4|\Cl_{\mathcal{O}_0}^{(2)}| \,(2-|X(\mathcal{O},\mathcal{O}_0)|)}{|\Aut(\mathcal{O})|}.
\end{align*}
If $\Disc(k)<0$, then we have
\[
\sum_{y \in \Gamma\backslash \hat{L}_2(\mathcal{O},\mathcal{O}_0)} \mu(y)
=2|\Cl_{\mathcal{O}_0}^{(2)}|\,(2-|X(\mathcal{O},\mathcal{O}_0)|).
\]
\end{cor}

We give the following lemma for later use.
\begin{lem}\label{lem:LHAT-CLN3}
Let $k$ be a Galois cubic field and $\mathcal{O}$ and $\mathcal{O}_0$ be 
orders of $k$ such that $\mathcal{O}\subset \mathcal{O}_0 \subset \mathcal{O}_k$. 
We assume that the index $f=(\mathcal{O}_k:\mathcal{O})$ is a square free odd integer
such that each prime number $p$ dividing $f$ 
ramifies in $k/\mathbb{Q}$. 
Then we have $|\Cl_{\mathcal{O}_0}^{(2)}|=|\Cl_k^{(2)}|$, $|\Cl_{\mathcal{O}_0,+}^{(2)}|=|\Cl_{k,+}^{(2)}|$ and 
$X(\mathcal{O},\mathcal{O}_0)$, $X_+(\mathcal{O},\mathcal{O}_0)$ are trivial.
\end{lem}
\Proof
Let $\{1,\omega,\theta\}$ be a normalized basis of $\mathcal{O}_k$ such that 
$\{1,f\omega,\theta\}$ is a basis of $\mathcal{O}$. Then $\mathfrak{f}=[f,f\omega,\theta]$ is the conductor of $\mathcal{O}$. 
We put $g=(\mathcal{O}_k:\mathcal{O}_0)$ and $h=(\mathcal{O}_0:\mathcal{O})$. 
Since each prime number $p|f$ ramifies in $k/\mathbb{Q}$, the prime ideal decomposition of $p\mathcal{O}_k$ is of the form 
$p\mathcal{O}_k=\mathfrak{p}^3$. Hence the $p$-part of $\mathfrak{f}$ is $\mathfrak{f}_p=\mathfrak{p}^2$. 
Then it follows from the Chinese remainder theorem that 
\[ (\mathcal{O}_k/\mathfrak{f})^\times \cong \prod_{p|f} (\mathcal{O}_k/\mathfrak{f}_p)^\times, \quad |(\mathcal{O}_k/\mathfrak{f})^\times|=\prod_{p|f} p(p-1).\]
Since $\mathcal{O}=\mathbb{Z}+\mathfrak{f}$ and $\mathfrak{f}\cap \mathbb{Z}=f\mathbb{Z}$, we have
\[ (\mathcal{O}/\mathfrak{f})^\times \cong (\mathbb{Z}/f\mathbb{Z})^\times, \quad |(\mathcal{O}/\mathfrak{f})^\times|=|(\mathbb{Z}/f\mathbb{Z})^\times|=\prod_{p|f} (p-1).\]
Hence $|(\mathcal{O}_k/\mathfrak{f})^\times|/|(\mathcal{O}/\mathfrak{f})^\times|=\prod_{p|f} p=f$. 
By \eqref{eq:ORDERCG} and \eqref{eq:ORDERCG-PLUS}, we have two exact sequence of finite abelian groups
\begin{align*}
 &1\longrightarrow \ker(\Cl_\mathcal{O}\rightarrow \Cl_k)\longrightarrow \Cl_\mathcal{O}\longrightarrow \Cl_k \longrightarrow 1, \\
 &1\longrightarrow \ker(\Cl_{\mathcal{O},+}\rightarrow \Cl_{k,+})\longrightarrow \Cl_{\mathcal{O},+} \longrightarrow \Cl_{k,+} \longrightarrow 1.
\end{align*}
By Corollary \ref{cor:DEDEKIND} and Corollary \ref{cor:CARDORDER-PLUS}, 
the orders of the kernels are $f/(U(\mathcal{O}_k):U(\mathcal{O}))$ and 
$f/(U_+(\mathcal{O}_k):U_+(\mathcal{O}))$,  which are odd integers. 
Hence $\Cl_\mathcal{O}/\Cl^2_\mathcal{O} \cong \Cl_k/\Cl^2_k$ and 
$\Cl_{\mathcal{O},+}/\Cl^2_{\mathcal{O},+} \cong \Cl_{k,+}/\Cl^2_{k,+}$. 
Thus we have $|\Cl^{(2)}_\mathcal{O}|=|\Cl^{(2)}_k|$ and $|\Cl^{(2)}_{\mathcal{O},+}|=|\Cl^{(2)}_{k,+}|$. 
By the same argument, we have $|\Cl^{(2)}_{\mathcal{O}_0}|=|\Cl^{(2)}_k|$ and $|\Cl^{(2)}_{\mathcal{O}_0,+}|=|\Cl^{(2)}_{k,+}|$. 
Since $X(\mathcal{O},\mathcal{O}_0)$ (resp. $X_+(\mathcal{O},\mathcal{O}_0)$) is a subgroup 
of $\Cl_{\mathcal{O}_0}/\Cl_{\mathcal{O}_0}^2$ (resp. $\Cl_{\mathcal{O}_0,+}/\Cl_{\mathcal{O}_0,+}^2$) 
generated by the ideal class of $\mathfrak{j}(\mathcal{O},\mathcal{O}_0)$ and 
since $\mathfrak{j}(\mathcal{O},\mathcal{O}_0)\mathcal{O}_k=\mathfrak{h}=\prod_{p|h} \mathfrak{p}^2$, 
$X(\mathcal{O},\mathcal{O}_0)$ (resp. $X_+(\mathcal{O},\mathcal{O}_0)$) is trivial. 
\qed

\section{Quartic rings}
We summarize the results in Bhargava \cite{BH3}. 
Let $A(v)$ be an integral ternary quadratic form in variables $v=(v_1,v_2,v_3)$.
We write
\[ A(v)=a_{11}v_1^2+a_{12}v_1v_2+a_{13}v_1v_3+a_{22}v_2^2+a_{23}v_2v_3+a_{33}v_3^2,\]
with $a_{ij} \in \mathbb{Z}$. We identify $A(v)$ by the symmetric matrix 
\[ A=\left(
\begin{array}{ccc}
a_{11} & \frac{1}{2}a_{12}& \frac{1}{2}a_{13} \\[5pt]
\frac{1}{2}a_{12} & a_{22} & \frac{1}{2}a_{23} \\[5pt]
\frac{1}{2}a_{13} & \frac{1}{2}a_{23}& a_{33}
\end{array}
\right).
\]
Let $(A,B)$ be a pair of integral ternary quadratic forms and write
\[ A(v)=\sum_{1\leq i\leq j\leq 3} a_{ij}v_iv_j, \quad B(v)=\sum_{1\leq i\leq j\leq 3} b_{ij}v_iv_j\]
with $a_{ij}, b_{ij} \in \mathbb{Z}$. 
Letting $a_{ji} = a_{ij}$ and $b_{ji} = b_{ij}$, we define constants 
$\lambda_{k\ell}^{ij}=\lambda_{k\ell}^{ij}(A,B)$ by
\begin{equation}\label{eq:LIJKL}
\lambda_{k\ell}^{ij}(A,B)=\left|
\begin{array}{cc}
a_{ij} & b_{ij} \\
a_{k\ell} & b_{k\ell}
\end{array}
\right|.
\end{equation}
For any permutation $(i, j, k)$ of $(1, 2, 3)$, 
we define constants $c_{ij}^k=c_{ij}^k(A,B)$ by
\begin{equation}\label{eq:CIJK}
\begin{split}
c_{ii}^i &=\pm \lambda_{ij}^{ik}+C_i, \\
c_{ii}^j &=\pm \lambda_{ik}^{ii}, \\
c_{ij}^i &=\pm \frac{1}{2}\lambda_{jj}^{ik}+\frac{1}{2}C_j, \\
c_{ij}^k &=\pm \lambda_{ii}^{jj},
\end{split}
\end{equation}
where we have used $\pm$ to denote the sign of the permutation $(i, j, k)$ of $(1, 2, 3)$,
and where the constants $C_i$ are given by
\begin{equation}\label{eq:CI}
C_1=\lambda_{11}^{23}, \quad C_2=-\lambda_{22}^{13}, \quad C_3=\lambda_{33}^{12}.
\end{equation}
Then we have 
\[ c_{12}^1=c_{12}^2=c_{13}^1=0,  \quad c_{13}^3=\lambda_{11}^{23}, \quad c_{23}^2=\lambda_{33}^{12},
\quad c_{23}^3=-\lambda_{22}^{13}.
\]
Hence the values of the $c_{ij}^k$ (for $k > 0$) are all integral.
We put
\begin{equation}\label{eq:CIJ0}
c_{ij}^0=\sum_{r=1}^3 (c_{jk}^rc_{ri}^k - c_{ij}^r c_{rk}^k)
\end{equation}
for any $k \in \{1, 2, 3\}\smallsetminus \{i\}$. 
It follows from \eqref{eq:CIJK} that the expression of $c_{ij}^0$ 
does not depend on the choice of $k$. 
We now define a quartic ring $Q(A,B)$ as follows. 
Let $Q(A,B)$ be a free $\mathbb{Z}$-module with basis $\{\alpha_0=1,\alpha_1,\alpha_2,\alpha_3\}$ 
and the multiplication of $Q(A,B)$ is given by 
\begin{equation}\label{eq:QUARTIC-RING-STR}
\alpha_i\alpha_j=\sum_{k=0}^3 c_{ij}^k \alpha_k \quad (i,j \in \{1,2,3\}).
\end{equation}
Then $Q(A,B)$ becomes a quartic ring. 
We call a basis $\{1,\alpha_1,\alpha_2,\alpha_3\}$ \textit{normalized} if $c_{12}^1=c_{12}^2=c_{13}^1=0$. 
We defined in \S 1 an integral binary cubic form $F_{(A,B)}(u)$ in $u=(u_1,u_2)$, 
which equals to $4\det (u_1A-u_2B)$. 
The discriminant $\Disc(Q(A,B))$ of $Q(A,B)$ is equal to $\Disc(A,B)=\Disc(F_{(A,B)})$.
Let $R(A,B)$ be the cubic ring corresponding to the binary cubic form $F_{(A,B)}(u)$.
Then the discriminant $\Disc(R(A,B))$ of $R(A,B)$ is also equal to $\Disc(A,B)=\Disc(F_{(A,B)})$.

For any quartic ring $Q$, we denote by $\bar{Q}$ the $S_4$-closure of $Q$.
For an element $x \in Q$, we denote by $x,x',x'',x'''$ the conjugates of $x$ in $\bar{Q}$. 
We put 
\[ \tilde{\phi}_{4,3}(x)=xx'+x''x'''. \]
Then all $\tilde{\phi}_{4,3}(x)$ are contained in the same cubic ring, namely 
the cubic subring of $\bar{Q}$ consisting of all elements which are fixed by the action of 
$D_4=\langle (1\,2),(1\,3\,2\,4) \rangle \subset S_4$. 
We define the invariant cubic ring $\CInv(Q)$ by
\[ \CInv(Q)=\mathbb{Z}\left[\left\{\tilde{\phi}_{4,3}(x)\,|\,x \in Q \right\} \right].\]
A \textit{cubic resolvent ring} of $Q$ is a cubic ring $R$ such that $\Disc(R)=\Disc(Q)$ and 
$R \supset \CInv(Q)$. 

Bhargava proved the following theorem (\cite[Theorem 1]{BH3}).
\begin{thm}\label{thm:BHARGAVA3}
The correspondence $(A,B) \mapsto (Q(A,B),R(A,B))$ induces a canonical bijection 
between the set of $\Gamma$-orbits of nondegenerate pairs of integral ternary quadratic forms 
and the set of isomorphism classes of pairs $(Q,R)$, 
where $Q$ is a nondegenerate quartic ring and $R$ is a cubic resolvent ring of $Q$.
\end{thm}

It is clear from formulae \eqref{eq:CIJK}, \eqref{eq:CI}, \eqref{eq:CIJ0}, and \eqref{eq:QUARTIC-RING-STR} 
that the content of a quartic ring $Q(A,B)$ is equal to the greatest common divisor of the fifteen 
$\SL_2$-invariants $\lambda^{ij}_{k\ell}(A,B)$. Thus the \textit{content} $\ct(A,B)$ of 
a pair $(A,B)$ of integral ternary quadratic forms is defined by 
\[ \ct(A,B) = \ct(Q(A,B)) = \gcd \{\lambda^{ij}_{k\ell}(A,B)\}.\]
We call $(A,B)$ \textit{primitive} if $\ct(A,B)=1$. 
The following corollary is \cite[Corollary 5]{BH3}
\begin{cor}\label{cor:UNIQUE-RESOLVENT}
Every quartic ring has a cubic resolvent ring. 
A primitive quartic ring has a unique cubic resolvent ring up to isomorphism. 
In particular, every maximal quartic ring has a unique cubic resolvent ring. 
\end{cor}

Let $(A,B)$ be a nondegenerate pair of integral ternary quadratic forms and 
$\gamma=(\gamma_1,\gamma_2)$ an element of $\Gamma$.  We put $(A',B')=\gamma\cdot(A,B)$. 
Then there exists a ring isomorphism $\psi:Q(A',B')\rightarrow Q(A,B)$ 
by Theorem \ref{thm:BHARGAVA3}.  
We now give the isomorphism $\psi$ explicitly. 

\begin{prop}\label{prop:QRINGHOMO}
Let $(A,B)$ be a nondegenerate pair of integral ternary quadratic forms and 
$\delta=(\delta_1,\delta_2)$ be an element of $M_3(\mathbb{Z})\times M_2(\mathbb{Z})$ 
such that $\det \delta_1\neq0$, $\det \delta_2 \neq 0$. We put $(A',B')=\delta\cdot(A,B)$ 
and denote by $\{\beta_i\}$ the basis $\{\alpha_i\}$ for the quartic ring $Q(A',B')$. 
Then there exists an injective ring homomorphism 
$\psi:Q(A',B')\rightarrow Q(A,B)$ satisfying
\begin{equation}\label{eq:TRANSFORM}
\left(
\begin{array}{c}
\psi(\beta_1)  \\
\psi(\beta_2)  \\
\psi(\beta_3) 
\end{array}
\right)
\equiv (\det \delta_1) (\det \delta_2) \delta_1
\left(
\begin{array}{c}
\alpha_1  \\
\alpha_2  \\
\alpha_3 
\end{array}
\right)
\pmod{\mathbb{Z}}
\end{equation} 
\end{prop}
\Proof
We write $A'(v)=\sum_{i\leq j} a'_{ij}v_iv_j$ and $B'(v)=\sum_{i\leq j} b'_{ij}v_iv_j$. 
We set
\begin{equation}\label{eq:DEFDIJPRIME}
\lambda^{ij}_{k\ell}=\left|
\begin{array}{cc}
a_{ij} &  b_{ij}\\
a_{k\ell} & b_{k\ell}
\end{array}
\right|,\qquad
\tilde{\lambda}^{ij}_{k\ell}=\left|
\begin{array}{cc}
a'_{ij} &  b'_{ij}\\
a'_{k\ell} &' b_{k\ell}
\end{array}
\right|.
\end{equation}
We also set $c_{ij}^k=c_{ij}^k(A,B)$ and $\tilde{c}_{ij}^k=c_{ij}^k(A',B')$.

Case 1.  $\delta_1=1_3$, $\delta_2=\left(
\begin{array}{cc}
p & q \\
r & s
\end{array}
\right)$. By definition, we have
\[\tilde{\lambda}^{ij}_{k\ell}
=\left|
\begin{array}{cc}
a_{ij} & b_{ij} \\
a_{k\ell} & b_{k\ell}
\end{array}
\right|
\left|
\begin{array}{cc}
p & r \\
q & s
\end{array}
\right|=(\det \delta_2)\lambda^{ij}_{k\ell}.
\]
Hence $\tilde{c}_{ij}^k=(\det \delta_2)c_{ij}^k$ for $k>0$. 
So we have an injective ring homomorphism $\psi:Q(A',B')\rightarrow  Q(A,B)$ such that 
$\psi(\beta_i)=(\det \delta_2)\alpha_i$, $i=1,2,3$. 

Case 2. $\delta_1=\diag[d_1,d_2,d_3]$, $d_i\neq 0$, $\delta_2=1_2$. 
Then we have $a'_{ij}=d_id_j a_{ij}, b'_{ij}=d_id_j b_{ij}$. 
Hence we have 
\[ \tilde{\lambda}^{ij}_{k\ell}=
\left|
\begin{array}{cc}
d_id_ja_{ij} & d_id_jb_{ij} \\
d_kd_\ell a_{k\ell} & d_kd_\ell b_{k\ell}
\end{array}
\right|
=d_id_jd_kd_\ell \lambda^{ij}_{k\ell}.
\]
Hence $d_k \tilde{c}_{ij}^k=(\det \delta_1)d_id_jc_{ij}^k$ for $k>0$. 
So we have an injective ring homomorphism $\psi:Q(A',B') \rightarrow Q(A,B)$ 
such that $\psi(\beta_i)=(\det \delta_1)d_i \alpha_i$ 
for $i=1,2,3$.

Case 3. $\delta_2=1_2$ and $\delta_1$ is one of the following three matrices:
\[ \sigma_{12}=\left(
\begin{array}{ccc}
0 & 1 & 0 \\
1 & 0 & 0 \\
0 & 0 & 1
\end{array}
\right), \quad 
\sigma_{23}=\left(
\begin{array}{ccc}
1 & 0 & 0 \\
0 & 0 & 1 \\
0 & 1 & 0
\end{array}
\right), \quad
\tau_{21}=\left(
\begin{array}{ccc}
1 & 0 & 0 \\
1 & 1 & 0 \\
0 & 0 & 1
\end{array}
\right).
\]

Case 3-1. $\delta_1=\sigma_{12}$. 
We put $1'=2$, $2'=1$, $3'=3$. Then we have $a'_{ij}=a_{i'j'}$, $b'_{ij}=b_{i'j'}$.
Hence $\tilde{\lambda}^{ij}_{k\ell}=\lambda^{i'j'}_{k'\ell'}$. 
For any permutation $(i\,j\,k)$ of $(1\,2\,3)$, $(i'\,j'\,k')$ is also a permutation 
of $(1\,2\,3)$ whose sign is $-1$ times that of $(i\,j\,k)$. 
So if we use $\pm$ to denote the sign of $(i\,j\,k)$, then that of $(i'\,j'\,k')$ is $\mp$. 
Hence we have
\[ \tilde{c}_{ii}^j=\pm \tilde{\lambda}^{ii}_{ik}=\pm \lambda^{i'i'}_{i'k'}=-c_{i'i'}^{j'}.\]
Similarly we have $\tilde{c}_{ij}^k=-c_{i'j'}^{k'}$, 
$\tilde{c}_{12}^1=c_{12}^1=\tilde{c}_{12}^2=c_{12}^2=\tilde{c}_{13}^1=c_{13}^1=0$,
\begin{alignat*}{3}
 \tilde{c}_{11}^1 &=-c_{22}^2,
   & \quad  \tilde{c}_{22}^2 &=-c_{11}^1,
   & \quad  \tilde{c}_{33}^3 &=-c_{33}^3+2c_{23}^2,  \\
 \tilde{c}_{23}^2 &=c_{23}^2,
   & \quad  \tilde{c}_{13}^3 &=-c_{23}^3, 
   & \quad  \tilde{c}_{23}^3 &=-c_{13}^3. 
\end{alignat*}
We define an isomorphism $\psi:Q(A',B')\rightarrow Q(A,B)$ of $\mathbb{Z}$-modules by
\[ \psi(\beta_0)=\alpha_0, \quad \psi(\beta_1)=-\alpha_2, \quad \psi(\beta_2)=-\alpha_1,\quad  
\psi(\beta_3)=c_{23}^2-\alpha_3.
\]
By the relations above, we have the congruences
\begin{equation}\label{eq:PSI-CONG-MODZ}
\psi(\beta_i\beta_j)=\sum_{k=0}^3 \tilde{c}_{ij}^k\psi(\beta_k) 
\equiv \psi(\beta_i)\psi(\beta_j) \pmod{\mathbb{Z}}
\end{equation}
for all $1\leq i\leq j \leq 3$. The congruences imply that $\{\psi(\beta_i)\}$ is a normalized basis of 
the quartic ring $Q(A,B)$. So the constants $\tilde{c}_{ij}^0$'s are 
determined by $\tilde{c}_{ij}^k$'s for $k>0$. 
Hence the congruences become the equations $\psi(\beta_i\beta_j)=\psi(\beta_i)\psi(\beta_j)$. 
Thus $\psi$ is a ring isomorphism.

Case 3-2. $\delta_1=\sigma_{23}$. 
We put $1'=1$, $2'=3$, $3'=2$. By the same argument as in Case 3-1, 
we have
\[ \tilde{c}_{ii}^j=-c_{i'i'}^{j'}.\]
Similarly we have $\tilde{c}_{ij}^k=-c_{i'j'}^{k'}$, 
$\tilde{c}_{12}^1= c_{12}^1=\tilde{c}_{12}^2=c_{12}^2=\tilde{c}_{13}^1=c_{13}^1=0$,
\begin{alignat*}{3}
 \tilde{c}_{11}^1 &=-c_{11}^1+2c_{13}^3,
   & \quad  \tilde{c}_{22}^2 &=-c_{33}^3,
   & \quad  \tilde{c}_{33}^3 &=-c_{22}^2,  \\
 \tilde{c}_{23}^2 &=-c_{23}^3,
   & \quad  \tilde{c}_{13}^3 &=c_{13}^3, 
   & \quad  \tilde{c}_{23}^3 &=-c_{23}^2. 
\end{alignat*}
We define an isomorphism $\psi:Q(A',B')\rightarrow Q(A,B)$ of $\mathbb{Z}$-modules by
\[ \psi(\beta_0)=\alpha_0, \quad \psi(\beta_1)=c_{13}^3-\alpha_1, \quad \psi(\beta_2)=-\alpha_3,\quad  
\psi(\beta_3)=-\alpha_2.
\]
By the relations above, we have the congruences \eqref{eq:PSI-CONG-MODZ} 
which become the equations $\psi(\beta_i\beta_j)=\psi(\beta_i)\psi(\beta_j)$
as in Case 3-1. Thus $\psi$ is a ring isomorphism.

Case 3-3. $\delta_1=\tau_{21}$. In this case, we have 
\begin{alignat*}{3}
 a'_{11}&=a_{11},  & \quad  a'_{12}&=a_{12}+2a_{11},  & \quad  a'_{13}&=a_{13}, \\
 a'_{22}&=a_{22}+a_{11}+a_{12},  & \quad  a'_{23}&=a_{23}+a_{13},  & \quad  a'_{33}&=a_{33},  \\
 b'_{11}&=b_{11},  & \quad  b'_{12}&=b_{12}+2b_{11},  & \quad  b'_{13}&=b_{13}, \\  
 b'_{22}&=b_{22}+b_{11}+b_{12},  & \quad  b'_{23}&=b_{23}+b_{13},  & \quad  b'_{33}&=b_{33}.
\end{alignat*}
Hence we have $\tilde{c}_{12}^1=c_{12}^1=\tilde{c}_{12}^2=c_{12}^2=\tilde{c}_{13}^1=c_{13}^1=0$, 
\begin{align*}
\tilde{c}_{11}^1&=c_{11}^1-3c_{11}^2, \quad \tilde{c}_{11}^2=c_{11}^2, \quad \tilde{c}_{11}^3=c_{11}^3, \quad \tilde{c}_{12}^3=c_{12}^3+c_{11}^3,\quad \tilde{c}_{13}^2=c_{13}^2, \\
\tilde{c}_{13}^3&=c_{13}^3-c_{11}^2,\quad \tilde{c}_{22}^1=c_{11}^1-c_{11}^2+c_{22}^1-c_{22}^2,\quad \tilde{c}_{22}^2=c_{22}^2-2c_{11}^1+3c_{11}^2,\\
\tilde{c}_{22}^3&=c_{22}^3+2c_{12}^3+c_{11}^3, \quad \tilde{c}_{23}^1=c_{23}^1-c_{23}^2-c_{13}^2,\quad \tilde{c}_{23}^2=c_{23}^2+2c_{13}^2, \\
\tilde{c}_{23}^3&=c_{23}^3+c_{11}^2-c_{11}^1+c_{13}^3,\quad \tilde{c}_{33}^1=c_{33}^1-c_{33}^2, \\
\tilde{c}_{33}^2&=c_{33}^2,\quad \tilde{c}_{33}^3=c_{33}^3+2c_{13}^2.
\end{align*}
We define an isomorphism $\psi:Q(A',B')\rightarrow Q(A,B)$ of $\mathbb{Z}$-modules by 
$\psi(\beta_0)=\alpha_0$, $\psi(\beta_1)=\alpha_1-c_{11}^2$, $\psi(\beta_2)=\alpha_1+\alpha_2+c_{11}^2-c_{11}^1$, 
$\psi(\beta_3)=\alpha_3+c_{13}^2$.
By the relations above, we have the congruences \eqref{eq:PSI-CONG-MODZ} 
which become the equations $\psi(\beta_i\beta_j)=\psi(\beta_i)\psi(\beta_j)$
as in Case 3-1. Thus $\psi$ is a ring isomorphism.

In all cases, we have an injective ring homomorphism 
$\psi:Q(A',B')\rightarrow Q(A,B)$ which satisfies the formula (\ref{eq:TRANSFORM}). 
Since the three matrices $\sigma_{12}$, $\sigma_{23}$ and $\tau_{21}$ generate the group $\GL_3(\mathbb{Z})$
and any $\delta_1 \in M_3(\mathbb{Z})$ with $\det \delta_1 \neq 0$
can be written as $\delta_1=\gamma_1\diag[d_1,d_2,d_3]\gamma'_1$ with $\gamma_1, \gamma'_1 \in \GL_3(\mathbb{Z})$, 
the formula (\ref{eq:TRANSFORM}) holds for any $\delta=(\delta_1,\delta_2) \in M_3(\mathbb{Z})\times M_2(\mathbb{Z})$ 
with $\det \delta_1 \neq 0$, $\det \delta_2 \neq 0$. 
In particular, for any $\delta \in \Gamma$, the ring isomorphism $\psi:Q(A',B') \cong Q(A,B)$
satisfies (\ref{eq:TRANSFORM}).
\qed

Let $(A,B)$ be a nondegenerate pair of integral ternary quadratic forms. 
We now consider the relation between the isotropy group $\Gamma_{(A,B)}$ of $(A,B)$ in $\Gamma$ 
and the automorphism group $\Aut(Q(A,B))$ of the quartic ring $Q(A,B)$. 

\begin{lem}\label{lem:A11A12-ZERO}
For any nondegenerate pair $(A,B)$ of integral ternary quadratic forms, 
there exists an element $\gamma \in \Gamma$ such that $(A',B')=\gamma \cdot (A,B)$ satisfies 
$a'_{11}=a'_{12}=0$, $b'_{11}\neq 0$, where $a'_{ij}$ and $b'_{ij}$ are the coefficients 
of $v_iv_j$ in $A'(v)$ and $B'(v)$, respectively.
\end{lem}
\Proof
Since $B\neq 0$, there exist integers $p,q,r$ such that $B(p,q,r)\neq 0$ and $\gcd(p,q,r)=1$.
We take a matrix $\gamma_1\in \GL_3(\mathbb{Z})$ whose first row is $(p,q,r)$ and put
$(A',B')=(\gamma_1,1_2)\cdot (A,B)$. Then $B'(1,0,0)=B(p,q,r)\neq 0$. 
So we may assume $b_{11}\neq 0$ from the beginning. 
We put $t=\gcd(a_{11},b_{11})$ and write $a_{11}=ta_1$, $b_{11}=tb_1$, $\gcd(a_1,b_1)=1$. 
We take integers $p$, $q$ such that $a_1p+b_1q=1$ and put 
\[\gamma_2=\left(
\begin{array}{cc}
b_1 & -a_1 \\
p & q
\end{array}
\right) \in \SL_2(\mathbb{Z}), \quad (A',B')=(1_3,\gamma_2)\cdot (A,B).
\]
Then we have $a'_{11}=0$ and $b'_{11}=t\neq 0$. 
Thus we may assume $a_{11}=0$ and $b_{11}\neq 0$ from the beginning. 
If $a_{12}=0$, then nothing remains to prove.
We assume $a_{12}\neq 0$. We put $t=\gcd(a_{12},a_{13})$ and write 
$a_{12}=ta_2$, $a_{13}=ta_3$, $\gcd(a_2,a_3)=1$. 
We take integers $p$, $q$ such that $a_2p+a_3q=1$ and put 
\[ \gamma_1=\left(
\begin{array}{ccc}
1 & 0 & 0 \\
0 & a_3 & -a_2 \\
0 & p& q
\end{array}
\right) \in \SL_3(\mathbb{Z}), \quad (A',B')=(\gamma_1,1_2)\cdot (A,B).
\]
Then we have $a'_{11}=a'_{12}=0$, $b'_{11}=b_{11}\neq 0$. 
\qed

\begin{prop}\label{prop:ISOTROPY-AUT}
For any nondegenerate pair $(A,B)$ of integral ternary quadratic forms, 
there exists an injective group homomorphism 
$\Gamma_{(A,B)} \rightarrow \Aut(Q(A,B))$. 
Further if $(A,B)$ is primitive, then the homomorphism is an isomorphism.
\end{prop}
\Proof
We denote by $\{\alpha_i\}$ the basis of $Q(A,B)$ as before. 
By Proposition \ref{prop:QRINGHOMO}, any $\gamma=(\gamma_1,\gamma_2) \in \Gamma_{(A,B)}$ 
gives a ring automorphism $\psi:Q(A,B)\rightarrow Q(A,B)$ such that
\[\left(
\begin{array}{c}
\psi(\alpha_1)  \\
\psi(\alpha_2)  \\
\psi(\alpha_3) 
\end{array}
\right)
=(\det \gamma_2)\gamma_1 \left(
\begin{array}{c}
\alpha_1  \\
\alpha_2  \\
\alpha_3 
\end{array}
\right)+\left(
\begin{array}{c}
e_1  \\
e_2  \\
e_3 
\end{array}
\right)
\]
for some integers $e_1,e_2,e_3$. 
Since $\{\alpha_i\}$ is a normalized basis, i. e. $c_{12}^1=c_{12}^2=c_{13}^1=0$, 
the integers $e_1,e_2,e_3$ are determined by $\gamma$. Hence the correspondence 
$\gamma\mapsto \psi^{-1}$ defines a natural group homomorphism 
$\Gamma_{(A,B)} \rightarrow \Aut(Q(A,B))$. 
If $\gamma=(\gamma_1,\gamma_2) \in \Gamma_{(A,B)}$ is in the kernel, 
we have $(\det \gamma_2) \gamma_1=1_3$, hence 
$\det \gamma_2=1$ and $\gamma_1=1_3$. 
Since $A$ and $B$ are linearly independent over $\mathbb{C}$, we have $\gamma_2=1_2$. 
We now assume that $(A,B)$ is primitive. 
We replace $(A,B)$ by a suitable $\Gamma$-equivalent pair by Lemma \ref{lem:A11A12-ZERO}, 
we may assume that $(A,B)$ is of the form 
\begin{align*}
A(v) &=a_{13}v_1v_3+a_{22}v_2^2+a_{23}v_2v_3+a_{33}v_3^2, \\
B(v) &=b_{11}v_1^2+b_{12}v_1v_2+b_{13}v_1v_3+b_{22}v_2^2+b_{23}v_2v_3+b_{33}v_3^2, \quad b_{11} \neq 0.
\end{align*}
Since $(A,B)$ is primitive, $\gcd(a_{13},a_{22},a_{23},a_{33})=1$. 
Let $\sigma$ be an element of $\Aut(Q(A,B))$. 
Then there exists an element $\gamma_1 \in \GL_3(\mathbb{Z})$ such that 
${}^t(\sigma(\alpha_i))\equiv \gamma_1 {}^t(\alpha_i) \pmod{\mathbb{Z}}$. 
We take an element $\gamma_2 \in \Gamma_2$ such that $\det \gamma_2=\det \gamma_1$ 
and put $(A',B')=(\gamma_1,\gamma_2)\cdot(A,B)$. 
We denote by $a'_{ij}$ and $b'_{ij}$ the coefficients of $v_iv_j$ in 
$A'(v)$ and $B'(v)$, respectively.
We can take $\gamma_2$ so that $a'_{11}=0$. 
We denote by $\{\alpha'_i\}$ the basis $\{\alpha_i\}$ for $Q(A',B')$. 
By Proposition \ref{prop:QRINGHOMO}, there exists a ring isomorphism 
$\psi:Q(A',B')\rightarrow Q(A,B)$ such that 
${}^t(\psi(\alpha'_i)) \equiv \gamma_1{}^t(\alpha_i) \pmod{\mathbb{Z}}$. 
Hence $\sigma(\alpha_i)\equiv \psi(\alpha'_i) \pmod{\mathbb{Z}}$. 
Since $\{\psi(\alpha'_i)\}$ and $\{\sigma(\alpha_i)\}$ are both normalized basis of 
the quartic ring $Q(A,B)$, $\sigma(\alpha_i)=\psi(\alpha'_i)$. 
So we have $c_{ij}^k(A',B')=c_{ij}^k(A,B)$, hence $\lambda^{ij}_{k\ell}(A',B')=\lambda^{ij}_{k\ell}(A,B)$ for all $i,j,k,\ell$. 
Since $a_{11}=a'_{11}=0$, we have $b_{11}a_{k\ell}=b'_{11}a'_{k\ell}$ for all 
$1\leq k \leq \ell \leq 3$. 
If $b'_{11}=0$, then $b_{11} \neq 0$ implies $a_{k\ell}=0$. This is a contradiction. 
So we have $b'_{11} \neq 0$. We put $t=b_{11}/b'_{11}\in \mathbb{Q}^\times$. Then we have 
$a'_{k\ell}=ta_{k\ell}$ for all $k,\ell$. In particular, we have $a'_{12}=ta_{12}=0$. 
Since $(A,B)$ is primitive, so is $(A',B')$. 
Since $\gcd(a_{13},a_{22},a_{23},a_{33})=1$, the equations $a'_{k\ell}=ta_{k\ell}$ for all $k,\ell$ 
imply $t \in \mathbb{Z}$. Similarly we have $t^{-1} \in \mathbb{Z}$. Hence $t=\pm 1$. 
We take $(i,j)$ such that $a_{ij}\neq 0$. 
Then it follows from the equations $\lambda^{ij}_{k\ell}(A,B)=\lambda^{ij}_{k\ell}(A',B')$ that 
$a_{ij}(b_{k\ell}-tb'_{k\ell})=(b_{ij}-tb'_{ij})a_{k\ell}$. 
If we put $s=(b_{ij}-tb'_{ij})/a_{ij} \in \mathbb{Q}$, then we have 
$b'_{k\ell}=-t^{-1}sa_{k\ell}+t^{-1}b_{k\ell}$. 
Since $t=\pm 1$ and $\gcd(a_{13},a_{22},a_{23},a_{33})=1$, we have $s \in \mathbb{Z}$.
We put 
\[\delta_2=\left(
\begin{array}{cc}
t & 0 \\
-t^{-1}s & t^{-1}
\end{array}
\right) \in \SL_2(\mathbb{Z}). 
\]
Then we have $(A',B')=(1_3,\delta_2)\cdot (A,B)$. 
We put $\epsilon=\det \gamma_1$, $\gamma'_1=\epsilon \gamma_1$ and $\gamma'_2=\delta_2^{-1}\gamma_2$. 
Then $(\gamma'_1,\gamma'_2)\cdot (A,B)=(A,B)$ and $(\det \gamma'_2)\gamma'_1=\epsilon \gamma'_1=\gamma_1$. 
Hence $(\gamma'_1,\gamma'_2)$ is an element of $\Gamma_{(A,B)}$ which gives the ring automorphism $\sigma$. 
This proves the surjectivity.
\qed

Let $(A,B)$ be a nondegenerate pair of integral ternary quadratic forms. 
We put $K=Q(A,B) \otimes_\mathbb{Z} \mathbb{Q}$, $k=R(A,B)\otimes_\mathbb{Z} \mathbb{Q}$ 
and $F=F_{(A,B)}$.  If the index $(\mathcal{O}_K:Q(A,B))$ is cube free, 
then $Q(A,B)$ is primitive, hence by Proposition \ref{prop:ISOTROPY-AUT}, 
$\Gamma_{(A,B)} \cong \Aut(Q(A,B))$. 
So we have $\mu(A,B)=1/|\Aut(Q(A,B))|$ for any primitive pair $(A,B)$. 

If $K$ is a quartic field, we denote by $\tilde{K}$ the Galois closure of $K$ over $\mathbb{Q}$ 
and put $G=\Gal(\tilde{K}/\mathbb{Q})$. 
We denote by $S_n$, $A_n$, $D_n$, $C_n$ and $V_4$, the symmetric group of degree $n$, 
the alternating group of degree $n$, the dihedral group of order $2n$, 
the cyclic group of order $n$ and the Klein four-group, respectively. 
If $G=S_4$ or $A_4$, then $K$ has no nontrivial subfield. Hence $\Aut(\mathcal{O}_K)$ is trivial 
and so is $\Aut(Q(A,B))$. In this case, 
$k$ is a non-Galois cubic field or Galois cubic field according as $G=S_4$ or $A_4$. 
If $G=D_4$, then $K$ has a unique quadratic subfield $k_2$. 
Hence $\Aut(\mathcal{O}_K)=C_2$ and $\Aut(Q(A,B))$ is trivial or equal to $C_2$. 
In this case, we can write $K=k_2(\sqrt{\theta})$ and $k_2=\mathbb{Q}(\theta)$ 
for some $\theta \in k_2$. Then $k=\mathbb{Q}\oplus k_1$ and 
$k_1=\mathbb{Q}\big(\sqrt{N_{k_2/\mathbb{Q}}(\theta)}\big)$ is a quadratic field. 
If $G=C_4$, then $\Aut(\mathcal{O}_K)=C_4$, hence 
$\Aut(Q(A,B))$ is trivial, $C_2$ or $C_4$. In this case, 
$k=\mathbb{Q}\oplus k_1$, where $k_1$ is the unique quadratic subfields of $K$. 
If $G=V_4$, then $\Aut(\mathcal{O}_K)=V_4$, hence 
$\Aut(Q(A,B))$ is trivial, $C_2$ or $V_4$. In this case, $k=\mathbb{Q}^3$.

It remains to deal with the case when $K$ is not a quartic field. 
If $K=\mathbb{Q} \oplus k_3$ where $k_3$ is a non-Galois cubic field, then 
$\Aut(\mathcal{O}_K)$ is trivial, hence so is $\Aut(Q(A,B))$. In this case, 
$k=k_3$ is a non-Galois cubic field. 
If $K=\mathbb{Q} \oplus k_3$ where $k_3$ is a Galois cubic field, then 
$\Aut(\mathcal{O}_K)=C_3$, hence $\Aut(Q(A,B))$ is trivial or $C_3$. 
In this case, $k=k_3$ is a Galois cubic field. 
If $K=k_2 \oplus k_3$ where $k_2$, $k_3$ are distinct quadratic fields, 
then $\Aut(\mathcal{O}_K)=V_4$, hence $\Aut(Q(A,B))$ is trivial, 
$C_2$ or $V_4$. In this case, $k=\mathbb{Q}\oplus k_1$, $k_1=\mathbb{Q}(\sqrt{\Disc(k_2)\Disc(k_3)})$. 
If $K=\mathbb{Q}^2 \oplus k_1$ where $k_1$ is a quadratic field, 
then $\Aut(Q(A,B))$ is the same as in the previous case and $k=\mathbb{Q}\oplus k_1$. 
If $K=k_1 \oplus k_1$ where $k_1$ is a quadratic field, 
then $\Aut(\mathcal{O}_K)=D_4$, hence $\Aut(Q(A,B))$ is trivial, 
$C_2$, $V_4$ or $D_4$. In this case, $k=\mathbb{Q}^3$.
If $K=\mathbb{Q}^4$, then $\Aut(\mathcal{O}_K)=S_4$, hence 
$\Aut(Q(A,B)))$ is trivial, $C_2$, $C_4$, $V_4$, $D_4$, $A_4$ or $S_4$. 
In this case, $k=\mathbb{Q}^3$.

So the cubic algebra $k$ is a field if and only if 
$K$ is a quartic field with $G=S_4$, $A_4$ or 
$K$ is a direct sum of $\mathbb{Q}$ and a cubic field. 
In these cases, $\Aut(Q(A,B))$ is trivial or $C_3$. 

\begin{prop}\label{prop:QUARTIC-SUBRING}
Let $(A,B)$, $(A',B')$ be primitive nondegenerate pairs of integral ternary quadratic forms 
such that $Q(A',B')\subset Q(A,B)$. Then there exist matrices 
$\delta_1 \in M_3(\mathbb{Z})$ and $\delta_2 \in M_2(\mathbb{Z})$ which satisfy
$(A',B')=(\delta_1,\delta_2^{-1})\cdot (A,B)$, 
$|\det \delta_1|=|\det \delta_2|=(Q(A,B):Q(A',B'))$.
\end{prop}
\Proof
We put $m=(Q(A,B):Q(A',B'))$ and denote by $\{\alpha'_i\}$ the basis $\{\alpha_i\}$ 
for $Q(A',B')$. There exists a matrix $\delta_1 \in M_3(\mathbb{Z})$ such that $|\det \delta_1|=m$ 
and ${}^t(\alpha'_1,\alpha'_2,\alpha'_3) \equiv \delta_1 {}^t(\alpha_1,\alpha_2,\alpha_3) \pmod{\mathbb{Z}}$. 
We put $(A'',B'')=(\delta_1,1_2)\cdot (A,B)$ and denote by $\{\alpha''_i\}$
the basis $\{\alpha_i\}$ for $Q(A'',B'')$. 
By Proposition \ref{prop:QRINGHOMO}, there exists an injective ring homomorphism 
$\psi:Q(A'',B'')\rightarrow Q(A,B)$ such that
\[ {}^t(\alpha''_i)\equiv (\det \delta_1)\delta_1\,{}^t(\alpha_i) \equiv (\det \delta_1)\,{}^t(\alpha'_i) \pmod{\mathbb{Z}}.\]
Hence $\psi$ induces a ring isomorphism $Q(A'',B'')\cong \mathbb{Z}+mQ(A',B')$. 
Since $Q(A',B')$ is primitive, we have $\gcd\{\lambda^{ij}_{k\ell}(A'',B'')\}=\ct(A'',B'')=m$. 
We first assume that $m$ divides all coefficients of $A''$. 
We put $\delta_2=\diag[\det \delta_1,1]$ and 
$(A''',B''')=(1_3,\delta_2^{-1})\cdot(A'',B'')=(\delta_1,\delta_2^{-1})\cdot (A,B)$. 
Then $(A''',B''')$ is a pair of integral ternary quadratic forms such that $\ct(A''',B''')=1$. 
We next assume that $m$ does not divide $a''_{ij}$ for some $(i,j)$, 
where we denote by $a''_{k\ell}$ and $b''_{k\ell}$ the coefficients of $v_kv_\ell$ in $A''(v)$ and 
$B''(v)$, respectively. We put $d=\gcd(a''_{ij},b''_{ij})$, $e=\gcd(d,m)$ 
and write $d=ed_1$, $m=em_1$, $\gcd(d_1,m_1)=1$. 
Since $e\neq m$, we have $m_1>1$ and $m_1|m$. 
We write $a''_{ij}=da$, $b''_{ij}=db$ and take $s,t\in \mathbb{Z}$ so that 
$\gamma_2=\left(
\begin{array}{cc}
b & -a \\
s & t
\end{array}
\right) \in \SL_2(\mathbb{Z})$. 
We put $(A''',B''')=(1_3,\gamma_2)\cdot (A'',B'')$ and denote by 
$a'''_{k\ell}$ and $b'''_{k\ell}$ the coefficients of $v_kv_\ell$ in $A'''(v)$ and 
$B'''(v)$, respectively. Then we  have
\[ a'''_{k\ell}=ba''_{k\ell}-ab''_{k\ell}=-d^{-1}(a''_{ij}b''_{k\ell}-b''_{ij}a''_{k\ell})
               =-d^{-1} \lambda^{ij}_{k\ell}(A'',B'').
\]
Since all $\lambda^{ij}_{k\ell}(A'',B'')$ are divisible by $m$, 
$a'''_{k\ell}$ are divisible by $m_1$. 
Hence $(A^{(4)},B^{(4)})=(1_3,\diag[m_1^{-1},1])\cdot (A''',B''')$ is a pair
of integral ternary quadratic forms such that $\ct(A^{(4)},B^{(4)})=m/m_1$.
Repeating the argument, we obtain a matrix $\delta_2\in M_2(\mathbb{Z})$
such that $\det \delta_2=m$ and 
$(A^{(n)},B^{(n)})=(1_3,\delta_2^{-1})\cdot (A'',B'')=(\delta_1,\delta_2^{-1})\cdot (A,B)$ 
is a primitive pair of integral ternary quadratic forms. 
We denote by $\{\alpha''_i\}$ and $\{\alpha^{(n)}_i\}$ 
the basis $\{\alpha_i\}$ for $Q(A'',B'')$ and $Q(A^{(n)},B^{(n)})$, respectively. 
By Proposition \ref{prop:QRINGHOMO}, we have an injective ring homomorphism 
$\phi:Q(A'',B'')\rightarrow Q(A^{(n)},B^{(n)})$ such that 
$\phi(\alpha''_i)=m \alpha^{(n)}_i$ which induces an isomorphism $Q(A'',B'')\cong \mathbb{Z}+mQ(A^{(n)},B^{(n)})$. 
Since $Q(A'',B'')\cong \mathbb{Z}+mQ(A',B')$, 
this implies $Q(A',B') \cong Q(A^{(n)},B^{(n)})$. 
Since $(A',B')$ and $Q(A^{(n)},B^{(n)})$ are primitive, 
the cubic resolvent rings $R(A',B')$ and $R(A^{(n)},B^{(n)})$ are isomorphic to 
each other by Corollary \ref{cor:UNIQUE-RESOLVENT}. 
So $(Q(A',B'),R(A',B'))$ is isomorphic to $(Q(A^{(n)},B^{(n)}),R(A^{(n)},B^{(n)}))$,  
hence there is an element $\gamma=(\gamma_1,\gamma_2) \in \Gamma$ such that 
$(A',B')=\gamma \cdot (A^{(n)},B^{(n)})=(\gamma_1\delta_1,\gamma_2\delta_2^{-1})\cdot (A,B)$ 
by Theorem \ref{thm:BHARGAVA3}. 
\qed

\section{The case $G=S_4$}
\subsection{Cubic resolvent ring of the maximal order}
Let $K$ be a quartic field and $\tilde{K}$ be the normal closure of 
$K$ over $\mathbb{Q}$. We put $G=\Gal(\tilde{K}/\mathbb{Q})$. 
We assume that $G=S_4$ or $A_4$.
We denote by $k$ one of the three conjugate cubic fields contained in $\tilde{K}$ if $G=S_4$, 
the unique cubic field contained in $\tilde{K}$ if $G=A_4$. 
We denote by $k_6$ the non-Galois sextic field such that $k \subset k_6 \subset \tilde{K}$ and 
$k_6=k(\sqrt{\alpha})$ for some $\alpha \in k^\times\smallsetminus (k^\times)^2$ which has square norm in $\mathbb{Q}^\times$. 
We first study a cubic resolvent ring of the maximal order $\mathcal{O}_K$,  
which is unique up to isomorphism by Corollary \ref{cor:UNIQUE-RESOLVENT}. 
Hence it is unique as a subring of $\mathcal{O}_k$ if $G=S_4$. 
By the basic properties of Artin's $L$-function, the following equation on Dedekind zeta functions
holds: 
\begin{equation}\label{eq:DEDEKIND-ZETA}
\zeta(s)\zeta_{k_6}(s)=\zeta_{k}(s)\zeta_{K}(s). 
\end{equation}
By the functional equations of Dedekind zeta functions and \eqref{eq:DEDEKIND-ZETA}, 
we have
\begin{equation}\label{eq:DISCRIMINANT-RELATION}
\Disc(K)=\Disc(k) N(\Disc(k_6/k)).
\end{equation}
By an elementary argument in Galois theory, we obtain the following proposition.
\begin{prop}\label{prop:DISCQC}
The norm of the relative discriminant $N(\Disc(k_6/k))$ is square of a positive integer $f$.
Further if $f$ is odd, then $f$ is square free.
\end{prop}

Let $(A,B)$ be a pair of integral ternary quadratic forms 
such that $Q(A,B) \cong \mathcal{O}_K$. Since $\mathcal{O}_K$ is 
primitive, the $\Gamma$-equivalence class of $(A,B)$ is 
uniquely determined. 
We consider the cubic ring $R(A,B) \subset \mathcal{O}_k$. 
We study the relative discriminant $\Disc(k_6/k)$ and 
the conductor of the cubic ring $R(A,B)$. 
By Proposition \ref{prop:DISCQC} and \eqref{eq:DISCRIMINANT-RELATION}, we have
\[ \Disc(R(A,B))=\Disc(Q(A,B))=\Disc(\mathcal{O}_K)=\Disc(\mathcal{O}_k) f^2, \]
where we put $N(\Disc(k_6/k))=f^2$. 
We now assume that $f$ is square free. Since $f=(\mathcal{O}_k:R(A,B))$, we can take a basis $\{1,\omega,\theta\}$ 
of $\mathcal{O}_k$ such that $\{1,f\omega,\theta\}$ is a basis of $R(A,B)$ and 
the conductor of $R(A,B)$ is $\mathfrak{f}=[f,f\omega,\theta]$ by Proposition \ref{prop:CUNDUCTOR3}. 
So we have $N(\mathfrak{f})=f^2=N(\Disc(k_6/k))$. 
Hence we can expect that the conductor $\mathfrak{f}$ equals the relative discriminant 
$\Disc(k_6/k)$. 

\begin{lem}\label{lem:COFACTOR-MATRIX}
Let $S=(s_{ij})$ be a matrix of degree three and $t_{ij}$ be the $(i,j)$-cofactor of $S$. 
Put $T={}^t(t_{ij})$. Then the $(i,j)$-cofactor of $T$ is $(\det S)s_{ji}$. 
In particular, $t_{ii}t_{jj}=t_{ij}^2$ if $S$ is symmetric and $\det S=0$.
\end{lem}
\Proof
We have $ST=(\det S)1_3$, $\det T=(\det S)^2$. 
We denote by $u_{ij}$ the $(i,j)$-cofactor of $T$ and put $U={}^t(u_{ij})$. 
Then $UT=(\det S)^21_3$, hence $\left(U-(\det S)S \right)T=0$. 
If $\det S\neq 0$, then $U-(\det S)S=0$, hence $u_{ij}=(\det S)s_{ji}$. 
Since these equations are polynomial identities in $s_{ij}$'s, $u_{ij}=(\det S)s_{ji}$ still hold 
when $\det S=0$. If $S$ is symmetric and $\det S=0$, then 
$t_{ii}t_{jj}-t_{ij}^2=u_{\ell\ell}=0$ for $\{i,j,\ell\}=\{1,2,3\}$.  
\qed

We identify $(A,B)$ by a pair of symmetric matrices of degree three whose entries are half integers. 
\begin{lem}\label{lem:COFACTOR-K6}
Let $\Delta_{ij}$ and $\Delta'_{ij}$ be the $(i,j)$-cofactors of $\omega A + aB$ and $dA-\theta B$, respectively. 
Then we have $\Delta_{ii}\Delta_{jj}=\Delta_{ij}^2$, $\Delta'_{ii}\Delta'_{jj}=(\Delta'_{ij})^2$ and 
$a^2\Delta'_{ii}=\theta^2 \Delta_{ii}$ for all $i, j$. 
We also have $\Delta_{ii} \neq 0$ for some $i$ and $k_6=k(\sqrt{-\Delta_{ii}})$. 
\end{lem}
\Proof
Since $4\det(\omega A+ aB)=F_{(A,B)}(\omega,-a)=0$, by Lemma \ref{lem:COFACTOR-MATRIX} we have  
$\Delta_{ii}\Delta_{jj}=\Delta_{ij}^2$. Similarly we have $\Delta'_{ii}\Delta'_{jj}=(\Delta'_{ij})^2$. 
If $\{i,j,\ell\}=\{1,2,3\}$, then we have
\begin{align*}
a^2\Delta'_{ii}
 &=a^2(da_{jj}-\theta b_{jj})(da_{\ell\ell}-\theta b_{\ell\ell})-\frac{a^2}{4}(da_{j\ell}-\theta b_{j\ell})^2 \\
 &=\theta^2(\omega a_{jj}+ab_{jj})(\omega a_{\ell\ell}+a b_{\ell\ell})
 -\frac{\theta^2}{4}(\omega a_{j\ell} + a b_{j\ell})^2
 =\theta^2\Delta_{ii}.
\end{align*}
Since $(2A,2B)$ is a nondegenerate pair of integral symmetric matrices of degree three and 
$\hat{F}_{(2A,2B)}(u)=2F_{(A,B)}(u)$, 
there exists an ideal $\mathfrak{a}=[\alpha_1,\alpha_2,\alpha_3]$ of $R(2F_{(A,B)})=[1,2\omega,2\theta]$ 
and $\kappa \in k^\times$ such that 
\[ 2A=-\left(\trace(\Delta^{-1}\kappa^{-1}\alpha_i\alpha_j) \right), \quad 
   2B=-\left(\trace(\Delta^{-1}\kappa^{-1}\alpha_i\alpha_j(-\omega/a)) \right),
\]
where we put $\Delta=3a(-\omega/a)^2+2b(-\omega/a)+c \in k$ as we saw just after \eqref{eq:UPLUS2-ORDER} in \S 4. 
We denote by $\alpha^{(i)}_j$ the conjugates of $\alpha_j$ and put $U=(\alpha^{(i)}_j) \in \GL_3(\mathbb{C})$. 
If we put $\beta=-\Delta^{-1}\kappa^{-1}$ and $\lambda=-\omega/a$, then we have 
\begin{align*}
2A&=\left(\trace(\beta\alpha_i\alpha_j) \right)
   ={}^t U\diag[\beta^{(1)},\beta^{(2)},\beta^{(3)}]\, U, \\
2B&=\left(\trace(\beta\alpha_i\alpha_j\lambda) \right)={}^t U\diag[\beta^{(1)}\lambda^{(1)},\beta^{(2)}\lambda^{(2)},\beta^{(3)}\lambda^{(3)}]\, U.
\end{align*}
Hence
\begin{align*}
2(\omega^{(1)}A+aB)
 &=a\,{}^t U \diag[0,\beta^{(2)}(\lambda^{(2)}-\lambda^{(1)}),\beta^{(3)}(\lambda^{(3)}-\lambda^{(1)})]\, U .
\end{align*}
This proves that $\omega^{(1)}A+aB$ has rank two. 
So we have $\Delta^{(1)}_{ij} \neq 0$ for some $(i,j)$. 
Since $\Delta_{ii}\Delta_{jj}=\Delta_{ij}^2$, we have $\Delta_{ii} \neq 0$. 
We may assume $\Delta_{33} \neq 0$. 
We put $S=\omega A+ aB$ and 
$P=\left(
\begin{array}{ccc}
1 & 0 & 0 \\
0 & 1 & 0 \\
\Delta_{13} & \Delta_{23} & \Delta_{33}
\end{array}
\right)$. 
Since $\Delta_{33} \neq 0$, we have $P \in \GL_3(k)$ and 
\[ PS\,{}^t\!P = \left(
\begin{array}{cc}
S_1 & 0 \\
0 & 0
\end{array}
\right), \quad S_1=\left(
\begin{array}{cc}
\omega a_{11}+ab_{11} & (\omega a_{12}+a b_{12})/2 \\
(\omega a_{12}+a b_{12})/2 & \omega a_{22}+ab_{22}
\end{array}
\right). 
\]
Since $\det S_1=\Delta_{33}\neq 0$, the quadratic form $\omega A(v)+ aB(v)$ 
is decomposed into a product of two distinct linear forms over $k(\sqrt{-\Delta_{33}})$.
Therefore $\tilde{K}$ is obtained by adjoining the square roots of the conjugates of $-\Delta_{33}$. 
We can check that the norm of $-4\Delta_{33}$ is square in $\mathbb{Q}^\times$. 
Hence $k_6=k(\sqrt{-\Delta_{33}})$. 
\qed

For a prime number $p$ and a positive integer $f$, 
we denote by $\pi_p(f)$ the number of monic irreducible polynomials of degree $f$ 
over $\mathbb{F}_p$. Then $\pi_p(f)$ for $f\leq 4$ are given by 
\begin{align*}
\pi_p(1) &=p, \quad 
\pi_p(2)=\frac{p^2-p}{2}, \quad
\pi_p(3)=\frac{p^3-p}{3}, \quad 
\pi_p(4)=\frac{p^4-p^2}{4}.
\end{align*}
Let $F$ be a number field of degree $n$ over $\mathbb{Q}$ 
and $\theta$ be a primitive element in $\mathcal{O}_F$. 
Then $\Disc(\mathbb{Z}[\theta])=(\mathcal{O}_F:\mathbb{Z}[\theta])^2\Disc(F)$. 
The greatest common divisor of indices $(\mathcal{O}_F:\mathbb{Z}[\theta])$ of all primitive elements 
$\theta \in \mathcal{O}_F$ is called the \textit{inessential discriminant divisor} of $F$ and denoted by $i(F)$. 
For a prime number $p$ and a positive integer $f$, we denote by $r_p(f)$ the number of prime ideal divisors of $p$ in $\mathcal{O}_F$ of degree $f$. 
Then $p|i(F)$ if and only if $r_p(f) > \pi_p(f)$ for some $f$ by the criterion in Hasse \cite[p. 456]{HA}. 
By this criterion, if $p|i(F)$, then $p<n$. 
If $p<n$ and $p$ splits completely in $F$, then $p|i(F)$. 
So we obtain the following lemma for our quartic field $K$. 
\begin{lem}\label{lem:INESSTIAL-DIVISOR}
If a prime number $p$ divides $i(K)$, then $p=2$ or $p=3$. 
The prime number $3$ divides $i(K)$ if and only if $3$ is of type $1111$ in $K$. 
The prime number $2$ divides $i(K)$ if and only if $2$ is of type $1111$, $111^2$, or $22$.
\end{lem}

From now on, we restrict ourselves to the case of $G=S_4$. We shall discuss the case of $G=A_4$ in the next section. 
We consider the following subgroups of $S_4$. 
\begin{align*}
S_3 &=\{1,(2\,3),(2\,4),(3\,4),(2\,3\,4),(2\,4\,3)\}, \\
A_3 &=\{1,(2\,3\,4),(2\,4\,3)\}, \\
D_4 &=\{1,(1\,2\,3\,4), (1\,3)(2\,4), (1\,4\,3\,2), (1\,3), (2\,4), (1\,2)(3\,4),(1\,4)(2\,3) \}, \\
C_4 &=\{1,(1\,2\,3\,4), (1\,3)(2\,4), (1\,4\,3\,2)\}, \\
V_4 &=\{1,(1\,3)(2\,4),(1\,2)(3\,4),(1\,4)(2\,3)\}, \\
B_4 &=\{1, (1\,3), (2\,4), (1\,3)(2\,4)\},  \\
B_2 &=\{1, (1\,3)\}, \\
C_2 &=\{1, (1\,3)(2\,4)\}.
\end{align*}
We may assume $\Gal(\tilde{K}/K)=S_3$, $\Gal(\tilde{K}/k)=D_4$, $\Gal(\tilde{K}/k_6)=B_4$. 
For a prime number $p$, we take a prime ideal divisor $\mathfrak{P}$ in $\mathcal{O}_{\tilde{K}}$ of $p$ 
and denote by $Z$, $T$ and $V$ the decomposition group, 
the inertia group and the first ramification group of $\mathfrak{P}$ in $\tilde{K}/\mathbb{Q}$, respectively. 
Then the prime ideal decompositions of $p$ in $K$, $k$ and $k_6$ are given by the Tables  \ref{tabular:T1}, \ref{tabular:TV1}, \ref{tabular:TVNE1} and  \ref{tabular:TVNE2}. 

\begin{table}
\begin{center}
\begin{tabular}{cccccc} 
$T$ & $Z$ & $K$ & $k$ & $k_6$ & $k_6/k$ \\ \hline
$1$ &$1$ & $1111$ & $111$ & $111111$ & $p\nmid f$ \\ 
$1$ &$B_2$ & $112$ & $12$ & $1122$ & $p\nmid f$\\ 
$1$ &$C_2$ & $22$ & $111$ & $1122$ & $p\nmid f$\\ 
$1$ &$A_3$ & $13$ & $3$ & $33$ & $p\nmid f$\\ 
$1$ &$C_4$ & $4$  & $12$ & $24$ & $p\nmid f$\\ 
\end{tabular}
\caption{prime ideal decompositions of $p$ for $T=1$}
\label{tabular:T1}
\end{center}
\end{table}

\begin{table}
\centering
\begin{tabular}{cccccc} 
$T$ & $Z$ & $K$ & $k$ & $k_6$ & $k_6/k$ \\ \hline
$B_2$ &$B_2$ & $111^2$ & $11^2$ & $111^21^2$ & $p\nmid f$\\ 
$B_2$ &$B_4$ & $21^2$ & $11^2$ & $112^2$ & $p\nmid f$\\ 
$C_2$ &$C_2$ & $1^21^2$ & $111$ & $111^21^2$ & $p|f$\\ 
$C_2$ &$B_4$ & $1^21^2$ & $12$ & $112^2$ & $p|f$\\ 
$C_2$ &$V_4$ & $2^2$    & $111$ & $21^21^2$ & $p|f$\\ 
$C_2$ &$C_4$ & $2^2$    & $12$ & $22^2$ & $p|f$\\ 
$A_3$ &$A_3$, $S_3$ & $11^3$    & $1^3$ & $1^31^3$ & $p\nmid f$\\ 
$C_4$ &$C_4$, $D_4$ & $1^4$    & $11^2$ & $1^21^4$ & $p|f$\\ 
\end{tabular}
\caption{prime ideal decompositions of $p$ for cyclic $T\neq 1$}
\label{tabular:TV1}
\end{table}

\begin{table}
\centering
\begin{tabular}{ccccccc} 
$V$ & $T$ & $Z$ & $K$ & $k$ & $k_6$ & $k_6/k$\\ \hline
$A_3$&$A_3$ &$A_3$, $S_3$ & $11^3$    & $1^3$ & $1^31^3$ & $3\nmid f$\\ 
$A_3$&$S_3$ &$S_3$ & $11^3$    & $1^3$ & $1^31^3$ & $3\nmid f$\\ 
\end{tabular}
\caption{prime ideal decompositions of $3$ for $V\neq 1$}
\label{tabular:TVNE1}
\end{table}

\begin{table}
\centering
\begin{tabular}{ccccccc} 
$V$ & $T$ & $Z$ & $K$ & $k$ & $k_6$ & $k_6/k$\\ \hline
$V_4$&$A_4$ &$S_4$ & $1^4$    & $1^3$ & $1^6$ & $2|f$ \\ 
$D_4$&$D_4$ &$D_4$ & $1^4$    & $11^2$ & $1^21^4$ & $2^2|f$ \\ 
$B_4$&$B_4$ &$B_4$ & $1^21^2$ & $11^2$ & $111^4$  & $2|f$ \\ 
$B_4$&$B_4$ &$D_4$ & $2^2$    & $11^2$ & $21^4$   & $2|f$ \\ 
$V_4$&$V_4$ &$V_4$ & $1^4$    & $111$  & $1^21^21^2$ & $2^3|f$ \\ 
$V_4$&$V_4$ &$D_4$ & $1^4$    & $12$ & $1^22^2$ & $2^3|f$ \\ 
$V_4$&$V_4$ &$A_4$ & $1^4$    & $3$ & $3^2$ & $2^3|f$ \\ 
\end{tabular}
\caption{prime ideal decompositions of $2$ for non-cyclic $T$}
\label{tabular:TVNE2}
\end{table}

\begin{lem}\label{lem:F-DEVIDED-EXACT2}
Let $\mathfrak{d}=\Disc(k_6/k)$ be the relative discriminant of $k_6/k$ and write 
$N(\mathfrak{d})=f^2$. We assume $f$ is square free. Let $p$ be a prime number satisfying $p|f$ and denote by $\mathfrak{d}_p$ 
the $p$-part of $\mathfrak{d}$. 
If $p\neq 2$, then $\mathfrak{d}_p$ is a square free ideal. 
If $p=2$, then $2$ is of type $11^2$ or $1^3$ in $k$ and 
$\mathfrak{d}_2=\mathfrak{p}^2$ where $\mathfrak{p}$ is the prime ideal dividing $2$ and ramified in $k/\mathbb{Q}$.
\end{lem}
\Proof
If $\mathfrak{p}$ is an odd prime ideal of $\mathcal{O}_k$ dividing $\mathfrak{d}$, 
then $\mathfrak{p}$ divides exactly $\mathfrak{d}$ by Kummer theory.  Hence $\mathfrak{d}_p$ is a square free ideal. 
Let $\mathfrak{p}$ be a prime ideal of $\mathcal{O}_k$ dividing $2$. 
We denote by $e$ and $n$ the ramification index 
and the residue degree of $\mathfrak{p}$ in $k/\mathbb{Q}$. 
If $\mathfrak{p}$ divides $\mathfrak{d}$, then $\mathfrak{p}^2$ divides $\mathfrak{d}$ by Kummer theory.
Hence $2^{2n}$ divides $f^2$. Since $f$ is square free, we have $n=1$. 
We also see that such prime ideal $\mathfrak{p}$ is unique. 
It follows from tables \ref{tabular:TV1} and \ref{tabular:TVNE2} that $2$ is of type $11^2$ or $1^3$ in $k$, $e=2$ and 
$\mathfrak{d}_2=\mathfrak{p}^2$.
\qed

Let $p$ be a prime number satisfying $p|f$. Since $f$ is square free, 
$p$ is of type $1^21^2$, $2^2$ or $1^4$ in $K$ by the tables. 
Hence $p \nmid i(K)$. 
Thus there exists an element $\rho \in \mathcal{O}_K$ such that $K=\mathbb{Q}(\rho)$ and 
$p \nmid (\mathcal{O}_K:\mathbb{Z}[\rho])$. 
Let $\varphi(x)=x^4+a_1x^3+a_2x^2+a_3x+a_4 \in \mathbb{Z}[x]$ be the minimal polynomial of $\rho$ over $\mathbb{Q}$. 
Since $p$ is of type $1^21^2$, $2^2$ or $1^4$ in $K$, we may assume that if $p\neq 2$, 
\[ \varphi(x)\equiv (x^2-a)^2 \pmod{p},\]
where $a \in \mathbb{Z}$ is a quadratic residue modulo $p$, 
a quadratic nonresidue modulo $p$ or $a=0$ according as 
$p$ is of type $1^21^2$, $2^2$ or $1^4$. 
If $p=2$, we may assume that 
\[ \varphi(x) \equiv \left\{
\begin{array}{ll}
x^2(x+1)^2 \pmod{2}, &\quad  \text{$2$ is of type $1^21^2$ in $K$,}\\
(x^2+x+1)^2 \pmod{2}, &\quad \text{$2$ is of type $2^2$ in $K$,} \\
x^4 \pmod{2}, &\quad \text{$2$ is of type $1^4$ in $K$.}
\end{array}
\right.
\]
Let $\rho_i$ $(i=1,2,3,4$) be the conjugates of $\rho$ and put
\[ \xi_1=(\rho_1+\rho_2)(\rho_3+\rho_4),\; \xi_2=(\rho_1+\rho_3)(\rho_2+\rho_4), \;
\xi_3=(\rho_1+\rho_4)(\rho_2+\rho_3).
\]
If we put $g(x)=(x-\xi_1)(x-\xi_2)(x-\xi_3)=x^3+b_1x^2+b_2x+b_3$, then 
\[ b_1=-2a_2, \quad b_2=a_2^2+a_1a_3-4a_4, \quad b_3=a_3^2-a_1a_2a_3+a_1^2a_4. \]
We put $\alpha_1=\rho$, $\alpha_2=\rho^2+a_1\rho+a_2$, 
$\alpha_3=\rho^3+a_1\rho^2+a_2\rho+a_3$. 
Then $\mathbb{Z}[\rho]=[1,\alpha_1,\alpha_2,\alpha_3]$. 
The ring structure of $\mathbb{Z}[\rho]$ is given by
\begin{align*}
\alpha_1^2 &=-a_2-a_1\alpha_1+\alpha_2, \\
\alpha_1\alpha_2 &=-a_3+\alpha_3, \\
\alpha_1\alpha_3 &=-a_4, \\
\alpha_2^2 &=-a_1a_3-a_4-a_3\alpha_1+a_2\alpha_2+a_1\alpha_3, \\
\alpha_2\alpha_3 &=-a_1a_4-a_4\alpha_1+a_2\alpha_3, \\
\alpha_3^2 &=-a_4\alpha_2+a_3\alpha_3.
\end{align*}
If we put
\begin{equation}\label{eq:BINARY-QUARTIC-PAIR}
\begin{split}
 A_1(v)&= -v_1v_3+v_2^2, \\
 B_1(v)&= v_1^2+a_1v_1v_2+a_2v_1v_3+a_3v_2v_3+a_4v_3^2,
\end{split}
\end{equation}
then we have $Q(A_1,B_1)=\mathbb{Z}[\rho]$ and $F_{(A_1,B_1)}(u_1,1)=g(-u_1)$. 
Moreover $\lambda^{11}_{13}(A_1,B_1)=1$ implies that $(A_1,B_1)$ is primitive. 
Since $Q(A,B)=\mathcal{O}_K$, $(A,B)$ is also primitive. 
It follows from Proposition \ref{prop:QUARTIC-SUBRING} that there exist matrices 
$\delta_1 \in M_3(\mathbb{Z})$ and $\delta_2 \in M_2(\mathbb{Z})$ such that 
$(A_1,B_1)=(\delta_1,\delta_2^{-1})\cdot (A,B)$ and $|\det \delta_1|=|\det \delta_2|=(\mathcal{O}_K:\mathbb{Z}[\rho]))$. 
We put $m=(\mathcal{O}_K:\mathbb{Z}[\rho]))$ 
and denote by $\mathfrak{f}$ and $\mathfrak{f}_1$ 
the conductor of $R(A,B)$ and that of $R(A_1,B_1)$, respectively. 
We put $\delta'_2=\diag[1,-1](\det \delta_2)\delta_2^{-1}\diag[1,-1]$. 
Then we have $\det \delta'_2=\det \delta_2$ and $F_{(A_1,B_1)}(u)=(\delta'_2 F_{(A,B)})(u)$. 
It follows from Proposition \ref{prop:CUBIC-SUBRING} that $R(A_1,B_1)$ is a subring of $R(A,B)$ 
with index $m=|\det \delta_2|$. 
Since $f=(\mathcal{O}_k:R(A,B))$, $p|f$ and $p \nmid m$, 
it follows from Lemma \ref{lem:CONDUCTOR1} that 
the $p$-part of $\mathfrak{f}_1$ equals the $p$-part of $\mathfrak{f}$.

We first assume that $p\neq 2$ or $p=2$ is of type $1^4$ in $K$. 
Since $a_1\equiv a_3\equiv 0 \pmod{p}$, $a_2\equiv -2a \pmod{p}$ and $a_4\equiv a^2 \pmod{p}$ for $p\neq 2$ 
and $a_i \equiv 0 \pmod{2}$ ($i=1,2,3,4$) for $p=2$, we have
\begin{equation}\label{eq:RESOLVENT-CUBIC-MODP}
 b_2 \equiv 0 \pmod{p}, \quad b_3 \equiv 0 \pmod{p^2}.
\end{equation}
Moreover we have $b_1\equiv b_2\equiv b_3\equiv 0 \pmod{4}$ if $p=2$. 
If we put 
\[ F_0(u)=pF_{(A_1,B_1)}(u_2,u_1/p)=(b_3/p^2)u_1^3-(b_2/p)u_1^2u_2+b_1u_1u_2^2-pu_2^3, \]
then $F_0(u)$ is an integral binary cubic form by the congruences \eqref{eq:RESOLVENT-CUBIC-MODP}. 
Let $\{1,\eta_1,\eta_2\}$ be the $\mathbb{Z}$-basis of $R(F_0)$ corresponding to $F_0$. 
It follows from Proposition \ref{prop:CUBIC-SUBRING} that 
$R(A_1,B_1)=R(F_{(A_1,B_1)})$ is a subring of $R(F_0)$ 
with $\mathbb{Z}$-basis $\{1,p\eta_1,\eta_2\}$. We put $\mathfrak{g}=[p,p\eta_1,\eta_2]$. 
Then it easy to see that $\mathfrak{g}$ is an $R(F_0)$-ideal contained in $R(A_1,B_1)$. 
Since the index $[R(A_1,B_1):\mathfrak{g}]=p$ is a prime number, we see that 
$\mathfrak{g}$ is the largest $R(F_0)$-ideal contained in $R(A_1,B_1)$. 
Since $(\mathcal{O}_k:R(F_0))=mf/p$ is prime to $p$ and $(R(F_0):\mathfrak{g})=p^2$, 
it follows from Lemma \ref{lem:CONDUCTOR2} that $\mathfrak{g}\mathcal{O}_k$ is the $p$-part of $\mathfrak{f}_1$. 
Hence the $p$-part of $\mathfrak{f}$ is given by $\mathfrak{f}_p=\mathfrak{g}\mathcal{O}_k$. 
Since $pR(F_0) \subset \mathfrak{g}$, we have $p\mathcal{O}_k \subset \mathfrak{g}\mathcal{O}_k$. 
Since $R(F_0)=[1,\eta_1,\eta_2]$ and 
$(\mathcal{O}_k:R(F_0))=mf/p$ is prime to $p$, we have $\eta_2 \notin p\mathcal{O}_k$. 

If $p=2$, then $F_0(u)\equiv (b_3/4) u_1^3 \pmod{2}$. 
We see that $b_3/4$ is odd since the index $(\mathcal{O}_k:R(F_0))=mf/2$ is odd. 
Hence we have
\[ \eta_1^2 \equiv \eta_2, \quad \eta_2^2 \equiv 0,\quad \eta_1\eta_2 \equiv 0 \pmod{2R(F_0)}.\]
If we put $\mathfrak{p}=[2,\eta_1,\eta_2]$, then $\mathfrak{p}$ is an ideal of $R(F_0)$. 
It is easy to see that $\mathfrak{p}^2=[2,2\eta_1,\eta_2]=\mathfrak{g}$ and $\mathfrak{p}^3=2R(F_0)$. 
Hence $\tilde{\mathfrak{p}}=\mathfrak{p}\mathcal{O}_k$ is the prime ideal of $\mathcal{O}_k$ 
such that $2\mathcal{O}_k=\tilde{\mathfrak{p}}^3$. We have $\mathfrak{f}_2=\mathfrak{g}\mathcal{O}_k=\tilde{\mathfrak{p}}^2$. 
Since $\tilde{\mathfrak{p}}$ ramifies in $k_6/k$ and $f$ is a square free even integer, 
the $2$-part of $\Disc(k_6/k)$ equals $\tilde{\mathfrak{p}}^2=\mathfrak{f}_2$. 

%
%
We now assume $p\neq 2$. 
We denote by $\Delta_{33}$ the $(3,3)$-cofactor of $\eta_2 A_1+B_1$ and put 
$\beta=-4\Delta_{33}$. Then we have $\beta=a_1^2-4\eta_2$. 
Since $a_1 \equiv 0 \pmod{p}$, we have 
$\beta=a_1^2-4\eta_2 \in \mathfrak{g} \subset \mathfrak{g}\mathcal{O}_k=\mathfrak{f}_p$. 
By Lemma \ref{lem:COFACTOR-K6}, we have $k_6=k(\sqrt{\beta})$. 
We denote by $\mathfrak{d}_p$ the $p$-part of the relative discriminant $\Disc(k_6/k)$. 
We now show $\mathfrak{d}_p=\mathfrak{f}_p$ using $\beta \in \mathfrak{f}_p$ and $N(\mathfrak{f}_p)=p^2$. 
If $p$ is of type $12$ in $k$, write $p\mathcal{O}_k=\mathfrak{p}_1\mathfrak{p}_2$, $N(\mathfrak{p}_i)=p^i$ $(i=1,2)$. 
Since $\mathfrak{d}_p$ is square free and has norm $p^2$, we have $\mathfrak{d}_p=\mathfrak{p}_2$. 
On the other hand, $N(\mathfrak{f}_p)=p^2$ implies $\mathfrak{f}_p=\mathfrak{p}_2$ or $\mathfrak{f}_p=\mathfrak{p}_1^2$. 
Since $p \in \mathfrak{g} \subset \mathfrak{f}_p$, the latter case is impossible. 
Hence we have $\mathfrak{f}_p=\mathfrak{p}_2=\mathfrak{d}_p$. 
If $p$ is of type $111$ in $k$, write $p\mathcal{O}_k=\mathfrak{p}_1\mathfrak{p}_2\mathfrak{p}_3$. 
Then we have $\mathfrak{d}_p=\mathfrak{p}_{i_1}\mathfrak{p}_{i_2}$ for some $i_1\neq i_2$. 
On the other hand, $\mathfrak{f}_p=\mathfrak{p}_{j_1}\mathfrak{p}_{j_2}$ for some $j_1\neq j_2$ 
or $\mathfrak{f}_p=\mathfrak{p}_j^2$. 
The latter case is impossible as in the previous case. 
Suppose $\{j_1,j_2\} \neq \{i_1,i_2\}$. Then we have $\beta \in \mathfrak{d}_p \cap \mathfrak{f}_p = p\mathcal{O}_k$. 
This contradicts $\eta_2 \notin p\mathcal{O}_k$. So we have $\mathfrak{d}_p=\mathfrak{f}_p$. 
If $p$ is of type $11^2$ in $k$, write $p\mathcal{O}_k=\mathfrak{p}_1\mathfrak{p}_2^2$. 
Then we have $\mathfrak{d}_p=\mathfrak{p}_1\mathfrak{p}_2$. 
On the other hand, $\mathfrak{f}_p=\mathfrak{p}_1\mathfrak{p}_2$ or $\mathfrak{f}_p=\mathfrak{p}_i^2$. 
Suppose $\mathfrak{f}_p=\mathfrak{p}_2^2$. 
Then we have $\beta \in \mathfrak{d}_p \cap \mathfrak{f}_p = p\mathcal{O}_k$. 
This contradicts $\eta_2 \notin p\mathcal{O}_k$. 
Suppose $\mathfrak{f}_p=\mathfrak{p}_1^2$. 
This contradicts $p \in \mathfrak{g} \subset \mathfrak{f}_p$. 
So we have $\mathfrak{f}_p=\mathfrak{p}_1\mathfrak{p}_2=\mathfrak{d}_p$. 
If $p$ is of type $3$ in $k$, there exists no ideal of $\mathcal{O}_k$ having norm $p^2$. 
If $p$ is of type $1^3$ in $k$, there exists no square free ideal of $\mathcal{O}_k$ having norm $p^2$. 
So we have proved $\mathfrak{d}_p=\mathfrak{f}_p$ for any prime number $p|f$ if $p\neq 2$ 
or $p=2$ is of type $1^4$ in $K$. 

We next assume $p=2$ is of type $1^21^2$ or $2^2$ in $K$. 
Since $a_1\equiv a_3 \equiv 0 \pmod{2}$ and $a_2 \equiv 1 \pmod{2}$, 
we have $b_1\equiv 2 \pmod{4}$, $b_2\equiv 1 \pmod{4}$ and $b_3 \equiv 0 \pmod{4}$. 
Hence 
\begin{align*}
g(1)
 &=1+b_1+b_2+b_3 \equiv 0 \pmod{4}, \\
g(-1)
 &=-1+b_1-b_2+b_3 \equiv 0 \pmod{4}.
\end{align*}
Moreover we have
\begin{align*}
g(1)+g(-1)
 &=2b_1+2b_3=2(-2a_2+a_3^2-a_1a_2a_3+a_1^2a_4) \\
 &\equiv -4a_2 \equiv 4 \pmod{8}.
\end{align*}
Hence we can chose $\varepsilon=\pm 1$ so that $g(\varepsilon)\equiv 4 \pmod{8}$ and 
$g(-\varepsilon)\equiv 0 \pmod{8}$. We put $\gamma_2=\left(
\begin{array}{cc}
\varepsilon &  -1\\
0 & 1
\end{array}
\right) \in \GL_2(\mathbb{Z})$ and 
\[ (A'_1,B'_1)=(1_3,\gamma_2)\cdot(A_1,B_1)=(\varepsilon A_1-B_1, B_1). \]
Then we have
\begin{align*}
\lefteqn{F_{(A'_1,B'_1)}(u)=F_{(A_1,B_1)}(\varepsilon u_1,u_1+u_2)} \\
 &=g(-\varepsilon)u_1^3+(b_1-2b_2 \varepsilon+3b_3)u_1^2u_2+(-b_2\varepsilon+3b_3)u_1u_2^2+b_3u_2^3  \\
 &\equiv (-\varepsilon) u_1u_2^2 \pmod{4}.
\end{align*}
If we write $F_{(A'_1,B'_1)}(u)=a'u_1^3+b'u_1^2u_2+c'u_1u_2^2+d'u_2^3$, then 
$a'\equiv 0 \pmod{8}$, $b'\equiv d' \equiv 0 \pmod{4}$ and $c'\equiv 1 \pmod{2}$. 
We put $a_0=a'/4$, $b_0=b'/2$, $c_0=c'$, $d_0=2d'$ and $F_0(u)=2F_{(A'_1,B'_1)}(u_1/2,u_2)$. 
Then $F_0(u)=a_0u_1^3+b_0u_1^2u_2+c_0u_1u_2^2+d_0u_2^3$ is an integral binary cubic form 
such that $F_0(u)\equiv u_1u_2^2 \pmod{2}$. 
Let $\{1,\eta_1,\eta_2\}$ be the $\mathbb{Z}$-basis of $R(F_0)$ corresponding to $F_0$. 
Then $R(A_1,B_1)=R(A'_1,B'_1)$ is a subring of $R(F_0)$ 
with $\mathbb{Z}$-basis $\{1,2\eta_1,\eta_2\}$. 
The largest $R(F_0)$-ideal contained in $R(A_1,B_1)$ is $\mathfrak{g}=[2,2\eta_1,\eta_2]$ 
and $\mathfrak{f}_2=\mathfrak{g}\mathcal{O}_k$ as we saw in the previous case. 
Since $2R(F_0) \subset \mathfrak{g}$, we have $2\mathcal{O}_k \subset \mathfrak{g}\mathcal{O}_k$. 
Since $R(F_0)=[1,\eta_1,\eta_2]$ and 
$(\mathcal{O}_k:R(F_0))=mf/2$ is odd, we have $\eta_2 \notin 2\mathcal{O}_k$. 
We see that $a_0\not\equiv 0 \pmod{4}$ since the index $(\mathcal{O}_k:R(F_0))=mf/2$ is odd. 
If we put $\mathfrak{p}_1=[2,\eta_1,\eta_2-1]$ and $\mathfrak{p}_2=[2,\eta_1,\eta_2]$, 
then we see that $\mathfrak{p}_1$ and $\mathfrak{p}_2$ are ideals of $R(F_0)$ satisfying
\[ \mathfrak{p}_1^2=[4,\eta_1,\eta_2+c_0], \quad \mathfrak{p}_2^2=[2,2\eta_1,\eta_2]=\mathfrak{g}, \quad \mathfrak{p}_1\mathfrak{p}_2^2=R(F_0). \]
Hence $\tilde{\mathfrak{p}}_i=\mathfrak{p}_i\mathcal{O}_k$ are the prime ideal of $\mathcal{O}_k$ 
such that $2\mathcal{O}_k=\tilde{\mathfrak{p}}_1\tilde{\mathfrak{p}}_2^2$. 
We have $\mathfrak{f}_2=\mathfrak{g}\mathcal{O}_k=\tilde{\mathfrak{p}}_2^2$. 
By Lemma \ref{lem:F-DEVIDED-EXACT2}, $\mathfrak{d}_2$ equals $\tilde{\mathfrak{p}}_2^2=\mathfrak{f}_2$. 

Hence we have 
\[ \mathfrak{f}=\prod_{p|f} \mathfrak{f}_p= \prod_{p|f} \mathfrak{d}_p = \Disc(k_6/k).\]
Thus we obtain the following proposition. 
\begin{prop}\label{prop:MAXIMAL-ORDER-CONDUCTOR}
Let $K$ be a quartic field with Galois group $G=S_4$ and $\Disc(K)=\Disc(k)f^2$. 
If $f$ is square free, then the relative discriminant $\Disc(k_6/k)$ equals the conductor of the unique cubic resolvent ring of $\mathcal{O}_K$. 
\end{prop}

\subsection{Number of $S_4$-quartic fields}
We next study the number of $S_4$-quartic fields with a fixed cubic resolvent field $k$. 
Here $k$ is a non-Galois cubic field. 
We summarize some basic facts on quadratic residues over $k$ by \cite[Chapter VIII]{HE}. 

We say that an algebraic integer $\alpha \in \mathcal{O}_k$ or an integral ideal $\mathfrak{a}$ of $\mathcal{O}_k$ 
is \textit{odd} if $\alpha$ or $\mathfrak{a}$ is relatively prime to $2$. 
Let $\mathfrak{p}$ be an odd prime ideal of $\mathcal{O}_k$ and $\alpha \in \mathcal{O}_k$ 
be an algebraic integer which is relatively prime to $\mathfrak{p}$. 
We say that $\alpha$ is a \textit{quadratic residue modulo $\mathfrak{p}$} 
if there exists an element $\xi \in \mathcal{O}_k$ such that $\alpha\equiv \xi^2 \pmod{\mathfrak{p}}$, 
otherwise $\alpha$ is a \textit{quadratic nonresidue modulo $\mathfrak{p}$}.  
We define the Legendre symbol $\legend{\alpha}{\mathfrak{p}}$ by 
$\legend{\alpha}{\mathfrak{p}}=1$ or $-1$ according as $\alpha$ is a quadratic residue or nonresidue modulo $\mathfrak{p}$. 
Further we put $\legend{\alpha}{\mathfrak{p}}=0$ if $\alpha \equiv 0 \pmod{\mathfrak{p}}$. 
If an odd integral ideal $\mathfrak{n}$ is written as $\mathfrak{n}=\mathfrak{p}_1\cdots \mathfrak{p}_r$ where $\mathfrak{p}_i$'s are prime ideals, 
then we define $\legend{\alpha}{\mathfrak{n}}$ for any $\alpha \in \mathcal{O}_k$ by 
\[ \legend{\alpha}{\mathfrak{n}}=\prod_{i=1}^r \legend{\alpha}{\mathfrak{p}_i}.\]
If $\mathfrak{n}=(\nu)$ is a principal ideal, we simply write $\legend{\alpha}{\mathfrak{n}}=\legend{\alpha}{\nu}$. 
For any $\alpha, \beta \in \mathcal{O}_k$, $\legend{\alpha}{\mathfrak{n}}=\legend{\beta}{\mathfrak{n}}$ if $\alpha\equiv \beta \pmod{\mathfrak{n}}$. 
Further we have $\legend{\alpha\beta}{\mathfrak{n}}=\legend{\alpha}{\mathfrak{n}}\legend{\beta}{\mathfrak{n}}$. 
We say that $\alpha \in \mathcal{O}_k$ is \textit{primary} 
if there exists an odd element $\xi \in \mathcal{O}_k$ such that $\alpha\equiv \xi^2 \pmod{4\mathcal{O}_k}$. 
We need a lemma on Gauss sums over the cubic field $k$. 
We denote by $\mathcal{D}_k$ the different of $k/\mathbb{Q}$. 
By definition,  $\mathcal{D}_k^{-1}=\{\eta \in k\,|\, \trace_{k/\mathbb{Q}}(\eta\mathcal{O}_k)\subset \mathbb{Z}\}$. 
For any $\lambda \in k^\times$, we write $\lambda\mathcal{D}_k=\mathfrak{b}\mathfrak{a}^{-1}$ 
with integral ideals $\mathfrak{a}$ and $\mathfrak{b}$ which are coprime to each other. 
Then the \textit{Gauss sum} $C(\lambda)$ is defined by 
\begin{equation}\label{eq:GAUSS-SUM-DEF}
C(\lambda)=\sum_{\mu \in \mathcal{O}_k/\mathfrak{a}} \exp\left(2\pi i \trace_{k/\mathbb{Q}}(\mu^2 \lambda) \right).
\end{equation}

Let $\{1,\omega,\theta\}$ be a $\mathbb{Z}$-basis of $\mathcal{O}_k$ such that the multiplication 
is given by \eqref{eq:CUBIC-RING-STR} and put $F_k(u)=au_1^3+bu_1^2u_2+cu_1u_2^2+du_2^3$. 
We may assume $F_k(u) \equiv u_1u_2^2 \pmod{2}$ or $F_k(u) \equiv u_1^3 \pmod{2}$ if $2$ is ramified in $k/\mathbb{Q}$. 
Then the ramified prime ideal dividing $2$ is given by $\mathfrak{l}_1=[2,\omega,\theta]$ and 
satisfies $\mathfrak{l}_1^2=[2,2\omega,\theta]$. 
If $2\mathcal{O}_k=\mathfrak{l}_1^2\mathfrak{l}_2$ is the prime ideal decomposition, then 
$\mathfrak{l}_2=[2,\omega,\theta+1]$. 
\begin{lem}\label{lem:GAUSS-SUM-DEN4}
We assume that $2$ is ramified in $k/\mathbb{Q}$. 
Let $\{\eta_0,\eta_1,\eta_2\}$ be the dual basis of $\{1,\omega,\theta\}$ and put $\eta=\eta_2-\eta_1$. 
Let $\lambda \in \mathcal{O}_k$ be odd with $\lambda\equiv 1 \pmod{\mathfrak{l}_1^2}$ 
and write $\lambda=1+2r+2s\omega+t\theta$, $r,s,t \in \mathbb{Z}$. Then
\begin{align*}
C(-\eta\lambda/4)
 &=8 \left(1+i^{a+b+2ar-adt} \right)\left(1+ i^{c+d+2cr-c^2t-cdt} \right) \\
 &\quad +8i^{2s-t}\left(1-i^{a+b+2ar-adt} \right)\left(1-(-1)^{ct} i^{c+d+2cr-c^2t-cdt} \right).
\end{align*}
\end{lem}
\Proof
We write $\mu=x+y\omega+z\theta$, $x,y,z \in \{0,1,2,3\}$. 
We have 
$\mu^2 \equiv -acy^2+(by^2+dz^2)\omega-(ay^2+cz^2)\theta  \pmod{4\mathcal{O}_k}$ for $x=0,2$, 
$\mu^2 \equiv 1-acy^2+(by^2+dz^2+2y)\omega -(ay^2+cz^2-2z)\theta \pmod{4\mathcal{O}_k}$ 
for $x=1,3$.
Hence 
\begin{align*}
\mu^2\lambda
 &\equiv -acy^2 +\left[by^2+dz^2-dt(ay^2+cz^2) \right]\omega\\
 &\quad  -\left[(ay^2+cz^2)(1+2r-ct)+acty^2 \right]\theta \pmod{4\mathcal{O}_k}
\end{align*}
for $x=0,2$, and 
\begin{align*}
\mu^2\lambda
 &\equiv 1+2r-acy^2+\left[by^2+dz^2-dt(ay^2+cz^2)+2y+2s \right]\omega \\
 &\quad  -\left[(ay^2+cz^2)(1+2r-ct)+t(acy^2-1)+2ctz-2z \right]\theta \\
 &\quad \pmod{4\mathcal{O}_k} 
\end{align*}
for $x=1,3$. 
Since $\trace_{k/\mathbb{Q}}(\eta)=0, \trace_{k/\mathbb{Q}}(\eta\omega)=-1, \trace_{k/\mathbb{Q}}(\eta\theta)=1$, 
\begin{align*}
\lefteqn{\trace_{k/\mathbb{Q}}(-\eta\mu^2\lambda/4)} \\
 &\equiv \frac{1}{4}\left[(a+b+2ar-adt)y^2 + (c+d+2cr-c^2t-cdt)z^2 \right] \pmod{\mathbb{Z}}
\end{align*}
for $x=0,2$, and 
\begin{align*}
\lefteqn{\trace_{k/\mathbb{Q}}(-\eta\mu^2\lambda/4)} \\
 &\equiv \frac{1}{4}(2s-t) + \frac{1}{4}\left[(a+b+2ar- adt)y^2 + 2y \right] \\ 
 &\quad + \frac{1}{4}\left[(c+d+2cr - cdt -c^2t)z^2  + 2ctz - 2z   \right]  \pmod{\mathbb{Z}}
\end{align*}
for $x=1,3$.  So we obtain 
\begin{align*}
\lefteqn{C(-\eta\lambda/4)}\\
 &=2\sum_{y=0}^3 e^{\frac{\pi i}{2}(a+b+2ar-adt)y^2} \sum_{z=0}^3 e^{\frac{\pi i}{2}(c+d+2cr-c^2t-cdt)z^2} \\
 &\quad+2i^{2s-t} \sum_{y=0}^3 (-1)^y e^{\frac{\pi i}{2}(a+b+2ar- adt)y^2} 
 \sum_{z=0}^3 (-1)^{ctz-z} e^{\frac{\pi i}{2}(c+d+2cr -c^2t - cdt)z^2}  \\
 &=8 \left(1+i^{a+b+2ar-adt} \right)\left(1+ i^{c+d+2cr-c^2t-cdt} \right) \\
 &\quad +8i^{2s-t}\left(1-i^{a+b+2ar-adt} \right)\left(1-(-1)^{ct} i^{c+d+2cr-c^2t-cdt} \right).
\end{align*}
\qed

\begin{lem}\label{lem:RECIPROCITYLAW}
Let $\alpha$, $\beta$ be two odd elements in $\mathcal{O}_k$ which are relatively prime to each other. 
We assume that $\alpha$ satisfies one of the following three conditions \rm{(P1)}, \rm{(P2)} and \rm{(P3)}:
\begin{itemize}
 \item[\rm{(P1)}] $\alpha$ is primary, i.e. $\alpha \equiv \xi^2 \pmod{4\mathcal{O}_k}$ for some $\xi \in \mathcal{O}_k$. 
 \item[\rm{(P2)}] $\alpha \equiv \xi^2 \pmod{\mathfrak{l}_1^4\mathfrak{l}_2^2}$ for some $\xi \in \mathcal{O}_k$,  if $2\mathcal{O}_k=\mathfrak{l}_1^2\mathfrak{l}_2$. 
 \item[\rm{(P3)}] $\alpha \equiv \xi^2 \pmod{\mathfrak{l}_1^4}$ for some $\xi \in \mathcal{O}_k$,  if $2\mathcal{O}_k=\mathfrak{l}_1^3$. 
\end{itemize}
Further we assume $\beta\equiv 1 \pmod{\mathfrak{l}_1^2}$ if $\alpha$ satisfies the conditions \rm{(P2)} or \rm{(P3)}. 
Then the \textit{quadratic reciprocity law} 
\begin{equation}\label{eq:QUADRATIC-RECIPROCITY}
\legend{\alpha}{\beta}\legend{\beta}{\alpha}=(-1)^{\nu(\alpha,\beta)}
\end{equation}
holds. Here $\nu(\alpha,\beta)=\sum_{j=1}^{r_1} (\sgn \alpha^{(j)}-1)/2 \cdot (\sgn \beta^{(j)}-1)/2$ and 
$\alpha^{(j)}, \beta^{(j)}$ ($j=1,\dots,r_1$) are real conjugates of $\alpha, \beta$.
\end{lem}
\Proof
If $\alpha$ is primary, \eqref{eq:QUADRATIC-RECIPROCITY} is \cite[Theorem 165]{HE}. 
We assume that $\alpha$ satisfies (P2) or (P3) and $\beta \equiv 1 \pmod{\mathfrak{l}_1^2}$.  
Replacing $F_k$ by $\gamma_2 F_k$ for a suitable $\gamma_2 \in \Gamma_2$, 
we may assume that $F_k(1,1)$ is coprime to the odd integer $m=N_{k/\mathbb{Q}}(\alpha\beta)$. 
We put $\xi=-\omega/a$ and $\Delta=3a\xi^2+2b\xi+c$. 
Then $\{-\theta/\Delta,(b-\omega)/\Delta,a/\Delta\}$ is the dual basis of 
$\{1,\xi,\xi^2\}$ with respect to the trace pairing as we saw in \S~4. 
This implies that the dual basis $\{\eta_0,\eta_1,\eta_2\}$ of $\{1,\omega,\theta\}$ is given by 
\[ \eta_0=\Delta^{-1}(a\xi^2+b\xi), \quad \eta_1=-\Delta^{-1}\xi, \quad \eta_2=-\Delta^{-1}.\]
Hence we have $\mathcal{D}_k^{-1}=[\eta_0,\eta_1,\eta_2]=\Delta^{-1}[1,\xi,a\xi^2]$. 
It is easy to see that $\mathfrak{a}=[a,\omega-b,\theta]$ is an ideal of $\mathcal{O}_k$ 
satisfying $\mathfrak{a}[1,\xi,a\xi^2]=\mathcal{O}_k$. 
This proves $\mathcal{D}_k=\Delta \mathfrak{a}$. 
We put $\eta=\eta_2-\eta_1 \in \mathcal{D}_k^{-1}$ and $\mathfrak{b}=\eta \mathcal{D}_k$. 
Then we have $\eta=-(\omega+a)/(a\Delta)$ and $\mathfrak{b}=[a+b+c+d,\omega+a,\theta-d]$. 
Since $N(\mathfrak{b})=|a+b+c+d|$ is odd,  $\mathfrak{b}$ is an odd ideal. 
Further $\mathfrak{b}$ is coprime to $\alpha\beta$ since $F_k(1,1)=a+b+c+d$ is coprime to $m$. 
By the equation (206) in \cite{HE} we have 
\begin{equation}\label{eq:RECIPROCITY-GAUSS-SUM}
\legend{\beta}{\alpha}\legend{\alpha}{\beta}
=(-1)^{\nu(\alpha,\beta)} \, \frac{C(-\eta\alpha/4)C(-\eta\beta/4)}{C(-\eta/4)C(-\eta\alpha\beta/4)}.
\end{equation}
We compute the four Gauss sums on the right hand side. 
By assumption, there exits an element $\xi \in \mathcal{O}_k$ such that 
$\alpha \equiv \xi^2 \pmod{\mathfrak{l}_1^2\mathfrak{l}_2^2}$ if $(2)=\mathfrak{l}_1^2\mathfrak{l}_2$ 
and $\alpha \equiv \xi^2 \pmod{\mathfrak{l}_1^4}$ if $(2)=\mathfrak{l}_1^3$. 
We take an element $\xi' \in \mathcal{O}_k$ such that $\xi'\xi \equiv 1 \pmod{4\mathcal{O}_k}$ 
and put $\alpha_1=\xi'^2 \alpha$. Then we have $C(-\eta \alpha/4)=C(-\eta \alpha_1/4)$ and 
$C(-\eta \alpha\beta/4)=C(-\eta \alpha_1\beta/4)$ by a basic property of Gauss sums. 

We first assume $(2)=\mathfrak{l}_1^2\mathfrak{l}_2$. 
Since $\mathfrak{l}_2=[2,\omega,\theta+1]$, $\mathfrak{l}_1^2\mathfrak{l}_2^2=2\mathfrak{l}_2$ 
and $\alpha_1 \equiv 1 \pmod{\mathfrak{l}_1^2\mathfrak{l}_2^2}$,  
we can write $\alpha_1=1+4r+2s\omega+2t(\theta+1)$, $r,s,t \in \mathbb{Z}$. 
Since $\beta \equiv 1 \pmod{\mathfrak{l}_1^2}$ and $\beta \equiv 1 \pmod{\mathfrak{l}_2}$, 
we have $\beta \equiv 1 \pmod{2\mathcal{O}_k}$, 
hence we can write $\beta=1+2u+2v\omega+2w\theta$, $u,v,w \in \mathbb{Z}$. 
Then we have
\[ \alpha_1\beta \equiv 1 + 2(t+u)+2(s+v)\omega+2(t+w)\theta \pmod{4\mathcal{O}_k}.\]
By Lemma \ref{lem:GAUSS-SUM-DEN4}, we have 
\begin{align*}
C(-\eta/4)
 &=16(1+i^{a+b+c+d}), \\
C(-\eta\alpha_1/4)
 &=\left\{
\begin{array}{ll}
16(1+i^{a+b+c+d}) , &\quad s \equiv t \pmod{2}, \\
16i^{a+b}(1+i^{a+b+c+d}) , &\quad s \not\equiv t \pmod{2},
\end{array}
\right. \\
C(-\eta\beta/4)
 &=\left\{
\begin{array}{ll}
16\left(1+(-1)^{u-w}i^{a+b+c+d} \right) , &\quad  v \equiv w \pmod{2}, \\
16i^{a+b} \left(1+(-1)^{u-w}i^{a+b+c+d} \right) , &\quad  v \not \equiv w \pmod{2},
\end{array}
\right. 
\end{align*}
and 
\begin{align*}
\lefteqn{C(-\eta\alpha_1\beta/4)} \\
 &=\left\{
\begin{array}{ll}
16\left(1+(-1)^{u-w}i^{a+b+c+d} \right) , &\quad  s+v \equiv t+w \pmod{2}, \\
16i^{a+b}\left(1+(-1)^{u-w}i^{a+b+c+d} \right) , &\quad  s+v \not\equiv t+w \pmod{2}.
\end{array}
\right.
\end{align*}
We note that the Gauss sums are nonzero since $a+b+c+d$ is odd. 
These equations imply
\[ \frac{C(-\eta\alpha_1/4)}{C(-\eta/4)}=\frac{C(-\eta\alpha_1\beta/4)}{C(-\eta\beta/4)}
=\left\{
\begin{array}{ll}
1 ,   &\quad  s \equiv t \pmod{2}, \\
i^{a+b} , &\quad s \not\equiv t \pmod{2}.
\end{array}
\right.
\]
The reciprocity law \eqref{eq:QUADRATIC-RECIPROCITY} now follows form \eqref{eq:RECIPROCITY-GAUSS-SUM} 
when $(2)=\mathfrak{l}_1^2\mathfrak{l}_2$.

We next assume $(2)=\mathfrak{l}_1^3$. 
Since $\alpha_1 \equiv 1 \pmod{\mathfrak{l}_1^4}$ and $\mathfrak{l}_1^4=2\mathfrak{l}_1$, 
we can write $\alpha_1=1+4r+2s\omega+2t\theta$, $r,s,t \in \mathbb{Z}$. 
Since $\beta \equiv 1 \pmod{\mathfrak{l}_1^2}$, we can write $\beta=1+2u+2v\omega+w\theta$, $u,v,w \in \mathbb{Z}$. 
Then we have
\[ \alpha_1\beta \equiv 1 + 2u+2(s+v)\omega+(2t+w)\theta \pmod{4\mathcal{O}_k}.\]
By Lemma \ref{lem:GAUSS-SUM-DEN4}, we have $C(-\eta/4)=16(1+i^{a+b+c+d})$, 
\begin{align*}
C(-\eta\alpha_1/4)
 &=\left\{
\begin{array}{ll}
16(1+i^{a+b+c+d}) , &\quad s \equiv t \pmod{2}, \\
16i^{c+d}(1+i^{a+b+c+d}) , &\quad s \not\equiv t \pmod{2}.
\end{array}
\right. 
\end{align*}
If $w=2w_1$ is even, then 
\begin{align*}
C(-\eta\beta/4)
 &=\left\{
\begin{array}{ll}
16\left(1+(-1)^u i^{a+b+c+d} \right) , &\quad  v \equiv w_1 \pmod{2}, \\
16i^{c+d} \left(1+(-1)^u i^{a+b+c+d} \right) , &\quad  v \not \equiv w_1 \pmod{2},
\end{array}
\right. 
\end{align*}
\begin{align*}
\lefteqn{C(-\eta\alpha_1\beta/4)} \\
 &=\left\{
\begin{array}{ll}
16\left(1+(-1)^u i^{a+b+c+d} \right) , &\quad s+v \equiv t+w_1 \pmod{2}, \\
16 i^{c+d}\left(1+(-1)^u i^{a+b+c+d} \right), &\quad s+v \not\equiv t+w_1 \pmod{2}.
\end{array}
\right.
\end{align*}
We put $\varepsilon=(-1)^{(c+d)/2}$. If $w=2w_1-1$ is odd, then 
\begin{align*}
\lefteqn{C(-\eta\beta/4)}\\
 &=\left\{
\begin{array}{ll}
8\left(1 + (-1)^u i^{a+b+c}\right) \left[1 + \varepsilon+(1-\varepsilon)i\right] , &  v \equiv w_1 \pmod{2}, \\
8\left(1 + (-1)^u i^{a+b+c}\right) \left[1 + \varepsilon-(1-\varepsilon)i\right] , &  v \not \equiv w_1 \pmod{2},
\end{array}
\right. \\
\lefteqn{C(-\eta\alpha_1\beta/4)} \\
 &=\left\{
\begin{array}{ll}
8\left(1 + (-1)^u i^{a+b+c}\right) \left[1 + \varepsilon+(1-\varepsilon)i\right] , &  s+v \equiv t+w_1 \pmod{2}, \\
8\left(1 + (-1)^u i^{a+b+c}\right) \left[1 + \varepsilon-(1-\varepsilon)i\right] , &  s+v \not \equiv t+w_1 \pmod{2}.
\end{array}
\right.
\end{align*}
These equations imply
\[ \frac{C(-\eta\alpha_1/4)}{C(-\eta/4)}=\frac{C(-\eta\alpha_1\beta/4)}{C(-\eta\beta/4)}
=\left\{
\begin{array}{ll}
1 ,   &\quad  s \equiv t \pmod{2}, \\
i^{c+d} , &\quad s \not\equiv t \pmod{2}.
\end{array}
\right.
\]
The reciprocity law \eqref{eq:QUADRATIC-RECIPROCITY} now follows form \eqref{eq:RECIPROCITY-GAUSS-SUM} 
when $(2)=\mathfrak{l}_1^3$. 
\qed


Let $K$ be an $S_4$-quartic field whose cubic resolvent field is $k$. 
Let $k_6$ be as in before. We first assume that $K$ is totally real if $\Disc(K)>0$. 
We have $N(\Disc(k_6/k))=f^2$. We assume that $f$ is square free. 
We denote by $\mathcal{O}$ $(\subset \mathcal{O}_k)$ the unique cubic resolvent ring of $\mathcal{O}_K$ and $\mathfrak{f}$ the 
conductor of $\mathcal{O}$. 
By Proposition \ref{prop:MAXIMAL-ORDER-CONDUCTOR}, $\Disc(k_6/k)=\mathfrak{f}$. 
We denote by $\mathfrak{f}_\mathrm{odd}$ the odd part of $\mathfrak{f}$. 
We take an element $\alpha \in \mathcal{O}_k$ such that $k_6=k(\sqrt{\alpha})$. 
Since $f$ is square free, we may assume that $\alpha$ is odd and satisfies one of the conditions (P1), (P2) and (P3) 
in Lemma \ref{lem:RECIPROCITYLAW} by ramification theory in relative quadratic extensions. 
Then $(\alpha)=\mathfrak{f}_\mathrm{odd}\mathfrak{a}^2$ for some integral ideal $\mathfrak{a}$ of $\mathcal{O}_k$. 
We may also assume that $\mathfrak{a}$ is relatively prime to $2\mathfrak{f}$. 
By an elementary argument in Galois theory, $N_{k/\mathbb{Q}}(\alpha)$ is square of a rational integer. 
Since $\alpha$ satisfies (P1), (P2), or (P3), we see that $N_{k/\mathbb{Q}}(\alpha)\equiv 1 \pmod{4}$. 
Then $(\alpha)=\mathfrak{f}_\mathrm{odd}\mathfrak{a}^2$ implies $N_{k/\mathbb{Q}}(\alpha)=f^2 N(\mathfrak{a})^2$. 
Hence $\alpha$ is totally positive if $\Disc(K)<0$. 
If $\Disc(K)>0$, then $K$ is totally real by assumption, hence 
$\alpha$ is totally positive. 

By class field theory, the quadratic extension $k_6/k$ corresponds to a subgroup $H$ of $I_k(\mathfrak{f})$ 
such that $H \supset P_{k,1}(\mathfrak{f})$ and $(I_k(\mathfrak{f}):H)=2$. 
We denote by $\chi$ the unique nontrivial character of $I_k(\mathfrak{f})$ such that $\ker \chi=H$. 
If $\mathfrak{q}$ is a prime ideal of $\mathcal{O}_k$ not dividing $\mathfrak{f}$, 
then $\chi(\mathfrak{q})=1$ if and only if $\mathfrak{q}$ splits in $k_6/k$. 
Further if $\mathfrak{q}$ is odd and $\mathfrak{q}\nmid (\alpha)$, then 
$\chi(\mathfrak{q})=\legend{\alpha}{\mathfrak{q}}$. 
Hence we have $\chi(\mathfrak{b})=\legend{\alpha}{\mathfrak{b}}$ for 
any odd integral ideal $\mathfrak{b}$ which is relatively prime to $(\alpha)$. 
We write simply $\chi((\gamma))=\chi(\gamma)$ for $(\gamma) \in P_k(\mathfrak{f})$. 
We now show that $H$ contains $\tilde{P}_\mathcal{O} I_k(\mathfrak{f})^2$, where $\tilde{P}_{\mathcal{O}}$ is defined by \eqref{eq:PO-POTILDE}.
Any element of $\tilde{P}_{\mathcal{O}} I_k(\mathfrak{f})^2$ is written as 
$(\beta) \mathfrak{c}^2$, where $\beta \in \mathcal{O}$ is relatively prime to $\mathfrak{f}$ and 
$\mathfrak{c} \in I_k(\mathfrak{f})$. Then we have $\chi((\beta)\mathfrak{c}^2)=\chi(\beta)$. 
We first assume that $\beta$ is odd and relatively prime to $(\alpha)$. 
By Proposition \ref{prop:CUNDUCTOR3}, we have $\mathcal{O}=\mathbb{Z}+\mathfrak{f}$, hence 
$\beta\equiv b \pmod{\mathfrak{f}}$ for some $b \in \mathbb{Z}$. 
Since $\alpha$ is totally positive, 
the reciprocity law \eqref{eq:QUADRATIC-RECIPROCITY} implies
\begin{equation}\label{eq:KERNEL-CHI1}
\chi(\beta)=\legend{\alpha}{\beta}=\legend{\beta}{\alpha}
=\legend{\beta}{\mathfrak{f}_\mathrm{odd}\mathfrak{a}^2}
=\legend{\beta}{\mathfrak{f}_\mathrm{odd}}=\legend{b}{\mathfrak{f}_\mathrm{odd}}.
\end{equation}
Let $p\neq 2$ be a prime divisor of $f$ and denote by $\mathfrak{f}_p$ the $p$-part of $\mathfrak{f}$. 
Since $p|f$, $p$ is neither of type $3$ nor $1^3$. 
If $p$ is of type $111$ or $11^2$ in $k$, $\mathfrak{f}_p=\mathfrak{p}_1\mathfrak{p}_2$, 
$\mathfrak{p}_1\neq \mathfrak{p}_2$ and $N(\mathfrak{p}_1)=N(\mathfrak{p}_2)=p$. 
Then we have 
\[ \legend{b}{\mathfrak{f}_p}=\legend{b}{\mathfrak{p}_1}\legend{b}{\mathfrak{p}_2}=\legend{b}{p}^2=1.\]
If $p$ is of type $12$, we write $(p)=\mathfrak{p}_1\mathfrak{p}_2$, $N(\mathfrak{p}_i)=p^i$ $(i=1,2)$. 
Then $\mathfrak{f}_p=\mathfrak{p}_2$. 
Since $(\mathcal{O}_k/\mathfrak{p}_2)^\times \cong \mathbb{F}_{p^2}^\times$ is a cyclic group of order $p^2-1$, 
$b\in \mathbb{Z}$ is a  quadratic residue modulo $\mathfrak{p}_2$. 
Therefore $\legend{b}{\mathfrak{f}_p}=\legend{b}{\mathfrak{p}_2}=1$. 
So we have proved $\legend{b}{\mathfrak{f}_p}=1$ for all $p|f$, $p\neq 2$. 
Hence $\chi(\beta)=\legend{b}{\mathfrak{f}_\mathrm{odd}}=1$. 
We next assume that $\beta$ is odd but is not relatively prime to $(\alpha)$. 
We take an odd integral ideal $\mathfrak{a}_1$ which belongs to the same ideal class of $\mathfrak{a}$ 
and is relatively prime to $(\beta)$. 
Then $\mathfrak{a}_1=\lambda \mathfrak{a}$ for some $\lambda \in k^\times$. 
If we put $\alpha_1=\alpha \lambda^2$, then $(\alpha_1)=\mathfrak{f}_\mathrm{odd}\mathfrak{a}_1^2$, hence $\alpha_1 \in \mathcal{O}_k$ and 
$\alpha_1$ is odd. Since $k_6=k(\sqrt{\alpha_1})$, $\alpha_1$ satisfies the same condition as $\alpha$. 
Since $\beta$ is odd and relatively prime to $\mathfrak{f}_\mathrm{odd}\mathfrak{a}_1^2=(\alpha_1)$, 
we have $\chi(\beta)=\legend{\alpha_1}{\beta}=1$ by the result of the first case. 
We finally assume that $\beta$ is not odd. We take an odd integral ideal $\mathfrak{b}$ 
which belongs to the same ray ideal class of $(\beta)$ in $\Cl_k(\mathfrak{f})$. 
Then $(\beta)=(\lambda_1/\lambda_2) \mathfrak{b}$ for some 
$\lambda_1, \lambda_2 \in \mathcal{O}_k$ which are relatively prime to $\mathfrak{f}$ and 
$\lambda_1\equiv \lambda_2 \pmod{\mathfrak{f}}$. 
If put $\xi=\beta\lambda_2/\lambda_1$, then $(\xi)=\mathfrak{b}$, 
hence $\xi \in \mathcal{O}_k$ is odd and relatively prime to $\mathfrak{f}$. 
We have 
\[ \lambda_1(\xi-\beta)=\beta(\lambda_2-\lambda_1)\equiv 0 \pmod{\mathfrak{f}}.\]
Since $\lambda_1$ is relatively prime to $\mathfrak{f}$, 
we have $\xi-\beta \in \mathfrak{f}$, hence $\xi \in \mathfrak{f}+\mathcal{O}=\mathcal{O}$. 
So we have $\chi(\xi)=1$ by the results of the second case. 
Since $(\lambda_1/\lambda_2) \in P_{k,1}(\mathfrak{f})$ and $P_{k,1}(\mathfrak{f}) \subset \ker \chi$ 
by class field theory, we have $\chi(\beta)=\chi(\xi)\chi(\lambda_1/\lambda_2)=1$. 
Thus we have shown that $H$ contains $\tilde{P}_\mathcal{O}I_k(\mathfrak{f})^2$. 

Conversely, given a subgroup $H$ of $I_k(\mathfrak{f})$ such that 
$(I_k(\mathfrak{f}):H)=2$ and $H\supset \tilde{P}_\mathcal{O}I_k(\mathfrak{f})^2$, 
let $k(\sqrt{\alpha})$ be the quadratic extension of $k$ corresponding to $H$ by class field theory. 
Here $\alpha$ is a totally positive element of $k^\times \smallsetminus (k^\times)^2$. 
Let $\chi$ be the character of $I_k(\mathfrak{f})$ such that $\ker \chi=H$. 
We denote by $\mathfrak{g}$ the conductor of $H$. 
Since $\mathfrak{g}|\mathfrak{f}$, we may assume that 
$\alpha \in \mathcal{O}_k$ is odd and satisfies one of (P1), (P2) and (P3). 
Hence $(\alpha)=\mathfrak{g}_\mathrm{odd} \mathfrak{a}^2$ for some 
odd integral ideal $\mathfrak{a}$. Multiplying $\alpha$ by square of an appropriate element, 
we may assume that $\mathfrak{a}$ is relatively prime to $2f$. 
For any prime number $p$, we denote by $\mathfrak{g}_p$ the $p$-part of $\mathfrak{g}$. 
We show that $\mathfrak{g}_p$ equals $\mathfrak{f}_p$ or $(1)$ for all $p|f$. 
Suppose this is not the case. Then there exists a prime divisor $p$ of $f$ 
and a prime ideal $\mathfrak{p}$ such that $p|f$, $N(\mathfrak{p})=p$ and $\mathfrak{g}_p=\mathfrak{p}$. 
If $p=2$, then we must have $\mathfrak{p}^2|\mathfrak{g}=\Disc(k(\sqrt{\alpha})/k)$, this contradicts $\mathfrak{g}_p=\mathfrak{p}$. 
Hence $p \neq 2$. 
We take a rational integer $b_p$ which is a quadratic nonresidue modulo $p$. 
For each prime number $l|f$, $l \neq 2, p$, we take a rational integer $b_l$ 
which is a quadratic residue modulo $l$. 
By Chinese remainder theorem, we can take a positive rational integer $b$ such that 
$b \equiv b_p \pmod{p}$, $b \equiv b_l \pmod{l}$ for all $l|f$ ($l \neq 2, p$) 
and $b \equiv 1 \pmod{2N(\mathfrak{a})}$. 
Let $\mathfrak{q}$ be an odd prime ideal 
which belongs to the ray class of $(b)$ in $\Cl_k(\mathfrak{f})$ and is relatively prime to $\alpha$. 
Then $\mathfrak{q}=(b)(\xi)$ for some $(\xi) \in P_{k,1}(\mathfrak{f})$. 
Since $P_{k,1}(\mathfrak{f}) \subset \ker \chi=H$, we have $\chi(\mathfrak{q})=\chi(b)$. 
The reciprocity law \eqref{eq:QUADRATIC-RECIPROCITY} implies 
$\chi(b)=\legend{\alpha}{b}=\legend{b}{\alpha}=\legend{b}{\mathfrak{g}_\mathrm{odd}}$. 
Since $\legend{b}{\mathfrak{g}_p}=\legend{b}{p}=-1$ and 
$\legend{b}{\mathfrak{g}_l}=1$ for $l|f$ ($l \neq 2,p$), we have 
$\legend{b}{\mathfrak{g}_\mathrm{odd}}=-1$. So we have $\chi(\mathfrak{q})=-1$, hence $\mathfrak{q}$ remains prime in $k(\sqrt{\alpha})/k$. 
This contradicts the fact that $\mathfrak{q} \in (b) P_{k,1}(\mathfrak{f}) \subset \tilde{P}_\mathcal{O}I_k(\mathfrak{f})^2 \subset H$. 
Therefore we have shown that $\mathfrak{g}_p$ equals $\mathfrak{f}_p$ or $(1)$ for all $p|f$. 
Hence $N(\mathfrak{g})=g^2$, $g|f$ and $N_{k/\mathbb{Q}}(\alpha)=\pm g^2 N(\mathfrak{a})^2$. 
We have $N_{k/\mathbb{Q}}(\alpha) \equiv 1 \pmod{4}$ as we have seen before. 
Hence $N_{k/\mathbb{Q}}(\alpha)=g^2 N(\mathfrak{a})^2$. This implies that 
$k(\sqrt{\alpha})$ can not be obtained by adjoining square root of a rational integer. 
So there exists an $S_4$-quartic field $K$ such that $k(\sqrt{\alpha})$ is the sextic filed $k_6$ for $K$. 
Then $\Disc(k_6/k)=\mathfrak{g}$ and $\Disc(K)=\Disc(k)g^2$. 
By Proposition \ref{prop:MAXIMAL-ORDER-CONDUCTOR}, the conductor of the unique cubic resolvent ring 
of $\mathcal{O}_K$ is $\mathfrak{g}$. 
By \eqref{eq:CLO-RAYCLASSGROUP}, we have $\Cl_\mathcal{O}/\Cl_\mathcal{O}^2 \cong I_k(\mathfrak{f})/\tilde{P}_\mathcal{O}I_k(\mathfrak{f})^2$. 
Hence the number of such subgroups $H$ is equal to 
$|\Cl_\mathcal{O}/\Cl_\mathcal{O}^2|-1=|\Cl_\mathcal{O}^{(2)}|-1$. 
Thus we have proved the former half of the following proposition. 
The latter half is proved similarly. 
\begin{prop}\label{prop:NUMBER-QUARTIC-FIELD}
Let $k$ be a non-Galois cubic field and $\mathcal{O}$ be an order of $k$ 
such that the index $f=(\mathcal{O}_k:\mathcal{O})$ is square free. 
We denote by $\mathfrak{f}$ the conductor of $\mathcal{O}$. 
For each divisor $g$ of $f$, we put $\mathfrak{g}=\prod_{p|g} \mathfrak{f}_p$ and 
$R_g=\mathbb{Z}+\mathfrak{g}$. We denote by $\mathcal{K}_k(\mathfrak{g})$ the set of isomorphism classes of quartic fields 
$K$ satisfying the following conditions:
\begin{itemize}
 \item[(a)] The normal closure $\tilde{K}$ of $K$ over $\mathbb{Q}$ has Galois group $S_4$ and contains $k$. 
 \item[(b)] The unique cubic resolvent ring of the maximal order $\mathcal{O}_K$ is isomorphic to $R_g$. 
 \item[(c)] $K$ is totally real if $\Disc(k)>0$. 
\end{itemize}
Then we have 
\[ \sum_{g|f} |\mathcal{K}_k(\mathfrak{g})|=|\Cl_\mathcal{O}^{(2)}|-1.\] 
If $\Disc(k)>0$, we denote by $\mathcal{K}_k(\mathfrak{g}\mathfrak{f}_\infty)$ the set of isomorphism classes of quartic fields 
$K$ satisfying the conditions \rm{(a)} and \rm{(b)} above. Then we have 
\[ \sum_{g|f} |\mathcal{K}_k(\mathfrak{g}\mathfrak{f}_\infty)|=|\Cl_{\mathcal{O},+}^{(2)}|-1.\] 
\end{prop}

\subsection{Number of quartic rings with a fixed cubic resolvent ring}
Let $k$, $\mathcal{O}$, $f$ and $\mathfrak{f}$ be as in Proposition \ref{prop:NUMBER-QUARTIC-FIELD}. 
We study the number of quartic rings whose cubic resolvent rings are isomorphic to $\mathcal{O}$. 

Let $K$ be an $S_4$-quartic field. We assume that an order $Q$ of $K$ has a cubic resolvent ring isomorphic to $\mathcal{O}$. 
Let $k_6$ be the non-Galois sextic subfield of the Galois closure $\tilde{K}$ as before. 
We denote by $R_K\subset \mathcal{O}_k$ the unique cubic resolvent ring of $\mathcal{O}_K$. 
If $R\subset \mathcal{O}_k$ is a cubic resolvent ring of $Q$, then obviously $R \subset R_K$. 
Since $k$ is a non-Galois cubic field and $R$ is isomorphic to $\mathcal{O}$, 
we have $\mathcal{O}=R \subset R_K \subset \mathcal{O}_k$. 
Since $\Disc(Q)=\Disc(R)$ and $\Disc(\mathcal{O}_K)=\Disc(R_K)$, 
we have $(R_K:\mathcal{O})=(\mathcal{O}_K:Q)$. 
We put $g=(\mathcal{O}_k:R_K)$ and $h=(R_K:\mathcal{O})$. 
Then $f=(\mathcal{O}_k:\mathcal{O})=gh$. 
We denote by $\mathfrak{g}$ the conductor of $R_K$. 
By Lemma \ref{lem:CONDUCTOR1}, $\mathfrak{g}|\mathfrak{f}$. 
We write $\mathfrak{f}=\mathfrak{g}\mathfrak{h}$. 
Since $f$ is square free, it follows from Proposition \ref{prop:CUNDUCTOR3} that 
$\mathcal{O}=\mathbb{Z}+\mathfrak{f}$ and $N(\mathfrak{f})=f^2$. 
Similarly we have $R_K=\mathbb{Z}+\mathfrak{g}$ and $N(\mathfrak{g})=g^2$. 
Since $f$ is square free, $g$ and $h$ are also square free and relatively prime to each other. 
In particular, $Q$ is primitive, hence $R=\mathcal{O}$ is the unique cubic resolvent ring of $Q$. 
We now study the number of such quartic rings $Q$. 

We denote by $a_K(\mathfrak{h})$ the number of quartic rings $Q$ with index $h$ in $\mathcal{O}_K$ such that 
the unique cubic resolvent ring of $Q$ has conductor $\mathfrak{f}=\mathfrak{g}\mathfrak{h}$. 
We have $\mathfrak{g}=\prod_{p|g} \mathfrak{f}_p$ and $\mathfrak{h}=\prod_{p|h} \mathfrak{f}_p$. 
For any order $R$ of $k$, we denote by $\mathfrak{f}_R$ the conductor of $R$. 
\begin{lem}\label{lem:AK-FORMULA}
$a_K(\mathfrak{h})=\prod_{p|h} a_K(\mathfrak{f}_p)$. 
\end{lem}
\Proof
We write $h=p_1 \cdots p_t$, where $p_i$'s are distinct prime numbers. 
Let $Q$ be a quartic ring with index $h$ in $\mathcal{O}_K$ such that 
the unique cubic resolvent ring of $Q$ has conductor $\mathfrak{g}\mathfrak{h}=\mathfrak{f}$, 
Since $h$ is square free, there exists a $\mathbb{Z}$-basis $\{1,\alpha_1,\alpha_2,\alpha_3\}$ of $\mathcal{O}_K$ 
such that $\{1,\alpha_1,\alpha_2,h\alpha_3\}$ is a $\mathbb{Z}$-basis of $Q$ by a theorem of elementary divisors. 
We put $Q_i=Q+p_i\mathcal{O}_K$ ($i=1,\dots,t$). 
Then $Q_i=[1,\alpha_1,\alpha_2,p_i\alpha_3]$, hence $Q_i$ is a quartic ring with index $p_i$ in $\mathcal{O}_K$. 
Since $Q\subset Q_1 \cap \cdots \cap Q_t$ and both of $Q$ and $Q_1 \cap \cdots \cap Q_t$ have index $h$ in $\mathcal{O}_K$, 
we have $Q=Q_1 \cap \cdots \cap Q_t$. 
Conversely, let $Q_i$ be quartic rings with index $p_i$ in $\mathcal{O}_K$ for $i=1,\dots,t$.  
We put $Q=Q_1 \cap \cdots \cap Q_t$. Then $Q$ is a quartic ring with index $h$ in $\mathcal{O}_K$. 
Since $Q+p_i\mathcal{O}_K \subset Q_i$ and both of $Q+p_i\mathcal{O}_K$ and $Q_i$ have 
index $p_i$ in $\mathcal{O}_K$, we have $Q+p_i\mathcal{O}_K=Q_i$. 
So the correspondence $Q\mapsto (Q_1,\dots,Q_t)$ defines a bijection between the set of 
quartic rings $Q$ with index $h$ in $\mathcal{O}_K$ 
and the set of $t$-tuples $(Q_1,\dots,Q_t)$, where $Q_i$'s are quartic rings with index $p_i$ in $\mathcal{O}_K$. 
We denote by $R_K$, $R$ and $R_i$ the unique cubic resolvent rings of $\mathcal{O}_K$, $Q$ and $Q_i$, respectively. 
Then we have $R \subset R_i \subset R_K$ 
$(R_K:R)=(\mathcal{O}_K:Q)=h$ and $(R_K:R_i)=(\mathcal{O}_K:Q_i)=p_i$. 
Hence we have $(R_i:R)=(Q_i:Q)=h/p_i$. 
Since $(\mathcal{O}_k:R_i)=p_ig$ is square free, we have $N(\mathfrak{f}_{R_i})=p_i^2g^2$
by Proposition \ref{prop:CUNDUCTOR3}. 
Since $p_i\nmid (R_i:R)$, Lemma \ref{lem:CONDUCTOR1} implies 
that the $p$-part of $\mathfrak{f}_{R_i}$ equals the $p$-part of $\mathfrak{f}_R$. 
The desired formula for $a_K(\mathfrak{h})$ now follows from the bijection $Q\mapsto (Q_1,\dots,Q_t)$. 
\qed

\begin{lem}\label{lem:PRIME-INDEX-ORDER}
Let $F$ be an \'{e}tale algebra of degree $n$ over $\mathbb{Q}$ and $p$ be a prime number. 
If $\{1,e_1,\dots,e_{n-1}\}$ is a basis of $\mathcal{O}_F/p\mathcal{O}_F$ over $\mathbb{F}_p$, 
then there exists a $\mathbb{Z}$-basis $\{1,\alpha_1,\dots,\alpha_{n-1}\}$ of $\mathcal{O}_F$ 
and a positive rational integer $\lambda$ with $p\nmid \lambda$ such that 
$\pi(\alpha_i)=e_i$ for $1\leq i<n-1$ and $\pi(\alpha_{n-1})=\lambda e_{n-1}$, 
where $\pi:\mathcal{O}_F\rightarrow \mathcal{O}_F/p\mathcal{O}_F$ is the natural homomorphism. 
In particular, when $n=3$ or $4$, if $\{1,e_1,\dots,e_{n-1}\}$ is normalized, then 
we can chose $\alpha_i$'s so that $\{1,\alpha_1,\dots,\alpha_{n-1}\}$ is also normalized. 
\end{lem}
\Proof
We take any $\mathbb{Z}$-basis $\{1,\beta_1,\dots,\beta_{n-1}\}$ of $\mathcal{O}_F$. 
Translating $\beta_i$'s by the appropriate elements of $\mathbb{Z}$,  
we have ${}^t(e_i)=\bar{h}\, {}^t (\pi(\beta_i))$ for some $\bar{h} \in \GL_{n-1}(\mathbb{F}_p)$. 
We take a positive rational integer $\lambda$ such that $\bar{\lambda}=\lambda \bmod{p}=(\det \bar{h})^{-1}$ 
and put $\bar{g}=\diag[1,\dots,1,\bar{\lambda}]\, \bar{h}$. 
Then we have $\bar{g} \in \SL_{n-1}(\mathbb{F}_p)$. Since the natural homomorphism 
$\SL_{n-1}(\mathbb{Z})\rightarrow \SL_{n-1}(\mathbb{F}_p)$ is surjective, there exists an element 
$g \in \SL_{n-1}(\mathbb{Z})$ such that $g \bmod{p}=\bar{g}$. 
If we put ${}^t(\alpha_i)=g\,{}^t(\beta_i)$, then 
$\{1,\alpha_1,\dots,\alpha_{n-1}\}$ is a desired $\mathbb{Z}$-basis of $\mathcal{O}_F$. 
We assume that $n=3$ or $4$ and $\{1,e_1,\dots,e_{n-1}\}$ is normalized. 
Translating $\alpha_i$'s by the appropriate elements of $p\mathbb{Z}$,  
the second statement follows. 
\qed

By Proposition \ref{prop:MAXIMAL-ORDER-CONDUCTOR}, the relative discriminant $\Disc(k_6/k)$ equals 
the conductor $\mathfrak{g}$ of $R_K$. 
We denote by $H=H_K$ the subgroup of $I_k(\mathfrak{g})$ corresponding to the quadratic extension $k_6/k$ by class field theory 
and $\chi=\chi_K$ the character of $I_k(\mathfrak{g})$ such that $\ker \chi=H$. 
%
\begin{prop}\label{prop:AK-FORMULA-P}
$a_K(\mathfrak{f}_p)=1+\chi(\mathfrak{f}_p)$ for $p|h$. 
\end{prop}
\Proof
Let $\pi:\mathcal{O}_K\rightarrow \mathcal{O}_K/p\mathcal{O}_K$ 
be the natural homomorphism. 

Case 1. $p$ is of type $1111$ in $K$. 
We have a ring isomorphism $\varphi:\mathcal{O}_K/p\mathcal{O}_K \cong \mathbb{F}_p^4$. 
We put $e_0=(1,1,1,1)$, $e_1=(0,1,0,0)$, $e_2=(0,0,1,0)$, $e_3=(0,0,0,1)$. 
Then $\{e_i\}$ is a normalized basis of $\mathbb{F}_p^4$. 
By Lemma \ref{lem:PRIME-INDEX-ORDER}, 
there exists a normalized $\mathbb{Z}$-basis $\{\alpha_i\}$ of $\mathcal{O}_K$ and a positive integer $\lambda$ with $p\nmid \lambda$ such that 
$\varphi(\pi(\alpha_i)) = e_i$ ($i=0,1,2$) and $\varphi(\pi(\alpha_3))=\lambda e_3$. 
We take a pair of integral ternary quadratic forms $(A,B)$ such that 
$Q(A,B)=\mathcal{O}_K$ and the multiplication of $\mathcal{O}_K$ is given by \eqref{eq:QUARTIC-RING-STR}. 
Since $e_i^2=e_i$ ($1\leq i \leq 3$) and $e_ie_j=0$ ($1\leq i<j\leq 3$), 
almost all of the fifteen $\SL_2$-invariants $\lambda^{ij}_{k\ell}=\lambda^{ij}_{k\ell}(A,B)$ are congruent to $0$ modulo $p$ 
except that $\lambda^{12}_{13} \equiv -1$,  $\lambda^{12}_{23} \equiv 1$ and $\lambda^{13}_{23} \equiv -\lambda \pmod{p}$.
Since $\lambda^{12}_{13} \equiv -1 \pmod{p}$, we can take an element $\gamma_2 \in \SL_2(\mathbb{Z})$ such that 
$\gamma_2 \equiv \left(\begin{array}{cc}
 b_{12} & -a_{12} \\
 -b_{13} & a_{13} 
\end{array} \right) \pmod{p}$. Then we have 
\begin{align*}
b_{12}A(v)-a_{12}B(v)
 &=\sum_{i\leq j} \lambda^{ij}_{12}v_iv_j \equiv v_1v_3-v_2v_3 \pmod{p}, \\
-b_{13}A(v)+a_{13}B(v)
 &=\sum_{i\leq j} \lambda^{13}_{ij}v_iv_j \equiv v_1v_2-\lambda v_2v_3 \pmod{p}. 
\end{align*}
Replacing $(A,B)$ by $(1_3,\gamma_2)\cdot (A,B)$, we may assume that 
\begin{equation}\label{eq:AB-MODP-CASE1}
A(v)\equiv v_1v_3-v_2v_3 \pmod{p},\quad  B(v)\equiv v_1v_2-\lambda v_2v_3 \pmod{p}.
\end{equation}
Then we have 
\[ F_{(A,B)}(u)=au_1^3+bu_1^2u_2+cu_1u_2^2+du_2^3 \equiv u_1u_2(u_1-\lambda u_2) \pmod{p}. \]
Hence $a\equiv d \equiv 0$, $b\equiv 1$ and $c \equiv -\lambda \pmod{p}$. 
Since the ring structure of $R_K=R(A,B)=[1,\omega,\theta]$ is given by \eqref{eq:CUBIC-RING-STR}, 
we have 
\begin{equation}\label{eq:CUBIC-MODP-CASE1}
\omega^2 \equiv \omega, \quad \theta^2 \equiv \lambda \theta, \quad
\omega\theta \equiv 0 \pmod{pR_K}.
\end{equation}
We put $\mathfrak{p}_1=[p,\omega,\theta-\lambda]$, 
$\mathfrak{p}_2=[p,\omega-1,\theta]$ and $\mathfrak{p}_3=[p,\omega,\theta]$. 
It is easy to see that $\mathfrak{p}_i$'s are $R_K$-ideals such that 
$\mathfrak{p}_2\mathfrak{p}_3=[p,p\omega,\theta]$, 
$\mathfrak{p}_1\mathfrak{p}_3=[p,\omega,p\theta]$, 
$\mathfrak{p}_1\mathfrak{p}_2=[p,p\omega,\theta+\lambda(\omega-1)]$ and 
$\mathfrak{p}_1\mathfrak{p}_2\mathfrak{p}_3=pR_K$. 
If we put $\tilde{\mathfrak{p}}_i=\mathfrak{p}_i\mathcal{O}_k$, then 
$p\mathcal{O}_k=\tilde{\mathfrak{p}}_1\tilde{\mathfrak{p}}_2\tilde{\mathfrak{p}}_3$. 
There exist exactly three cubic rings contained in $R_K$ with index $p$ in $R_K$, namely 
\[ R_1=[1,p\omega,\theta], \quad R_2=[1,\omega,p\theta],\quad R_3=[1, p\omega,\theta+\lambda(\omega-1)].\]
Hence $\mathfrak{f}_i=\prod_{j\neq i} \mathfrak{p}_j$ is the largest $R_K$-ideal contained in $R_i$ for $i=1,2,3$. 
By Lemma \ref{lem:CONDUCTOR2}, the $p$-part of the conductor of $R_i$ equals $\tilde{\mathfrak{f}}_i=\prod_{j\neq i} \tilde{\mathfrak{p}}_j$. 
On the other hand, for each $1\leq i<j\leq 4$, we put 
\[ B_{ij}=\{(b_1,b_2,b_3,b_4) \in \mathbb{F}_p^4\,|\, b_i=b_j\}. \]
Then it is easy to see that the six subrings $B_{ij}$'s are the index $p$ subrings of $\mathbb{F}_p^4$ with unity.  
Then $(\varphi\circ \pi)^{-1}(B_{ij})$'s are the quartic rings with index $p$ in $\mathcal{O}_K$. 
They are given by 
\begin{alignat*}{2}
 Q_1&=[1,\alpha_1,p\alpha_2,\alpha_3],  & \quad  Q_2&=[1,p\alpha_1,\alpha_2,\alpha_3+\lambda\alpha_1],  \\
 Q_3&=[1,\alpha_1,\alpha_2,p\alpha_3],  & \quad  Q_4&=[1,p\alpha_1,\alpha_2+\alpha_1,\alpha_3],  \\
 Q_5&=[1,p\alpha_1,\alpha_2,\alpha_3],  & \quad  Q_6&=[1,\alpha_1,p\alpha_2,\alpha_3+\lambda\alpha_2]. 
\end{alignat*}
We now compute the cubic resolvent rings of $Q_i$'s. 
We put $\delta_{11}=\diag[1,p,1]$, $\delta_{12}=\diag[1,p]$, $\delta_1=(\delta_{11},\delta_{12}^{-1})$ 
and  $(A_1,B_1)=\delta_1 \cdot (A,B)$. 
Then we have $A_1(v)=A(v_1,pv_2,v_3)$ and 
\[ B_1(v) 
=b_{12}v_1v_2+pb_{22}v_2^2+b_{23}v_2v_3
+p^{-1}(b_{11}v_1^2+b_{13}v_1v_3+b_{33}v_3^2).
\]
It follows from \eqref{eq:AB-MODP-CASE1} that $(A_1,B_1)$ is integral. 
We have $Q(A_1,B_1)=Q_1$ and 
$F_{(A_1,B_1)}(u)=p^{-1}F_{(A,B)}(pu_1,u_2)=(\delta_{12}'\cdot F_{(A,B)})(u)$, 
where $\delta_{12}'=\diag[p,1]$. By Proposition \ref{prop:CUBIC-SUBRING}, 
$R(A_1,B_1)=[1,p\omega,\theta]=R_1$. 

We put 
$\delta_{21}=\left(
\begin{array}{ccc}
p & 0 & 0 \\
0 & 1 & 0 \\
\lambda & 0 & 1
\end{array}
\right)$, 
$\delta_2=(\delta_{21},\delta_{12}^{-1})$ 
and  $(A_2,B_2)=\delta_2 \cdot (A,B)$. 
Then we have $A_2(v)=A(pv_1+\lambda v_3,v_2,v_3)$ and 
\begin{align*}
B_2(v)
 &=pb_{11}v_1^2+b_{12}v_1v_2+(2b_{11}\lambda+b_{13})v_1v_3 \\
 &\quad+p^{-1}\left[b_{22}v_2^2+(b_{12}\lambda+b_{23})v_2v_3+(b_{11}\lambda^2+b_{13}\lambda+b_{33})v_3^2 \right].
\end{align*}
It follows from \eqref{eq:AB-MODP-CASE1} that $(A_2,B_2)$ is integral. 
We have $Q(A_2,B_2)=Q_2$ and $F_{(A_2,B_2)}(u)=F_{(A_1,B_1)}(u)$, 
hence $R(A_2,B_2)=R_1$.

We put $\delta_{31}=\diag[1,1,p]$, $\delta_{32}=\diag[p,1]$, $\delta_3=(\delta_{31},\delta_{32}^{-1})$ 
and  $(A_3,B_3)=\delta_3 \cdot (A,B)$. 
Then we have $B_3(v)=B(v_1,v_2,pv_3)$ and 
\[ A_3(v)=a_{13}v_1v_3+a_{23}v_2v_3+pa_{33}v_3^2+p^{-1}(a_{11}v_1^2+a_{12}v_1v_2+a_{22}v_2^2).\]
It follows from \eqref{eq:AB-MODP-CASE1} that $(A_3,B_3)$ is integral. 
We have $Q(A_3,B_3)=Q_3$ and 
$F_{(A_3,B_3)}(u)=p^{-1}F_{(A,B)}(u_1,pu_2)=(\delta_{32}'\cdot F_{(A,B)})(u)$, 
where $\delta_{32}'=\diag[1,p]$. By Proposition \ref{prop:CUBIC-SUBRING}, 
$R(A_3,B_3)=[1,\omega,p\theta]=R_2$. 

We put 
$\delta_{41}=\left(
\begin{array}{ccc}
p & 0 & 0 \\
1 & 1 & 0 \\
0 & 0 & 1
\end{array}
\right)$, $\delta_4=(\delta_{41},\delta_{32}^{-1})$ 
and  $(A_4,B_4)=\delta_4 \cdot (A,B)$. 
Then we have $B_4(v)=B(pv_1+v_2,v_2,v_3)$ and 
\begin{align*}
A_4(v)
 &=pa_{11}v_1^2+(2a_{11}+a_{12})v_1v_2+a_{13})v_1v_3 \\
 &\quad+p^{-1}\left[(a_{11}+a_{12}+a_{22})v_2^2+(a_{13}+a_{23})v_2v_3+a_{33}v_3^2 \right].
\end{align*}
It follows from \eqref{eq:AB-MODP-CASE1} that $(A_4,B_4)$ is integral. 
We have $Q(A_4,B_4)=Q_4$ and $F_{(A_4,B_4)}(u)=F_{(A_3,B_3)}(u)$, 
hence $R(A_4,B_4)=R_2$. 

We take a rational integer $\lambda'$ such that $\lambda\lambda'\equiv 1 \pmod{p}$. 
We put $\delta_{51}=\diag[p,1,1]$, $\delta_{52}=\left(\begin{array}{cc}
 p & \lambda' \\
 0 & 1 
\end{array} \right)$, $\delta_5=(\delta_{51},\delta_{52}^{-1})$ 
and  $(A_5,B_5)=\delta_5 \cdot (A,B)$. 
Then we have $B_5(v)=B(pv_1,v_2,v_3)$ and 
\begin{align*}
A_5(v)&=p(a_{11}-b_{11}\lambda')v_1^2+(a_{12}-b_{12}\lambda')v_1v_2+(a_{13}-b_{13}\lambda')v_1v_3 \\
&\quad +p^{-1}\left[(a_{22}-b_{22}\lambda')v_2^2+(a_{23}-b_{23}\lambda')v_2v_3+(a_{33}-b_{33}\lambda')v_3^2 \right].
\end{align*}
It follows from \eqref{eq:AB-MODP-CASE1} and $\lambda\lambda'\equiv 1 \pmod{p}$ that $(A_5,B_5)$ is integral. 
We have $Q(A_5,B_5)=Q_5$ and 
\[ F_{(A_5,B_5)}(u)=p^{-1}F_{(A,B)}(u_1,\lambda'u_1+pu_2)=(\delta_{52}'\cdot F_{(A,B)})(u), 
\]
where $\delta_{52}'=\left(
\begin{array}{cc}
1 & \lambda' \\
0 & p
\end{array}
\right)$. Since $\left(
\begin{array}{cc}
p & 0 \\
\lambda & 1
\end{array}
\right)(\delta_{52}')^{-1} \in \SL_2(\mathbb{Z})$, Proposition \ref{prop:CUBIC-SUBRING} implies 
$R(A_5,B_5)=[1,\omega+\lambda'\theta,p\theta]=[1,p\omega, \theta+\lambda\omega]=R_3$.

We put 
$\delta_{61}=\left(
\begin{array}{ccc}
1 & 0 & 0 \\
0 & p & 0 \\
0 & \lambda & 1
\end{array}
\right)$, $\delta_6=(\delta_{61},\delta_{52}^{-1})$ 
and  $(A_6,B_6)=\delta_6 \cdot (A,B)$. 
Then we have $B_6(v)=B(v_1,pv_2+\lambda v_3,v_3)$ and 
\begin{align*}
A_6(v)
&=(a_{12}-b_{12}\lambda')v_1v_2+p(a_{22}-b_{22}\lambda')v_2^2\\
&\quad +(2a_{22}\lambda+a_{23}-2b_{22}\lambda\lambda'-b_{23}\lambda')v_2v_3 \\
&\quad +p^{-1}\left[(a_{11}-b_{11}\lambda')v_1^2+(a_{12}\lambda+a_{13}-b_{12}\lambda\lambda'-b_{13}\lambda')v_1v_3 \right]\\
&\qquad +p^{-1}\left[(a_{22}\lambda^2+a_{23}\lambda+a_{33})-\lambda'(b_{22}\lambda^2+b_{23}\lambda+b_{33}) \right]v_3^2.
\end{align*}
It follows from \eqref{eq:AB-MODP-CASE1} and $\lambda\lambda'\equiv 1 \pmod{p}$ that $(A_6,B_6)$ is integral. 
We have $Q(A_6,B_6)=Q_6$ and $F_{(A_6,B_6)}(u)=F_{(A_5,B_5)}(u)$, 
hence $R(A_6,B_6)=R_3$. 
The only ideals $\mathfrak{f}_p$ such that $p\mathcal{O}_k \subset \mathfrak{f}_p$ and $N(\mathfrak{f}_p)=p^2$ are 
$\tilde{\mathfrak{f}}_1$, $\tilde{\mathfrak{f}}_2$ and $\tilde{\mathfrak{f}}_3$  
which are the $p$-part of the conductors of $R_1$, $R_2$ and $R_3$, respectively. 
Since $\tilde{\mathfrak{p}}_i$'s split in $k_6/k$, we have $\chi(\tilde{\mathfrak{p}}_i)=1$. 
So we have proved that $a_K(\tilde{\mathfrak{f}}_i)=2=1+\chi(\tilde{\mathfrak{f}}_i)$ for $i=1,2,3$.

Case 2. $p$ is of type $22$ in $K$. 
We have a ring isomorphism $\varphi:\mathcal{O}_K/p\mathcal{O}_K \cong \mathbb{F}_{p^2}\oplus \mathbb{F}_{p^2}$. 
We write $\mathbb{F}_{p^2}=\mathbb{F}_p(\xi)$. Let $x^2+\bar{r}x+\bar{s}$ be the minimal polynomial of $\xi$ over $\mathbb{F}_p$, 
where $r, s \in \mathbb{Z}$, $\bar{r}=r \bmod{p}$ and $\bar{s}=s \bmod{p}$. 
We put $e_0=(1,1)$, $e_1=(\xi,0)$, $e_2=(0,1)$ and $e_3=(0,\xi)$. 
Then $\{e_i\}$ is a normalized basis of $\mathbb{F}_{p^2}\oplus \mathbb{F}_{p^2}$. 
By Lemma \ref{lem:PRIME-INDEX-ORDER}, 
there exists a normalized $\mathbb{Z}$-basis $\{\alpha_i\}$ of $\mathcal{O}_K$ and a positive integer $\lambda$ with $p\nmid \lambda$ such that 
$\varphi(\pi(\alpha_i)) = e_i$ ($i=0,1,2$) and $\varphi(\pi(\alpha_3))=\lambda e_3$. 
We take a pair of integral ternary quadratic forms $(A,B)$ such that 
$Q(A,B)=\mathcal{O}_K$ and the multiplication of $\mathcal{O}_K$ is given by \eqref{eq:QUARTIC-RING-STR}. 
Since $e_1^2=-se_0-re_1+se_2$, $e_1e_2=e_1e_3=0$, $e_2^2=e_2$, $e_2e_3=e_3$  
and $e_3^2=-se_2-re_3$, 
we have 
$\lambda^{11}_{13} \equiv s$,  $\lambda^{12}_{13} \equiv r$, $\lambda^{13}_{22} \equiv -1$, 
$\lambda^{13}_{23} \equiv r\lambda$, $\lambda^{13}_{33} \equiv -s\lambda^2 \pmod{p}$ 
and remaining ten $\lambda^{ij}_{k\ell}$'s are congruent to $0$ modulo $p$. 
Since $\lambda^{13}_{22} \equiv -1 \pmod{p}$, we can take an element $\gamma_2 \in \SL_2(\mathbb{Z})$ such that 
$\gamma_2 \equiv \left(\begin{array}{cc}
 b_{13} & -a_{13} \\
-b_{22} & a_{22} 
\end{array} \right) \pmod{p}$. 
Then we have 
\begin{align*}
b_{13}A(v)-a_{13}B(v)
 &=\sum_{i\leq j} \lambda^{ij}_{13}v_iv_j \\
 &\equiv sv_1^2+rv_1v_2+v_2^2-r\lambda v_2v_3+s\lambda^2 v_3^2 \pmod{p}, \\
-b_{22}A(v)+a_{22}B(v)
 &=\sum_{i\leq j} \lambda^{22}_{ij}v_iv_j
 \equiv v_1v_3 \pmod{p}. 
\end{align*}
We put $\gamma'_2=\left(
\begin{array}{cc}
1 & -2s\lambda \\
0 & 1
\end{array}
\right) \in \SL_2(\mathbb{Z})$. 
Replacing $(A,B)$ by $(1_3,\gamma'_2\gamma_2)\cdot (A,B)$, we may assume that 
\begin{equation}\label{eq:AB-MODP-CASE2}
\begin{split}
A(v) &\equiv sv_1^2+rv_1v_2-2s\lambda v_1v_3\\
&\quad +v_2^2-r\lambda v_2v_3+s\lambda^2 v_3^2 \pmod{p}, \\
B(v) &\equiv v_1v_3 \pmod{p}. 
\end{split}
\end{equation}
Then $F_{(A,B)}(u)=au_1^3+bu_1^2u_2+cu_1u_2^2+du_2^3$ satisfies 
\[ F_{(A,B)}(u) \equiv u_1u_2\left[\lambda(r^2-4s)u_1-u_2 \right] \pmod{p}. \]
Hence $a\equiv d \equiv 0$, $b\equiv \lambda(r^2-4s)$ and $c \equiv -1 \pmod{p}$. 
Since the ring structure of $R_K=R(A,B)=[1,\omega,\theta]$ is given by \eqref{eq:CUBIC-RING-STR}, 
we have 
\begin{equation}\label{eq:CUBIC-MODP-CASE2}
\omega^2 \equiv \lambda(r^2-4s)\omega, \quad \theta^2 \equiv \theta, \quad
\omega\theta \equiv 0 \pmod{pR_K}.
\end{equation}
We put $\mathfrak{p}_1=[p,\omega,\theta-1]$, 
$\mathfrak{p}_2=[p,\omega-\lambda(r^2-4s),\theta]$ and $\mathfrak{p}_3=[p,\omega,\theta]$. 
It is easy to see that $\mathfrak{p}_i$'s are $R_K$-ideals such that 
$\mathfrak{p}_2\mathfrak{p}_3=[p,p\omega,\theta]$, 
$\mathfrak{p}_1\mathfrak{p}_3=[p,\omega,p\theta]$, 
$\mathfrak{p}_1\mathfrak{p}_2=[p,\omega-\lambda(r^2-4s),p\theta]$ and 
$\mathfrak{p}_1\mathfrak{p}_2\mathfrak{p}_3=pR_K$. 
If we put $\tilde{\mathfrak{p}}_i=\mathfrak{p}_i\mathcal{O}_k$, then 
$p\mathcal{O}_k=\tilde{\mathfrak{p}}_1\tilde{\mathfrak{p}}_2\tilde{\mathfrak{p}}_3$. 
There exist exactly three cubic rings contained in $R_K$ with index $p$ in $R_K$, namely 
\[ R_1=[1,p\omega,\theta], \quad R_2=[1,\omega,p\theta],\quad R_3=[1, \omega-\lambda(r^2-4s),p\theta].\]
Hence $\mathfrak{f}_i=\prod_{j\neq i} \mathfrak{p}_j$ is the largest $R_K$-ideal contained in $R_i$ for $i=1,2,3$. 
By Lemma \ref{lem:CONDUCTOR2}, the $p$-part of the conductor of $R_i$ equals $\tilde{\mathfrak{f}}_i=\prod_{j\neq i} \tilde{\mathfrak{p}}_j$. 
On the other hand, it is easy to see that there exist exactly two index $p$ subrings of $\mathbb{F}_{p^2}\oplus \mathbb{F}_{p^2}$ with unity, 
namely $B_1=\mathbb{F}_p \oplus \mathbb{F}_{p^2}$ and $B_2=\mathbb{F}_{p^2}\oplus \mathbb{F}_p$. 
Then $(\varphi\circ \pi)^{-1}(B_i)$'s are the quartic rings with index $p$ in $\mathcal{O}_K$. 
They are given by 
\[ Q_1=[1, p\alpha_1,\alpha_2,\alpha_3], \quad Q_2=[1,\alpha_1,\alpha_2,p\alpha_3]. \]
We now compute the cubic resolvent rings of $Q_i$'s. 
We put $\delta_{11}=\diag[p,1,1]$, $\delta_{12}=\diag[1,p]$, $\delta_1=(\delta_{11},\delta_{12}^{-1})$ 
and  $(A_1,B_1)=\delta_1 \cdot (A,B)$. 
Then we have $A_1(v)=A(pv_1,v_2,v_3)$ and 
\[ B_1(v)=pb_{11}v_1^2+b_{12}v_1v_2+b_{13}v_1v_3 + p^{-1}(b_{22}v_2^2+b_{23}v_2v_3+b_{33}v_3^2).\]
It follows from \eqref{eq:AB-MODP-CASE2} that $(A_1,B_1)$ is integral. 
We have $Q(A_1,B_1)=Q_1$ and 
$F_{(A_1,B_1)}(u)=p^{-1}F_{(A,B)}(pu_1,u_2)=(\delta_{12}'\cdot F_{(A,B)})(u)$, 
where $\delta_{12}'=\diag[p,1]$. By Proposition \ref{prop:CUBIC-SUBRING}, 
$R(A_1,B_1)=[1,p\omega,\theta]=R_1$. 

We put $\delta_{21}=\diag[1,1,p]$, $\delta_2=(\delta_{21},\delta_{12}^{-1})$ 
and  $(A_2,B_2)=\delta_2 \cdot (A,B)$. 
Then we have $A_2(v)=A(v_1,v_2,pv_3)$ and 
\[ B_2(v)=p^{-1}(b_{11}v_1^2+b_{12}v_1v_2+b_{22}v_2^2) + b_{13}v_1v_3+b_{23}v_2v_3+pb_{33}v_3^2.\]
It follows from \eqref{eq:AB-MODP-CASE2} that $(A_2,B_2)$ is integral. 
We have $Q(A_2,B_2)=Q_2$ and $F_{(A_2,B_2)}(u)=F_{(A_1,B_1)}(u)$, hence $R(A_2,B_2)=R_1$. 
The only ideals $\mathfrak{f}_p$ such that $p\mathcal{O}_k \subset \mathfrak{f}_p$ and $N(\mathfrak{f}_p)=p^2$ are 
$\tilde{\mathfrak{f}}_1$, $\tilde{\mathfrak{f}}_2$ and $\tilde{\mathfrak{f}}_3$  
which are the $p$-part of the conductors of $R_1$, $R_2$ and $R_3$, respectively. 
So we have $a_K(\tilde{\mathfrak{f}}_1)=2$ and $a_K(\tilde{\mathfrak{f}}_2)=a_K(\tilde{\mathfrak{f}}_3)=0$. 

We denote by $\Delta_{ij}$ and $\Delta'_{ij}$ the $(i,j)$-cofactors of $\omega A + aB$ and $dA-\theta B$, respectively. 
Then it follows from \eqref{eq:AB-MODP-CASE2} and $p|d$ that 
\[ -4\Delta_{33} \equiv (r^2-4s)\omega^2 \pmod{p\mathcal{O}_k}, \quad
   -4\Delta'_{22} \equiv \theta^2  \pmod{p\mathcal{O}_k}.
\]
By Lemma \ref{lem:COFACTOR-K6}, $k_6=k(\sqrt{-4\Delta_{33}}\,)=k(\sqrt{-4\Delta'_{22}}\,)$. 
Since $\omega \equiv \lambda(r^2-4s) \pmod{\tilde{\mathfrak{p}}_2}$ and $\theta \equiv 1 \pmod{\tilde{\mathfrak{p}}_1}$, 
we have  $\chi(\tilde{\mathfrak{p}}_2)=\legend{-4\Delta_{33}}{\tilde{\mathfrak{p}}_2}=\legend{r^2-4s}{p}=-1$ and 
$\chi(\tilde{\mathfrak{p}}_1)=\legend{-4\Delta'_{22}}{\tilde{\mathfrak{p}}_1}=1$ if $p\neq 2$. 
If $p=2$, then it follows from \eqref{eq:AB-MODP-CASE2}, $2|a$ and $2|d$ that 
\[ -4\Delta_{33} \equiv -3\omega^2 \pmod{8\mathcal{O}_k}, \quad
-4\Delta'_{22} \equiv \theta^2  \pmod{8\mathcal{O}_k}.
\]
Since $\omega \equiv 1 \pmod{\tilde{\mathfrak{p}}_2}$ and $\theta \equiv 1 \pmod{\tilde{\mathfrak{p}}_1}$, 
we have $\chi(\tilde{\mathfrak{p}}_2)=-1$ and 
$\chi(\tilde{\mathfrak{p}}_1)=1$. 
Hence $\tilde{\mathfrak{p}}_2$ remains prime, while $\tilde{\mathfrak{p}}_1$ splits in $k_6/k$. 
By Table \ref{tabular:T1}, $\tilde{\mathfrak{p}}_3$ remains prime in $k_6/k$. 
So we have 
\begin{align*}
1+\chi(\tilde{\mathfrak{f}}_1) &=1+\chi(\tilde{\mathfrak{p}}_2\tilde{\mathfrak{p}}_3)=1+(-1)(-1)=2=a_K(\tilde{\mathfrak{f}}_1), \\
1+\chi(\tilde{\mathfrak{f}}_2) &=1+\chi(\tilde{\mathfrak{p}}_1\tilde{\mathfrak{p}}_3)=1+1(-1)=0=a_K(\tilde{\mathfrak{f}}_2), \\
1+\chi(\tilde{\mathfrak{f}}_3) &=1+\chi(\tilde{\mathfrak{p}}_1\tilde{\mathfrak{p}}_2)=1+1(-1)=0=a_K(\tilde{\mathfrak{f}}_3).
\end{align*}

Case 3. $p$ is of type $112$ in $K$. 
We have a ring isomorphism $\varphi:\mathcal{O}_K/p\mathcal{O}_K \cong \mathbb{F}_p\oplus \mathbb{F}_p\oplus \mathbb{F}_{p^2}$. 
Let $\xi$, $r$ and $s$ be as in Case 2. 
We put $e_0=(1,1,1)$, $e_1=(0,1,0)$, $e_2=(0,0,1)$ and $e_3=(0,0,\xi)$. 
Then $\{e_i\}$ is a normalized basis of $\mathbb{F}_p\oplus \mathbb{F}_p\oplus \mathbb{F}_{p^2}$. 
By Lemma \ref{lem:PRIME-INDEX-ORDER}, 
there exists a normalized $\mathbb{Z}$-basis $\{\alpha_i\}$ of $\mathcal{O}_K$ and a positive integer $\lambda$ with $p\nmid \lambda$ such that 
$\varphi(\pi(\alpha_i)) = e_i$ ($i=0,1,2$) and $\varphi(\pi(\alpha_3))=\lambda e_3$. 
We take a pair of integral ternary quadratic forms $(A,B)$ such that 
$Q(A,B)=\mathcal{O}_K$ and the multiplication of $\mathcal{O}_K$ is given by \eqref{eq:QUARTIC-RING-STR}. 
Since $e_1^2=e_1$, $e_1e_2=e_1e_3=0$, $e_2^2=e_2$, $e_2e_3=e_3$  
and $e_3^2=-se_2-re_3$, 
we have 
$\lambda^{12}_{13} \equiv -1$, $\lambda^{13}_{22} \equiv -1$, 
$\lambda^{13}_{23} \equiv r\lambda$, $\lambda^{13}_{33} \equiv -s\lambda^2 \pmod{p}$ 
and remaining eleven $\lambda^{ij}_{k\ell}$'s are congruent to $0$ modulo $p$. 
Since $\lambda^{12}_{13} \equiv -1 \pmod{p}$, we can take an element $\gamma_2 \in \SL_2(\mathbb{Z})$ such that 
$\gamma_2 \equiv \left(\begin{array}{cc}
 b_{13} & -a_{13} \\
 b_{12} & -a_{12} 
\end{array} \right) \pmod{p}$. 
Then we have 
\begin{align*}
b_{13}A(v)-a_{13}B(v)
 &=\sum_{i\leq j} \lambda^{ij}_{13}v_iv_j \\
 &\equiv -v_1v_2+v_2^2-r\lambda v_2v_3+s\lambda^2 v_3^2 \pmod{p}, \\
b_{12}A(v)-a_{12}B(v)
 &=\sum_{i\leq j} \lambda^{ij}_{12}v_iv_j 
 \equiv v_1v_3 \pmod{p}. 
\end{align*}
Replacing $(A,B)$ by $(1_3,\gamma_2)\cdot (A,B)$, we may assume that 
\begin{equation}\label{eq:AB-MODP-CASE3}
\begin{split}
A(v) &\equiv -v_1v_2+v_2^2-r\lambda v_2v_3+s\lambda^2 v_3^2 \pmod{p}, \\
B(v) &\equiv v_1v_3 \pmod{p}. 
\end{split}
\end{equation}
Then $F_{(A,B)}(u)=au_1^3+bu_1^2u_2+cu_1u_2^2+du_2^3$ satisfies 
\[ F_{(A,B)}(u) \equiv -u_1(s\lambda^2 u_1^2 + r\lambda u_1u_2 + u_2^2) \pmod{p}. \]
Hence $a\equiv -s\lambda^2$, $b\equiv -r\lambda$, $c\equiv -1$ and $d\equiv 0 \pmod{p}$. 
Since the ring structure of $R_K=R(A,B)=[1,\omega,\theta]$ is given by \eqref{eq:CUBIC-RING-STR}, 
we have 
\begin{equation}\label{eq:CUBIC-MODP-CASE3}
\omega^2 \equiv -s\lambda^2 -r\lambda\omega+s\lambda^2\theta, \quad \theta^2 \equiv \theta, \quad
\omega\theta \equiv 0 \pmod{pR_K}.
\end{equation}
We put $\mathfrak{p}_1=[p,\omega,\theta-1]$ and $\mathfrak{p}_2=[p,p\omega,\theta]$. 
It is easy to see that $\mathfrak{p}_1$ and $\mathfrak{p}_2$ are $R_K$-ideals such that 
$\mathfrak{p}_1\mathfrak{p}_2=pR_K$. 
If we put $\tilde{\mathfrak{p}}_i=\mathfrak{p}_i\mathcal{O}_k$, then 
$p\mathcal{O}_k=\tilde{\mathfrak{p}}_1\tilde{\mathfrak{p}}_2$. 
There exists exactly one cubic ring contained in $R_K$ with index $p$ in $R_K$, namely 
$R_1=[1,p\omega,\theta]=\mathbb{Z}+\mathfrak{p}_2$. 
Hence $\mathfrak{p}_2$ is the largest $R_K$-ideal contained in $R_1$. 
By Lemma \ref{lem:CONDUCTOR2}, the $p$-part of the conductor of $R_1$ equals $\tilde{\mathfrak{p}}_2$. 
On the other hand, it is easy to see that there exist exactly two index $p$ subrings of $\mathbb{F}_p \oplus \mathbb{F}_p\oplus \mathbb{F}_{p^2}$ with unity, 
namely 
\[ B_1=\{(b_1,b_1,b_3+b_4\xi)\,|\, b_1,b_3,b_4 \in \mathbb{F}_p\}, \quad
   B_2=\mathbb{F}_p \oplus \mathbb{F}_p \oplus \mathbb{F}_p.
\]
Then $(\varphi\circ \pi)^{-1}(B_i)$'s are the quartic rings with index $p$ in $\mathcal{O}_K$. 
They are given by 
\[ Q_1=[1, p\alpha_1,\alpha_2,\alpha_3], \quad Q_2=[1,\alpha_1,\alpha_2,p\alpha_3]. \]
We now compute the cubic resolvent rings of $Q_i$'s. 
We put $\delta_{11}=\diag[p,1,1]$, $\delta_{12}=\diag[1,p]$, $\delta_1=(\delta_{11},\delta_{12}^{-1})$ 
and  $(A_1,B_1)=\delta_1 \cdot (A,B)$. 
Then we have $A_1(v)=A(pv_1,v_2,v_3)$ and 
\[ B_1(v) =pb_{11}v_1^2+b_{12}v_1v_2+b_{13}v_1v_3 + p^{-1}(b_{22}v_2^2+b_{23}v_2v_3+b_{33}v_3^2).\]
It follows from \eqref{eq:AB-MODP-CASE3} that $(A_1,B_1)$ is integral. 
We have $Q(A_1,B_1)=Q_1$ and 
$F_{(A_1,B_1)}(u)=p^{-1}F_{(A,B)}(pu_1,u_2)=(\delta_{12}'\cdot F_{(A,B)})(u)$, 
where $\delta_{12}'=\diag[p,1]$. By Proposition \ref{prop:CUBIC-SUBRING}, 
$R(A_1,B_1)=[1,p\omega,\theta]=R_1$. 

We put $\delta_{21}=\diag[1,1,p]$, $\delta_2=(\delta_{21},\delta_{12}^{-1})$ 
and  $(A_2,B_2)=\delta_2 \cdot (A,B)$. 
Then we have $A_2(v)=A(v_1,v_2,pv_3)$ and 
\[ B_2(v)=p^{-1}(b_{11}v_1^2+b_{12}v_1v_2+b_{22}v_2^2) + b_{13}v_1v_3+b_{23}v_2v_3+pb_{33}v_3^2.\]
It follows from \eqref{eq:AB-MODP-CASE3} that $(A_2,B_2)$ is integral. 
We have $Q(A_2,B_2)=Q_2$ and $F_{(A_2,B_2)}(u)=F_{(A_1,B_1)}(u)$, hence $R(A_2,B_2)=R_1$. 
The $p$-part of the conductor of $R_1$ is $\tilde{\mathfrak{p}}_2$, which is 
the only ideal $\mathfrak{f}_p$ such that $p\mathcal{O}_k \subset \mathfrak{f}_p$ and $N(\mathfrak{f}_p)=p^2$. 
So we have $a_K(\tilde{\mathfrak{p}}_2)=2$. 
On the other hand, both of $\tilde{\mathfrak{p}}_1$ and $\tilde{\mathfrak{p}}_2$ split in $k_6/k$ by Table \ref{tabular:T1}. 
Hence we have $1+\chi(\tilde{\mathfrak{p}}_2)=2=a_K(\tilde{\mathfrak{p}}_2)$. 

Case 4. $p$ is of type $111^2$ in $K$. We put $S=\mathbb{F}_p[x]/(x^2)$ and $\varepsilon=x \bmod{x^2} \in S$. 
We have a ring isomorphism $\varphi:\mathcal{O}_K/p\mathcal{O}_K \cong \mathbb{F}_p\oplus \mathbb{F}_p\oplus S$. 
We put $e_0=(1,1,1)$, $e_1=(0,1,0)$, $e_2=(0,0,1)$ and $e_3=(0,0,\varepsilon)$. 
Then $\{e_i\}$ is a normalized basis of $\mathbb{F}_p\oplus \mathbb{F}_p\oplus S$. 
By Lemma \ref{lem:PRIME-INDEX-ORDER}, 
there exists a normalized $\mathbb{Z}$-basis $\{\alpha_i\}$ of $\mathcal{O}_K$ and a positive integer $\lambda$ with $p\nmid \lambda$ such that 
$\varphi(\pi(\alpha_i)) = e_i$ ($i=0,1,2$) and $\varphi(\pi(\alpha_3))=\lambda e_3$. 
We take a pair of integral ternary quadratic forms $(A,B)$ such that 
$Q(A,B)=\mathcal{O}_K$ and the multiplication of $\mathcal{O}_K$ is given by \eqref{eq:QUARTIC-RING-STR}. 
Since $e_1^2=e_1$, $e_1e_2=e_1e_3=0$, $e_2^2=e_2$, $e_2e_3=e_3$  
and $e_3^2=0$, 
we have $\lambda^{12}_{13} \equiv \lambda^{13}_{22} \equiv -1\pmod{p}$, 
and remaining thirteen $\lambda^{ij}_{k\ell}$'s are congruent to $0$ modulo $p$. 
Since $\lambda^{12}_{13} \equiv -1 \pmod{p}$, we can take an element $\gamma_2 \in \SL_2(\mathbb{Z})$ such that 
$\gamma_2 \equiv \left(\begin{array}{cc}
 -b_{12} & a_{12} \\
  b_{13} & -a_{13} 
\end{array} \right) \pmod{p}$. 
Then we have 
\begin{align*}
-b_{12}A(v)+a_{12}B(v)
 &=\sum_{i\leq j} \lambda^{12}_{ij}v_iv_j
 \equiv -v_1v_3 \pmod{p}, \\
b_{13}A(v)-a_{13}B(v)
 &=\sum_{i\leq j} \lambda^{ij}_{13}v_iv_j \equiv -v_1v_2 +v_2^2 \pmod{p}. 
\end{align*}
Replacing $(A,B)$ by $(1_3,\gamma_2)\cdot (A,B)$, we may assume that 
\begin{equation}\label{eq:AB-MODP-CASE4}
A(v) \equiv -v_1v_3 \pmod{p}, \quad 
B(v) \equiv -v_1v_2+v_2^2 \pmod{p}. 
\end{equation}
Then $F_{(A,B)}(u)=au_1^3+bu_1^2u_2+cu_1u_2^2+du_2^3$ satisfies 
\[ F_{(A,B)}(u) \equiv u_1^2 u_2 \pmod{p}. \]
Hence $a\equiv c\equiv d \equiv 0$ and $b\equiv 1 \pmod{p}$. 
Since the ring structure of $R_K=R(A,B)=[1,\omega,\theta]$ is given by \eqref{eq:CUBIC-RING-STR}, 
we have 
\begin{equation}\label{eq:CUBIC-MODP-CASE4}
\omega^2 \equiv \omega, \quad \theta^2 \equiv 0, \quad
\omega\theta \equiv 0 \pmod{pR_K}.
\end{equation}
We put $\mathfrak{p}_1=[p,\omega,\theta]$ and $\mathfrak{p}_2=[p,\omega-1,\theta]$. 
Suppose $p^2|d$. Then $R_0=[1,\omega, \theta/p]$ becomes a cubic ring such that 
$R_K \subset R_0 \subset \mathcal{O}_k$. Since $p \nmid g=(\mathcal{O}_k:R_K)$, 
this is a contradiction. Hence $p^2 \nmid d$. 
Now it is easy to see that $\mathfrak{p}_1$ and $\mathfrak{p}_2$ are $R_K$-ideals such that 
$\mathfrak{p}_1^2=[p,\omega,p\theta]$ and $\mathfrak{p}_1^2\mathfrak{p}_2=pR_K$. 
If we put $\tilde{\mathfrak{p}}_i=\mathfrak{p}_i\mathcal{O}_k$, then 
$p\mathcal{O}_k=\tilde{\mathfrak{p}}_1^2\tilde{\mathfrak{p}}_2$. 
There exist exactly two cubic rings contained in $R_K$ with index $p$ in $R_K$, namely 
$R_1=[1,p\omega,\theta]$ and $R_2=[1,\omega,p\theta]$. 
We put $\mathfrak{f}_1=\mathfrak{p}_1\mathfrak{p}_2=[p,p\omega,\theta]$ 
and $\mathfrak{f}_2=\mathfrak{p}_1^2=[p,\omega,p\theta]$. 
Then $\mathfrak{f}_i$ is the largest $R_K$-ideal contained in $R_i$ for $i=1,2$. 
By Lemma \ref{lem:CONDUCTOR2}, the $p$-part of the conductor of $R_i$ equals $\tilde{\mathfrak{f}}_i=\mathfrak{f}_i\mathcal{O}_k$. 
We define subrings $B_i$ of $\mathbb{F}_p\oplus \mathbb{F}_p \oplus S$ by  
\begin{alignat*}{2}
B_1&=\mathbb{F}_pe_0+\mathbb{F}_pe_1+\mathbb{F}_pe_3, & \quad 
B_2&=\mathbb{F}_pe_0+\mathbb{F}_p(e_1+e_2)+\mathbb{F}_pe_3,\\ 
B_3&=\mathbb{F}_pe_0+\mathbb{F}_pe_1+\mathbb{F}_pe_2,  & \quad  
B_4&=\mathbb{F}_pe_0+\mathbb{F}_pe_2+\mathbb{F}_pe_3.
\end{alignat*}
Then it is easy to see that the four subrings $B_i$'s are the index $p$ subrings of $\mathbb{F}_p\oplus \mathbb{F}_p \oplus S$ with unity. 
Hence $(\varphi\circ \pi)^{-1}(B_i)$'s are the quartic rings with index $p$ in $\mathcal{O}_K$. 
They are given by 
\begin{alignat*}{2}
 Q_1&=[1,\alpha_1,p\alpha_2,\alpha_3],  & \quad  Q_2&=[1,p\alpha_1,\alpha_1+\alpha_2,\alpha_3],  \\
 Q_3&=[1,\alpha_1,\alpha_2,p\alpha_3],  & \quad  Q_4&=[1,p\alpha_1,\alpha_2,\alpha_3].
\end{alignat*}
We now compute the cubic resolvent rings of $Q_i$'s. 
We put $\delta_{11}=\diag[1,p,1]$, $\delta_{12}=\diag[1,p]$, $\delta_1=(\delta_{11},\delta_{12}^{-1})$ 
and  $(A_1,B_1)=\delta_1 \cdot (A,B)$. 
Then we have $A_1(v)=A(v_1,pv_2,v_3)$ and 
\[ B_1(v) =b_{12}v_1v_2++pb_{22}v_2^2+b_{23}v_2v_3+p^{-1}(b_{11}v_1^2+b_{13}v_1v_3+b_{33}v_3^2).\]
It follows from \eqref{eq:AB-MODP-CASE4} that $(A_1,B_1)$ is integral. 
We have $Q(A_1,B_1)=Q_1$ and 
$F_{(A_1,B_1)}(u)=p^{-1}F_{(A,B)}(pu_1,u_2)=(\delta_{12}'\cdot F_{(A,B)})(u)$, 
where $\delta_{12}'=\diag[p,1]$. By Proposition \ref{prop:CUBIC-SUBRING}, 
$R(A_1,B_1)=[1,p\omega,\theta]=R_1$. 

We put $\delta_{21}=\left(
\begin{array}{ccc}
p & 0 & 0 \\
1 & 1 & 0 \\
0 & 0 & 1
\end{array}
\right)$ and $\delta_2=(\delta_{21},\delta_{12}^{-1})$ 
and  $(A_2,B_2)=\delta_2 \cdot (A,B)$. 
Then we have $A_2(v)=A(pv_1+v_2,v_1+v_2,v_3)$ and 
\begin{align*}
B_2(v)
 &=pb_{11}v_1^2+(2b_{11}+b_{12})v_1v_2+b_{13}v_1v_3 \\
 &\quad +p^{-1}\left[(b_{11}+b_{12}+b_{22})v_2^2
+(b_{13}+b_{23})v_2v_3+b_{33}v_3^2 \right].
\end{align*}
It follows from \eqref{eq:AB-MODP-CASE4} that $(A_2,B_2)$ is integral. 
We have $Q(A_2,B_2)=Q_2$ and $F_{(A_2,B_2)}(u)=F_{(A_1,B_1)}(u)$, hence $R(A_2,B_2)=R_1$. 

We put $\delta_{31}=\diag[1,1,p]$, $\delta_{32}=\diag[p,1]$, $\delta_3=(\delta_{31},\delta_{32}^{-1})$ 
and  $(A_3,B_3)=\delta_3 \cdot (A,B)$. 
Then we have $B_3(v)=B(v_1,v_2,pv_3)$ and 
\[ A_3(v)=a_{13}v_1v_3+a_{23}v_2v_3+pa_{33}v_3^2+p^{-1}(a_{11}v_1^2+a_{12}v_1v_2+a_{22}v_2^2).\]
It follows from \eqref{eq:AB-MODP-CASE4} that $(A_3,B_3)$ is integral. 
We have $Q(A_3,B_3)=Q_3$ and 
$F_{(A_3,B_3)}(u)=p^{-1}F_{(A,B)}(u_1,pu_2)=(\delta_{32}'\cdot F_{(A,B)})(u)$, 
where $\delta_{32}'=\diag[1,p]$. By Proposition \ref{prop:CUBIC-SUBRING}, 
$R(A_3,B_3)=[1,\omega,p\theta]=R_2$. 

We put 
$\delta_{41}=\diag[p,1,1]$, $\delta_4=(\delta_{41},\delta_{32}^{-1})$ 
and  $(A_4,B_4)=\delta_4 \cdot (A,B)$. 
Then we have $B_4(v)=B(pv_1,v_2,v_3)$ and 
\[ A_4(v)=pa_{11}v_1^2+a_{12}v_1v_2+a_{13})v_1v_3+p^{-1}(a_{22}v_2^2+a_{23}v_2v_3+a_{33}v_3^2).\]
It follows from \eqref{eq:AB-MODP-CASE4} that $(A_4,B_4)$ is integral. 
We have $Q(A_4,B_4)=Q_4$ and $F_{(A_4,B_4)}(u)=F_{(A_3,B_3)}(u)$, 
hence $R(A_4,B_4)=R_2$. 
The only ideals $\mathfrak{f}_p$ such that $p\mathcal{O}_k \subset \mathfrak{f}_p$ and $N(\mathfrak{f}_p)=p^2$ are 
$\tilde{\mathfrak{f}}_1$ and $\tilde{\mathfrak{f}}_2$  
which are the $p$-part of the conductors of $R_1$ and $R_2$, respectively. 
So we have $a_K(\tilde{\mathfrak{f}}_1)=a_K(\tilde{\mathfrak{f}}_2)=2$. 
On the other hand, both of $\tilde{\mathfrak{p}}_1$ and $\tilde{\mathfrak{p}}_2$ split in $k_6/k$ by Table \ref{tabular:TV1}. 
Hence we have $1+\chi(\tilde{\mathfrak{f}}_i)=2=a_K(\tilde{\mathfrak{f}}_i)$ for $i=1,2$. 

Case 5. $p$ is of type $21^2$ in $K$.  Let $\xi$, $r$ and $s$ be as in Case 2 and let $S$ and $\varepsilon$ be as in Case 4. 
We have a ring isomorphism $\varphi:\mathcal{O}_K/p\mathcal{O}_K \cong \mathbb{F}_{p^2}\oplus S$. 
We put $e_0=(1,1)$, $e_1=(\xi,0)$, $e_2=(0,1)$ and $e_3=(0,\varepsilon)$. 
Then $\{e_i\}$ is a normalized basis of $\mathbb{F}_{p^2}\oplus S$. 
By Lemma \ref{lem:PRIME-INDEX-ORDER}, 
there exists a normalized $\mathbb{Z}$-basis $\{\alpha_i\}$ of $\mathcal{O}_K$ and a positive integer $\lambda$ with $p\nmid \lambda$ such that 
$\varphi(\pi(\alpha_i)) = e_i$ ($i=0,1,2$) and $\varphi(\pi(\alpha_3))=\lambda e_3$. 
We take a pair of integral ternary quadratic forms $(A,B)$ such that 
$Q(A,B)=\mathcal{O}_K$ and the multiplication of $\mathcal{O}_K$ is given by \eqref{eq:QUARTIC-RING-STR}. 
Since $e_1^2=-se_0-re_1+se_2$, $e_1e_2=e_1e_3=0$, $e_2^2=e_2$, $e_2e_3=e_3$  
and $e_3^2=0$, 
we have $\lambda^{11}_{13} \equiv s$, $\lambda^{12}_{13} \equiv r$, $\lambda^{13}_{22} \equiv -1\pmod{p}$, 
and remaining twelve $\lambda^{ij}_{k\ell}$'s are congruent to $0$ modulo $p$. 
Since $\lambda^{13}_{22} \equiv -1 \pmod{p}$, we can take an element $\gamma_2 \in \SL_2(\mathbb{Z})$ such that 
$\gamma_2 \equiv \left(\begin{array}{cc}
 b_{22} & -a_{22} \\
 b_{13} & -a_{13} 
\end{array} \right) \pmod{p}$. 
Then we have 
\begin{align*}
b_{22}A(v)-a_{22}B(v)
 &=\sum_{i\leq j} \lambda^{ij}_{22}v_iv_j
 \equiv -v_1v_3 \pmod{p}, \\
b_{13}A(v)-a_{13}B(v)
 &=\sum_{i\leq j} \lambda^{ij}_{13}v_iv_j \equiv sv_1^2 +rv_1v_2+v_2^2 \pmod{p}. 
\end{align*}
Replacing $(A,B)$ by $(1_3,\gamma_2)\cdot (A,B)$, we may assume that 
\begin{equation}\label{eq:AB-MODP-CASE5}
A(v) \equiv -v_1v_3 \pmod{p}, \quad 
B(v) \equiv sv_1^2 +rv_1v_2+v_2^2 \pmod{p}. 
\end{equation}
Then $F_{(A,B)}(u)=au_1^3+bu_1^2u_2+cu_1u_2^2+du_2^3$ satisfies 
\[ F_{(A,B)}(u) \equiv u_1^2 u_2 \pmod{p}. \]
Then there exist exactly two cubic rings contained in $R_K$ with index $p$ in $R_K$, namely 
$R_1=[1,p\omega,\theta]$ and $R_2=[1,\omega,p\theta]$ as in Case 4. 
The $p$-part of the conductor of $R_i$ equals $\tilde{\mathfrak{f}}_i=\mathfrak{f}_i\mathcal{O}_k$ 
where $\mathfrak{f}_i$ and $\mathfrak{p}_i$ ($i=1,2$) are as in Case 4. 
On the other hand, it is easy to see that there exist exactly two index $p$ subrings of $\mathbb{F}_{p^2}\oplus S$ with unity, 
namely $B_1=\mathbb{F}_p \oplus S$ and $B_2=\mathbb{F}_{p^2} \oplus \mathbb{F}_p$. 
Then $(\varphi\circ \pi)^{-1}(B_i)$'s are the quartic rings with index $p$ in $\mathcal{O}_K$. 
They are given by 
\[ Q_1=[1, p\alpha_1,\alpha_2,\alpha_3], \quad Q_2=[1,\alpha_1,\alpha_2,p\alpha_3]. \]
We put $\delta_{11}=\diag[p,1,1]$, $\delta_{12}=\diag[p,1]$, $\delta_1=(\delta_{11},\delta_{12}^{-1})$ 
and  $(A_1,B_1)=\delta_1 \cdot (A,B)$. 
Then we have $B_1(v)=B(pv_1,v_2,v_3)$ and 
\[ A_1(v) =pa_{11}v_1^2+a_{12}v_1v_2+a_{13}v_1v_3 + p^{-1}(a_{22}v_2^2+a_{23}v_2v_3+a_{33}v_3^2).\]
It follows from \eqref{eq:AB-MODP-CASE5} that $(A_1,B_1)$ is integral. 
We have $Q(A_1,B_1)=Q_1$ and 
$F_{(A_1,B_1)}(u)=p^{-1}F_{(A,B)}(u_1,pu_2)=(\delta_{12}'\cdot F_{(A,B)})(u)$, 
where $\delta_{12}'=\diag[1,p]$. By Proposition \ref{prop:CUBIC-SUBRING}, 
$R(A_1,B_1)=[1,\omega,p\theta]=R_2$. 

We put $\delta_{21}=\diag[1,1,p]$, $\delta_2=(\delta_{21},\delta_{12}^{-1})$ 
and  $(A_2,B_2)=\delta_2 \cdot (A,B)$. 
Then we have $B_2(v)=A(v_1,v_2,pv_3)$ and 
\[ A_2(v)=p^{-1}(a_{11}v_1^2+a_{12}v_1v_2+a_{22}v_2^2) + a_{13}v_1v_3+a_{23}v_2v_3+pa_{33}v_3^2.\]
It follows from \eqref{eq:AB-MODP-CASE5} that $(A_2,B_2)$ is integral. 
We have $Q(A_2,B_2)=Q_2$ and $F_{(A_2,B_2)}(u)=F_{(A_1,B_1)}(u)$, hence $R(A_2,B_2)=R_2$. 
The $p$-part of the conductor of $R_2$ is $\tilde{\mathfrak{f}}_2=\tilde{\mathfrak{p}}_1^2$. 
So we have $a_K(\tilde{\mathfrak{f}}_1)=0$ and $a_K(\tilde{\mathfrak{f}}_2)=2$. 
On the other hand, $\tilde{\mathfrak{p}}_1$ remains prime, while $\tilde{\mathfrak{p}}_2$ splits in $k_6/k$ by Table \ref{tabular:TV1}. 
Hence we have 
\begin{align*}
1+\chi(\tilde{\mathfrak{f}}_1)
 &=1+\chi(\tilde{\mathfrak{p}}_1)\chi(\tilde{\mathfrak{p}}_2)=1+(-1)1=0=a_K(\tilde{\mathfrak{f}}_1), \\
1+\chi(\tilde{\mathfrak{f}}_2)
 &=1+\chi(\tilde{\mathfrak{p}}_1^2)=1+1=2=a_K(\tilde{\mathfrak{f}}_2).
\end{align*}

Case 6. $p$ is of type $11^3$ in $K$. We put $S=\mathbb{F}_p[x]/(x^3)$ and $\varepsilon=x \bmod{x^3} \in S$. 
We have a ring isomorphism $\varphi:\mathcal{O}_K/p\mathcal{O}_K \cong \mathbb{F}_p\oplus S$. 
We put $e_0=(1,1)$, $e_1=(1,0)$, $e_2=(0,\varepsilon)$ and $e_3=(0,\varepsilon^2)$. 
Then $\{e_i\}$ is a normalized basis of $\mathbb{F}_p \oplus S$. 
By Lemma \ref{lem:PRIME-INDEX-ORDER}, 
there exists a normalized $\mathbb{Z}$-basis $\{\alpha_i\}$ of $\mathcal{O}_K$ and a positive integer $\lambda$ with $p\nmid \lambda$ such that 
$\varphi(\pi(\alpha_i)) = e_i$ ($i=0,1,2$) and $\varphi(\pi(\alpha_3))=\lambda e_3$. 
We take a pair of integral ternary quadratic forms $(A,B)$ such that 
$Q(A,B)=\mathcal{O}_K$ and the multiplication of $\mathcal{O}_K$ is given by \eqref{eq:QUARTIC-RING-STR}. 
Since $e_1^2=e_1$, $e_1e_2=e_1e_3=0$, $e_2^2=e_3$ and $e_2e_3=e_3^2=0$, 
we have $\lambda^{12}_{13} \equiv -1$, $\lambda^{12}_{22} \equiv -\lambda' \pmod{p}$, 
and remaining thirteen $\lambda^{ij}_{k\ell}$'s are congruent to $0$ modulo $p$. 
Here $\lambda'$ is a rational integer such that $\lambda\lambda'\equiv 1 \pmod{p}$. 
Since $\lambda^{12}_{13} \equiv -1 \pmod{p}$, we can take an element $\gamma_2 \in \SL_2(\mathbb{Z})$ such that 
$\gamma_2 \equiv \left(\begin{array}{cc}
  b_{12} & -a_{12} \\
 -b_{13} &  a_{13} 
\end{array} \right) \pmod{p}$. 
Then we have 
\begin{align*}
 b_{12}A(v)-a_{12}B(v)
 &=\sum_{i\leq j} \lambda^{ij}_{12}v_iv_j
 \equiv v_1v_3 + \lambda' v_2^2 \pmod{p}, \\
-b_{13}A(v)+a_{13}B(v)
 &=\sum_{i\leq j} \lambda^{13}_{ij}v_iv_j \equiv v_1v_2 \pmod{p}. 
\end{align*}
Replacing $(A,B)$ by $(1_3,\gamma_2)\cdot (A,B)$, we may assume that 
\begin{equation}\label{eq:AB-MODP-CASE6}
A(v) \equiv v_1v_3 + \lambda' v_2^2 \pmod{p}, \quad 
B(v) \equiv v_1v_2 \pmod{p}. 
\end{equation}
Then $F_{(A,B)}(u)=au_1^3+bu_1^2u_2+cu_1u_2^2+du_2^3$ satisfies 
\[ F_{(A,B)}(u) \equiv -\lambda'u_1^3 \pmod{p}. \]
Hence $b \equiv c\equiv d \equiv 0$ and $a\equiv -\lambda' \pmod{p}$. 
So we have 
\begin{equation}\label{eq:CUBIC-MODP-CASE6}
\omega^2 \equiv \lambda'\theta, \quad \theta^2 \equiv 0, \quad
\omega\theta \equiv 0 \pmod{pR_K}.
\end{equation}
We also have $p^2 \nmid d$ as in Case 4. We put $\mathfrak{p}_1=[p,\omega,\theta]$. 
Then $\mathfrak{p}_1$ is an $R_K$-ideal such that $\mathfrak{p}_1^2=[p,p\omega,\theta]$ and $\mathfrak{p}_1^3=pR_K$. 
If we put $\tilde{\mathfrak{p}}_1=\mathfrak{p}_1\mathcal{O}_k$, then 
$p\mathcal{O}_k=\tilde{\mathfrak{p}}_1^3$. 
There exists exactly one cubic ring contained in $R_K$ with index $p$ in $R_K$, namely 
$R_1=[1,p\omega,\theta]$. Then $\mathfrak{f}_1=\mathfrak{p}_1^2$ is the largest $R_K$-ideal contained in $R_1$. 
By Lemma \ref{lem:CONDUCTOR2}, the $p$-part of the conductor of $R_1$ equals $\tilde{\mathfrak{f}}_1=\tilde{\mathfrak{p}}_1^2$. 
We define subrings $B_i$ of $\mathbb{F}_p\oplus S$ by  
\[ B_1=\mathbb{F}_pe_0+\mathbb{F}_pe_2+\mathbb{F}_pe_3, \quad  
   B_2=\mathbb{F}_pe_0+\mathbb{F}_pe_1+\mathbb{F}_pe_3.
\]
Then it is easy to see that the two subrings $B_i$'s are the index $p$ subrings of $\mathbb{F}_p\oplus \mathbb{F}_p \oplus S$ with unity. 
Hence $(\varphi\circ \pi)^{-1}(B_i)$'s are the quartic rings with index $p$ in $\mathcal{O}_K$. 
They are given by $Q_1=[1,p\alpha_1,\alpha_2,\alpha_3]$ and $Q_2=[1,\alpha_1,p\alpha_2,\alpha_3]$. 
We put $\delta_{11}=\diag[p,1,1]$, $\delta_{12}=\diag[1,p]$, $\delta_1=(\delta_{11},\delta_{12}^{-1})$ 
and  $(A_1,B_1)=\delta_1 \cdot (A,B)$. 
Then we have $A_1(v)=A(pv_1,v_2,v_3)$ and 
\[ B_1(v)=pb_{11}v_1^2+b_{12}v_1v_2+b_{13}v_1v_3+p^{-1}(b_{22}v_2^2+b_{23}v_2v_3+b_{33}v_3^2).\]
It follows from \eqref{eq:AB-MODP-CASE6} that $(A_1,B_1)$ is integral. 
We have $Q(A_1,B_1)=Q_1$, $F_{(A_1,B_1)}(u)=p^{-1}F_{(A,B)}(pu_1,u_2)$ and  
$R(A_1,B_1)=[1,p\omega,\theta]=R_1$. 

We put $\delta_{21}=\diag[1,p,1]$, $\delta_2=(\delta_{21},\delta_{12}^{-1})$ 
and  $(A_2,B_2)=\delta_2 \cdot (A,B)$. 
Then we have $A_2(v)=A(v_1,pv_2,v_3)$ and 
\[ B_2(v)=pb_{22}v_2^2+b_{12}v_1v_2+b_{23}v_2v_3+p^{-1}(b_{11}v_1^2+b_{13}v_1v_3+b_{33}v_3^2).\]
It follows from \eqref{eq:AB-MODP-CASE6} that $(A_2,B_2)$ is integral. 
We have $Q(A_2,B_2)=Q_2$ and $R(A_2,B_2)=R_1$. 
So have $a_K(\tilde{\mathfrak{f}}_1)=2=1+\chi(\tilde{\mathfrak{p}}_1)^2=1+\chi(\tilde{\mathfrak{f}}_1)$. 

If $p$ is type $13$ in $K$, then there is no index $p$ subrings of $\mathcal{O}_K$ and 
there is no integral ideal of $\mathcal{O}_k$ which contains $p\mathcal{O}_k$ and has norm $p^2$ 
since $p$ is of type $3$ in $k$. 

If $p$ is type $4$ in $K$, then there is no index $p$ subrings of $\mathcal{O}_K$. 
In this case, Table \ref{tabular:T1} implies that 
$p\mathcal{O}_k=\mathfrak{p}_1\mathfrak{p}_2$ with $N(\mathfrak{p}_i)=p^i$ ($i=1,2$) and 
$\mathfrak{p}_1$ and $\mathfrak{p}_2$ remain prime in $k_6/k$. 
The ideal $\mathfrak{p}_2$ is the only ideal of $\mathcal{O}_k$ which contains $p\mathcal{O}_k$ and has norm $p^2$. 
We have $a_K(\mathfrak{p}_2)=0=1+\chi(\mathfrak{p}_2)$. 

If $p$ is of type $1^21^2$, $2^2$ or $1^4$, then 
$p$ divides $N(\Disc(k_6/k))=g$ by Table \ref{tabular:TV1}, hence $p \nmid h$. 
This completes the proof of the proposition. 
\qed

By Lemma \ref{lem:AK-FORMULA} and Proposition \ref{prop:AK-FORMULA-P}, we have
\begin{equation}\label{eq:AK-FORMULA-GENERAL}
a_K(\mathfrak{h})=\prod_{p|h} (1+\chi(\mathfrak{f}_p)).
\end{equation}
We use this formula to obtain the number of quartic rings $Q$ contained in some $S_4$-quartic fields 
with fixed cubic resolvent ring $\mathcal{O}$. 
If $\Disc(k)<0$, then the number of such quartic rings is obviously given by the sum 
\begin{equation}\label{eq:NUM-QRING-S4-DISCNEG1}
\sum_{g|f} \sum_{K \in \mathcal{K}_k(\mathfrak{g})} a_K(\mathfrak{h})=\sum_{g|f} \sum_{K \in \mathcal{K}_k(\mathfrak{g})} \prod_{p|h} (1+\chi_K(\mathfrak{f}_p)).
\end{equation}
If $\Disc(k)>0$, then the number of such quartic rings contained in some totally real $S_4$-quartic fields 
is given by the sum \eqref{eq:NUM-QRING-S4-DISCNEG1}. 
For any positive rational integer $n$, we denote by $\omega(n)$ the number of prime divisors of $n$. 
Since $\chi_K(\mathfrak{f}_p)=\pm 1$, it follows from \eqref{eq:AK-FORMULA-GENERAL}, 
that $a_K(\mathfrak{h})=2^{\omega(h)}$ or $0$ for any $K \in \mathcal{K}_k(\mathfrak{g})$. 
Since $\ker \chi_K=H_K$, $a_K(\mathfrak{h})=2^{\omega(h)}$ if and only if 
$\mathfrak{f}_p \in H_K$ for all $p|h$. 
We now simplify the sum \eqref{eq:NUM-QRING-S4-DISCNEG1}. 
Since $\mathcal{O}$ is a subring of $\mathcal{O}_k$ with index $f$ which is square free, 
there exists a unique subring $R_g$ of $\mathcal{O}_k$ such that $\mathcal{O} \subset R_g$ and 
$(\mathcal{O}_k:R_g)=g$ for each positive divisor $g$ of $f$. 
Then the conductor of $R_g$ is $\mathfrak{g}=\prod_{p|g} \mathfrak{f}_p$. 
We write $f=gh$ and $\mathfrak{h}=\prod_{p|h} \mathfrak{f}_p$. 
For any positive divisor $c$ of $g$, 
we denote by $Y_{k,c}(\mathfrak{h})$ the subgroup of $\Cl_{R_c}/\Cl_{R_c}^2$ 
generated by the ideal classes of $\mathfrak{f}_p\cap R_c$ for all $p|h$. 
For any $K \in \mathcal{K}_k(\mathfrak{g})$, 
$H_K$ can be identified with a subgroup of $\Cl_{R_g}/\Cl_{R_g}$ of index two 
as we have shown in the proof of Proposition \ref{prop:NUMBER-QUARTIC-FIELD}. 
Hence $a_K(\mathfrak{h})=2^{\omega(f)}$ if and only if $H_K \supset Y_{k,g}(\mathfrak{h})$. 
The number of such subgroups $H_K$ equals 
\[ (\Cl_{R_g}/\Cl_{R_g}^2:Y_{k,g}(\mathfrak{h}))-1=|\Cl_{R_g}^{(2)}|/|Y_{k,g}(\mathfrak{h})|-1.\]
This also equals the number of $S_4$-quartic fields $K$ such that 
$K \in \cup_{c|g} \mathcal{K}_k(\mathfrak{c})$ and $a_K(\mathfrak{h})=2^{\omega(h)}$, where we put $\mathfrak{c}=\prod_{p|c} \mathfrak{f}_p$ for $c|g$. 
Then M\"{o}bius inversion formula implies 
\[\#\{K\in \mathcal{K}_k(\mathfrak{g})\,|\, a_K(\mathfrak{h})=2^{\omega(h)}\}
=\sum_{c|g}\mu(g/c)\left(|\Cl_{R_c}^{(2)}|/|Y_{k,c}(\mathfrak{h})|-1 \right).
\]
Since $\sum_{c|g} \mu(g/c)=0$ if $g>1$, we have 
\begin{align}\label{eq:NUM-QRING-S4-DISCNEG2}
\lefteqn{\sum_{g|f} \sum_{K \in \mathcal{K}_k(\mathfrak{g})} a_K(\mathfrak{h})}\\
 &=\sum_{gh=f} 2^{\omega(h)} \sum_{c|g}\mu(g/c)|\Cl_{R_c}^{(2)}|/|Y_{k,c}(\mathfrak{h})|
 -2^{\omega(f)}. \nonumber 
\end{align}
If we rewrite the sum on the right hand side by setting $g=cd$, 
it is equal to 
\[ \sum_{c|f} |\Cl_{R_c}^{(2)}| \sum_{dh=f/c} \mu(d) 2^{\omega(h)}/|Y_{k,c}(\mathfrak{h})|.\]
We put $s=\omega(h)$ and write $h=p_1\cdots p_s$. 
For any integral ideal $\mathfrak{a}$ of $\mathcal{O}_k$ which is relatively prime to $\mathfrak{c}$, 
we denote by $[\mathfrak{a}]_c$ the ideal class of $\mathfrak{a} \cap R_c$ in $\Cl_{R_c}$. 
We define the mapping $\rho:\mathbb{F}_2^s \rightarrow Y_{k,c}(\mathfrak{h})$ 
by $\rho(a_1,\dots,a_s)=[\mathfrak{f}_{p_1}]_c^{a_1} \cdots [\mathfrak{f}_{p_s}]_c^{a_s}\Cl_{R_c}^2$. 
Then $\rho$ is a surjective group homomorphism, hence
$|\ker \rho|=2^{\omega(h)}/|Y_{k,c}(\mathfrak{h})|$. 
Moreover $\ker \rho$ is identified with the set 
$T_{c,h}=\{t\in \mathbb{N}\,:\, t|h,\; [\prod\nolimits_{p|t} \mathfrak{f}_p]_c \in \Cl_{R_c}^2\}$. 
For any positive divisor $t$ of $h$, we put 
$\delta_c(t)=1$ if $[\prod\nolimits_{p|t} \mathfrak{f}_p]_c \in \Cl_{R_c}^2$, 
otherwise $\delta_c(t)=0$. Then we have 
\[ \sum_{t|h} \delta_c(t)=|T_{c,h}|=|\ker \rho|=2^{\omega(h)}/|Y_{k,c}(\mathfrak{h})|.\]
By M\"{o}bius inversion formula, we have 
\[ \sum_{dh=f/c} \mu(d) 2^{\omega(h)}/|Y_{k,c}(\mathfrak{h})|
=\delta_c(f/c).
\]
Hence we can rewrite \eqref{eq:NUM-QRING-S4-DISCNEG2} as follows. 
\begin{equation}\label{eq:NUM-QRING-S4-DISCNEG3}
\sum_{g|f} \sum_{K \in \mathcal{K}_k(\mathfrak{g})} a_K(\mathfrak{h})
=\sum_{c|f} |\Cl_{R_c}^{(2)}| \,\delta_c(f/c)-2^{\omega(f)}.
\end{equation}
If we rewrite the sum on the right hand side of \eqref{eq:NUM-QRING-S4-DISCNEG3} by setting $c=g$ and $f=gh$, 
it is equal to $\sum_{g|f} |\Cl_{R_g}^{(2)}| \,\delta_g(h)$. 
By the definition of $\delta_g(h)$, we have 
$\delta_g(h)=1$ if $[\mathfrak{h}]_g \in \Cl_{R_g}^2$, otherwise 
$\delta_g(h)=0$. In \S 4, we defined the ideal $\mathfrak{j}(\mathcal{O},R_g)$ 
which is characterized by the largest $R_g$-ideal contained in $\mathcal{O}$. 
Since $\mathfrak{j}(\mathcal{O},R_g)\mathcal{O}_k=\mathfrak{h}$, 
we have $\mathfrak{h}\cap R_g=\mathfrak{j}(\mathcal{O},R_g)$ by Lemma \ref{lem:EXTEND-RESTRICT-IDEALS}. 
Hence $[\mathfrak{h}]_g$ is the ideal class of $\mathfrak{j}(\mathcal{O},R_g)$ in $\Cl_{R_g}$. 
Therefore the subgroup $X(\mathcal{O},R_g)$ defined in \S 4 is trivial if and only if $\delta_g(h)=1$. 
Thus the equation $\delta_g(h)=2-|X(\mathcal{O},R_g)|$ holds for all $g|f$. 
We finally obtain the following formula:
\begin{equation}\label{eq:NUM-QRING-S4-DISCNEG4}
\sum_{g|f} \sum_{K \in \mathcal{K}_k(\mathfrak{g})} a_K(\mathfrak{h})
=\sum_{g|f} |\Cl_{R_g}^{(2)}| (2-|X(\mathcal{O},R_g)|)-2^{\omega(f)}.
\end{equation}

If $\Disc(k)>0$, then the number of quartic rings contained in some quartic fields with fixed cubic resolvent ring $\mathcal{O}$ 
is given by 
\begin{equation}\label{eq:NUM-QRING-S4-DISCPOS1}
\sum_{g|f} \sum_{K \in \mathcal{K}_k(\mathfrak{g}\mathfrak{f}_\infty)} a_K(\mathfrak{h})=\sum_{g|f} \sum_{K \in \mathcal{K}_k(\mathfrak{g}\mathfrak{f}_\infty)} \prod_{p|h} (1+\chi_K(\mathfrak{f}_p)).
\end{equation}
By the same argument as above, we obtain the following formula:
\begin{equation}\label{eq:NUM-QRING-S4-DISCPOS2}
\sum_{g|f} \sum_{K \in \mathcal{K}_k(\mathfrak{g}\mathfrak{f}_\infty)} a_K(\mathfrak{h})
=\sum_{g|f} |\Cl_{R_g,+}^{(2)}| (2-|X_+(\mathcal{O},R_g)|)-2^{\omega(f)}.
\end{equation}

It remains to count the number of quartic rings contained in the quartic algebra 
$K=\mathbb{Q}\oplus k$ with fixed cubic resolvent ring $\mathcal{O}$. 
Let $\{1,\omega,\theta\}$ be a normalized basis of $\mathcal{O}_k$ such that the multiplication of $\mathcal{O}_k$ is 
given by \eqref{eq:CUBIC-RING-STR}. We put $F_k(u)=au_1^3+bu_1^2u_2+cu_1u_2^2+du_2^3$. 
We put $\alpha_0=(1,1)$, $\alpha_1=(1,0)$, $\alpha_2=(0,-\omega)$ and $\alpha_3=(0,-\theta)$. 
Then $\{\alpha_i\}$ is a normalized basis of $\mathcal{O}_K=\mathbb{Z}\oplus \mathcal{O}_k$ 
and the multiplication of $\mathcal{O}_K$ is given by 
\begin{gather*}
\alpha_1^2=\alpha_1, \quad \alpha_1\alpha_2=\alpha_1\alpha_3=0, \quad
\alpha_2^2 =-ac+ac\alpha_1-b\alpha_2+a\alpha_3, \\
\alpha_2\alpha_3
=-ad+ad\alpha_1,\quad \alpha_3^2=-bd+bd\alpha_1-d\alpha_2+c\alpha_3. 
\end{gather*}
We put 
\[ A(v)=-v_1v_3-av_2^2-bv_2v_3, \quad B(v)=-v_1v_2+cv_2v_3+dv_3^2.\]
Then we have $Q(A,B)=\mathcal{O}_K$, $F_{(A,B)}(u)=F_k(u)$ and $R(A,B)=\mathcal{O}_k$. 
For each prime divisor $p$ of $f$, we denote by $a_K(\mathfrak{f}_p)$ 
the number of quartic rings $Q$ with index $p$ in $\mathcal{O}_K$ such that 
the unique cubic resolvent ring of $Q$ has conductor $\mathfrak{f}_p$. 
By the same argument as in the proof of Proposition \ref{prop:AK-FORMULA-P}, 
we have $a_K(\mathfrak{f}_p)=2$ for all $p|f$, hence $a_K(\mathfrak{f})=2^{\omega(f)}$. 
Hence the number of quartic rings contained in $\mathbb{Q} \oplus k$ 
or in some $S_4$-quartic fields, which are totally real if $\Disc(k)>0$, 
is given by the sum 
\begin{equation}\label{eq:NUM-QRING-S4-DISCNEG5}
\sum_{g|f} |\Cl_{R_g}^{(2)}| (2-|X(\mathcal{O},R_g)|). 
\end{equation}
If $\Disc(k)>0$, the number of quartic rings contained in $\mathbb{Q} \oplus k$ 
or in some $S_4$-quartic fields is 
given by the sum 
\begin{equation}\label{eq:NUM-QRING-S4-DISCPOS3}
\sum_{g|f} |\Cl_{R_g,+}^{(2)}| (2-|X_+(\mathcal{O},R_g)|).
\end{equation}
We defined the subset $L(\mathcal{O})$ of $L$ in \S 1. 
We recall $L_i(\mathcal{O})=L(\mathcal{O}) \cap V_i$ for $i=1,2,3$. 
Since $k$ is a non-Galois cubic field, 
$\Aut(Q(A,B))$ is trivial for any $(A,B) \in L(\mathcal{O})$ whether 
$Q(A,B)$ is contained in $\mathbb{Q}\oplus k$ or contained in an $S_4$-quartic field. 
We write $x=(A,B) \in L(\mathcal{O})$ and $\mu(x)=1/|\Gamma_{(A,B)}|$. 
Then $\mu(x)=1/|\Aut(Q(A,B))|=1$ by Proposition \ref{prop:ISOTROPY-AUT}. 
Now it follows from Theorem \ref{thm:BHARGAVA3}, Corollary \ref{cor:UNIQUE-RESOLVENT}, 
\eqref{eq:NUM-QRING-S4-DISCNEG5} and \eqref{eq:NUM-QRING-S4-DISCPOS3} that 
\begin{equation}\label{eq:L-S4CLNPOS}
\begin{split}
\sum_{x \in \Gamma\backslash L_1(\mathcal{O})} \mu(x) 
&=\sum_{g|f} |\Cl_{R_g}^{(2)}| (2-|X(\mathcal{O},R_g)|), \\
\sum_{x \in \Gamma\backslash L_1(\mathcal{O})} \mu(x) 
+\sum_{x \in \Gamma\backslash L_3(\mathcal{O})} \mu(x) 
&=\sum_{g|f} |\Cl_{R_g,+}^{(2)}| (2-|X_+(\mathcal{O},R_g)|)
\end{split}
\end{equation}
if $\Disc(k)>0$, and 
\begin{equation}\label{eq:L-S4CLNNEG}
\sum_{x \in \Gamma\backslash L_2(\mathcal{O})} \mu(x) 
=\sum_{g|f} |\Cl_{R_g}^{(2)}| (2-|X(\mathcal{O},R_g)|)
\end{equation}
if $\Disc(k)<0$. 
So we finally complete the proof of Theorem \ref{thm:MAIN-S4A4} for a non-Galois cubic field $k$ by Corollary \ref{cor:LHAT-CLN2}, \eqref{eq:L-S4CLNPOS} and \eqref{eq:L-S4CLNNEG}. 

\section{The case $G=A_4$}
Let $k$ be a Galois cubic field and $\mathcal{O}$ be an order of $k$ 
such that the index $f=(\mathcal{O}_k:\mathcal{O})$ is square free. 
In this section, we study the number of quartic rings whose cubic resolvent rings are isomorphic to $\mathcal{O}$. 
We use the same notations as in the previous section.
The argument in $G=S_4$ case also works well in $G=A_4$ case. However we must be careful because $\Gal(k/\mathbb{Q})$ is nontrivial. 
Let $\sigma$ be a generator of $\Gal(k/\mathbb{Q})$.  
We denote by $\mathfrak{f}$ the conductor of $\mathcal{O}$ and denote by $\mathfrak{f}_p$ the $p$-part of $\mathfrak{f}$ for 
any prime divisor $p$ of $f$. Since $p\mathcal{O}_k \subset \mathfrak{f}_p$ and $N(\mathfrak{f}_p)=p^2$, 
$\mathfrak{f}_p$ is a product of two distinct prime ideals if $p$ is of type $111$ in $k$,  
or is a square of a prime ideal if $p$ is of type $1^3$ in $k$. 
We say that $\mathfrak{f}$ is the \textit{first kind} if $f$ has a prime divisor $p$ of type $111$ in $k$, 
otherwise we say that $\mathfrak{f}$ is the \textit{second kind}.
If $\mathfrak{f}$ is the first kind, then $\sigma(\mathfrak{f})\neq \mathfrak{f}$ and $\sigma(\mathcal{O})\neq \mathcal{O}$, 
hence $\Aut(\mathcal{O})$ is trivial. 
If $\mathfrak{f}$ is the second kind, then $\sigma(\mathfrak{f})=\mathfrak{f}$ and $\sigma(\mathcal{O})=\mathcal{O}$, 
hence $|\Aut(\mathcal{O})|=3$.

For each positive divisor $g$ of $f$, we put $\mathfrak{g}=\prod_{p|g} \mathfrak{f}_p$ and $R_g=\mathbb{Z}+\mathfrak{g}$. 
We denote by $\mathcal{K}_k(\mathfrak{g})$ the set of isomorphism classes of totally real $A_4$-quartic fields 
$K$ such that the unique cubic resolvent ring of the maximal order $\mathcal{O}_K$ 
is isomorphic to $R_g=\mathbb{Z}+\mathfrak{g}$.

Let $K$ be an $A_4$-quartic field whose cubic resolvent field is $k$ 
and assume that $\mathcal{O}_K$ contains a quartic ring $Q$ whose cubic resolvent ring is isomorphic to $\mathcal{O}$. 
We denote by $k_6$ one of the three conjugate sextic fields contained in $\tilde{K}$. 
The maximal order $\mathcal{O}_K$ has a unique cubic resolvent ring up to isomorphism. 
Since $Q$ has a cubic resolvent ring which is isomorphic to $\mathcal{O}$, 
we can take a cubic resolvent ring $R_K$ of $\mathcal{O}_K$ such that $\mathcal{O} \subset R_K$. 
We fix $R_K$. Then we have $f=gh$, $g=(\mathcal{O}_k:R_K)$, $h=(\mathcal{O}_K:Q)$, $N(\Disc(k_6/k))=g^2$ 
and the conductor of $R_K$ is $\mathfrak{g}$. Hence $R_K=\mathbb{Z}+\mathfrak{g}=R_g$. 
By the same argument as in the proof of Proposition \ref{prop:MAXIMAL-ORDER-CONDUCTOR}, 
we see that $\mathfrak{g}$ is a conjugate of $\Disc(k_6/k)$. 

\begin{lem}\label{lem:GODD}
Let $k$, $k_6$, $\mathfrak{g}$ and $g$ be as above. Then $g$ is odd, hence $\mathfrak{g}$ is a square free ideal. 
\end{lem}
\Proof
Suppose $g$ is even. If $2$ is of type $3$ in $k$, then $2^3|g$. 
By computing double coset decompositions of $A_4$, 
we see that two prime ideals dividing $2$ ramifies in $k_6/k$ if $2$ is of type $111$, hence $2^2|g$. 
Similarly we see that $g$ is odd if $2$ is of type $1^3$.
\qed

We first assume $g>1$. Since $\mathfrak{g}=\Disc(k_6/k)$ is a square free ideal of $\mathcal{O}_k$ and has norm $g^2$, 
the $p$-part $\mathfrak{f}_p$ of $\mathfrak{g}$ is a product of two distinct prime ideals for each 
prime divisor $p$ of $g$. Hence $\sigma(\mathfrak{g})\neq \mathfrak{g}$, $\sigma(R_K) \neq R_K$ 
and $\Aut(R_K)$ is trivial. 
Since $Q$ has a unique cubic resolvent ring $R$ which is contained in $R_K$ and is isomorphic to $\mathcal{O}$, 
we have $R=\sigma^i(\mathcal{O})$ for some $i$. 
Hence the conductor of $R$ is $\sigma^i(\mathfrak{f})$. 
Since $R\subset R_K$, The conductor $\mathfrak{g}$ of $R_K$ divides $\sigma^i(\mathfrak{f})$. 
So we have $\mathfrak{g}=\sigma^i(\mathfrak{g})$. This implies $i=0$ and $R=\mathbb{Z}+\mathfrak{f}=\mathcal{O}$. 
By choosing $k_6$ suitably, we may assume that $\Disc(k_6/k)=\mathfrak{g}$. 
We denote by $a_K(\mathfrak{h})$ the number of quartic rings $Q$ with index $h$ in $\mathcal{O}_K$ 
whose cubic resolvent rings are isomorphic to $\mathcal{O}$.
We denote by $H=H_K$ the subgroup of $I_k(\mathfrak{g})$ corresponding to 
the quadratic extension $k_6/k$ by class field theory 
and $\chi=\chi_K$ the character of $I_k(\mathfrak{g})$ such that $\ker \chi=H$. 
Then Lemma \ref{lem:AK-FORMULA} and Proposition \ref{prop:AK-FORMULA-P} are still valid in this case, 
hence we have $a_K(\mathfrak{h})=\prod_{p|h} (1+\chi(\mathfrak{f}_p))$. 
By the same argument as in the previous section, we have 
\begin{equation}\label{eq:NUM-QRING-A4-GGT1}
\sum_{K \in \mathcal{K}_k(\mathfrak{g})} a_K(\mathfrak{h})
=2^{\omega(h)} \sum_{c|g}\mu(g/c)|\Cl_{R_c}^{(2)}|/|Y_{k,c}(\mathfrak{h})|
\end{equation}
if $g>1$, where $\mathcal{K}_k(\mathfrak{g})$, $\omega(h)$ and $Y_{k,c}(\mathfrak{h})$ are the same as before. 

We next assume that $\mathfrak{f}$ is the first kind and $g=1$. 
Then we have $R_K=\mathcal{O}_k$. 
Let $p_1$ be a prime divisor of $g$ such that $p_1$ is of type $111$ in $k$ and 
$p_1\mathcal{O}_k = \mathfrak{p}_1 \sigma(\mathfrak{p}_1) \sigma^2(\mathfrak{p}_1)$ be the prime ideal decomposition. 
We may assume $\mathfrak{f}_{p_1}=\mathfrak{p}_1 \sigma(\mathfrak{p}_1)$. 
Since $R_{p_1}=\mathbb{Z}+\mathfrak{f}_{p_1}$, there exist exactly three index $p$ subrings of $\mathcal{O}_k$, 
namely $\sigma^i(R_{p_1})$ for $i=0,1,2$. 
We denote by $H=H_K$ the subgroup of $I_k$ corresponding to 
the unramified quadratic extension $k_6/k$ by class field theory 
and $\chi=\chi_K$ the character of $I_k$ such that $\ker \chi=H$. 
Hence $H$ contains $P_k$ and $(I_k:H)=2$. 
By the same argument as in the proof of Proposition \ref{prop:AK-FORMULA-P}, 
we see that there exist exactly $1+\chi(\sigma^i(\mathfrak{f}_{p_1}))$ index $p$ subrings 
of $\mathcal{O}_K$ whose cubic resolvent rings are $\sigma^i(R_{p_1})$. 
Since $\prod_{i=0}^2 \sigma^i(\mathfrak{f}_{p_1})=p^2\mathcal{O}_k$, 
the number of $i \in \{0,1,2\}$ such that $\chi(\sigma^i(\mathfrak{f}_{p_1}))=1$ 
is equal to $3$ or $1$. 
We take an $i$ such that $\chi(\sigma^i(\mathfrak{f}_{p_1}))=1$. 
We put $b_{k_6}(\sigma^i(\mathfrak{f}))=\prod_{p|f} (1+\chi(\sigma^i(\mathfrak{f}_p)))$. 
By Lemma \ref{lem:AK-FORMULA} and Proposition \ref{prop:AK-FORMULA-P}, 
the number of quartic rings $Q$ with index $f$ in $\mathcal{O}_K$ whose cubic resolvent rings are $\sigma^i(\mathcal{O})$ 
is given by 
\[ 2\prod_{p|f,\, p\neq p_1} (1+\chi(\sigma^i(\mathfrak{f}_p)))
=\prod_{p|f} (1+\chi(\sigma^i(\mathfrak{f}_p)))=b_{k_6}(\sigma^i(\mathfrak{f})).
\]
Hence $b_{k_6}(\sigma^i(\mathfrak{f}))=2^{\omega(f)}$ or $0$. 
If we identify $H$ with a subgroup of $\Cl_k/\Cl_k^2$ with index $2$, then 
$b_{k_6}(\sigma^i(\mathfrak{f}))=2^{\omega(f)}$ if and only if 
$H\supset Y_{k,1}(\sigma^i(\mathfrak{f}))$, 
where $Y_{k,1}(\mathfrak{f})$ is the subgroup of $\Cl_k/\Cl_k^2$ generated by 
the ideal classes of $\mathfrak{f}_p$ for all $p|f$. 
We denote by $\mathcal{H}$ the set of subgroups $H$ of $\Cl_k/\Cl_k^2$ 
satisfying $(\Cl_k/\Cl_k^2:H)=2$ and $H \supset Y_{k,1}(\mathfrak{f})$. 
Then we have $|\mathcal{H}|=|(\Cl_k/\Cl_k^2)/Y_{k,1}(\mathfrak{f})|-1$. 
For any $H \in \mathcal{H}$, $\sigma^i(H)$ ($i=0,1,2$) are distinct subgroups of $\Cl_k/\Cl_k^2$. 
To prove this, suppose $\sigma(H)=H$. 
We take an element $c \in \Cl_k/\Cl_k^2$ such that $c \notin H$. 
Then $\sigma(c), \sigma^2(c) \notin H$. 
Hence each $\sigma^i(c)H$ is a generator of $(\Cl_k/\Cl_k^2)/H$. 
So $\prod_{i=0}^2 \sigma^i(c)H$ is not identity. 
This contradicts the fact that the norm of an ideal of $\mathcal{O}_k$ is a principal ideal. 
We now decompose $\mathcal{H}$ into two disjoint subsets $\mathcal{H}_0$ and $\mathcal{H}_1$ 
by setting
\[ \mathcal{H}_0=\mathcal{H} \cap \sigma(\mathcal{H}) \cap \sigma^2(\mathcal{H}),\quad 
   \mathcal{H}_1=\mathcal{H}\smallsetminus \mathcal{H}_0.
\]
We assume $H \in \mathcal{H}_0$. Then $H\supset \sigma^i(Y_{k,1}(\mathfrak{f}))$ for $i=0,1,2$. 
It follows from the definition of $Y_{k,1}(\mathfrak{f})$ that $\sigma^i(Y_{k,1}(\mathfrak{f}))=Y_{k,1}(\sigma^i(\mathfrak{f}))$ . 
So we have $b_{k_6}(\sigma^i(\mathfrak{f}))=2^{\omega(f)}$ for each $i$. 
Hence we have $a_K(\mathfrak{f})=3 \cdot 2^{\omega(f)}$. 
We assume $H \in \mathcal{H}_1$. Then it is easy to see that $H\supset Y_{k,1}(\mathfrak{f})$, 
but $H\not\supset \sigma^i(Y_{k,1}(\mathfrak{f}))$ for $i=1,2$. 
So we have $b_{k_6}(\mathfrak{f})=2^{\omega(f)}$ and $b_{k_6}(\sigma^i(\mathfrak{f}))=0$ for $i=1, 2$. 
Hence we have $a_K(\mathfrak{f})=2^{\omega(f)}$. 
If $H \in \mathcal{H}_0$, then the three conjugates subgroups $\sigma^i(H)$ in $\mathcal{H}_0$ correspond to 
one isomorphism class of $A_4$-quartic field $K$.  
If $H \in \mathcal{H}_1$, then $\sigma^i(H) \notin \mathcal{H}$ for $i=1,2$ and  
only one  $H$ in $\mathcal{H}_1$ correspond to one isomorphism class of $A_4$-quartic field $K$. 
Therefore we obtain the following formula:
\begin{equation}\label{eq:NUM-QRING-A4-GEQ1K1}
\sum_{K \in \mathcal{K}_k(1)} a_K(\mathfrak{f})
=2^{\omega(f)} \left(|\Cl_k^{(2)}|/|Y_{k,1}(\mathfrak{f})|-1 \right).
\end{equation}
By \eqref{eq:NUM-QRING-A4-GGT1} and \eqref{eq:NUM-QRING-A4-GEQ1K1}, we have
\[\sum_{g|f} \sum_{K \in \mathcal{K}_k(\mathfrak{g})} a_K(\mathfrak{h})
=\sum_{g|f} 2^{\omega(h)} \sum_{c|g}\mu(g/c)|\Cl_{R_c}^{(2)}|/|Y_{k,c}(\mathfrak{h})|-2^{\omega(f)}.
\]
This is the same formula as \eqref{eq:NUM-QRING-S4-DISCNEG2}. 
Hence we can rewrite this as 
\begin{equation}\label{eq:NUM-QRING-A4-K1}
\sum_{g|f} \sum_{K \in \mathcal{K}_k(\mathfrak{g})} a_K(\mathfrak{h})
=\sum_{g|f} |\Cl_{R_g}^{(2)}| (2-|X(\mathcal{O},R_g)|)-2^{\omega(f)}
\end{equation}
provided that $\mathfrak{f}$ is the first kind. 

We finally assume that $\mathfrak{f}$ is the second kind and $g=1$. 
For each prime divisor $p$ of $f$, $p\mathcal{O}_k$ is cube of a prime ideal $\mathfrak{p}$ 
and the $p$-part of $\mathfrak{f}$ is $\mathfrak{f}_p=\mathfrak{p}^2$. 
There exist exactly two index $p$ subrings of $\mathcal{O}_K$ whose cubic resolvent rings are 
$R_p=\mathbb{Z}+\mathfrak{p}^2$ as we saw in Case 6 of the proof of Proposition \ref{prop:AK-FORMULA-P}. 
Hence $a_K(\mathfrak{f})=2^{\omega(f)}$. 
The number of unramified totally real quadratic extensions of $k$ is equal to 
$|\Cl_k/\Cl_k^2|-1$ by class field theory. 
Three conjugates of such quadratic extensions $k_6$ correspond to 
one isomorphism class of $A_4$-quartic field $K \in \mathcal{K}_k(1)$. 
Therefore we obtain the following formula:
\begin{equation}\label{eq:NUM-QRING-A4-GEQ1K2}
\sum_{K \in \mathcal{K}_k(1)} a_K(\mathfrak{f})
=\frac{2^{\omega(f)}}{3} \left(|\Cl_k^{(2)}|-1 \right).
\end{equation}
For any positive divisor $g$ of $f$ with $g>1$, 
$\mathfrak{g}=\prod_{p|g} \mathfrak{f}_p$ is square of an ideal, 
hence $\mathcal{K}_k(\mathfrak{g})=\emptyset$ by Lemma \ref{lem:GODD}. 
Thus we have
\begin{equation}\label{eq:NUM-QRING-A4-K2}
\sum_{g|f} \sum_{K \in \mathcal{K}_k(\mathfrak{g})} a_K(\mathfrak{h})
=\frac{2^{\omega(f)}}{3} \left(|\Cl_k^{(2)}|-1 \right)
\end{equation}
provided that $\mathfrak{f}$ is the second kind. 

We denote by $\mathcal{K}_k(\mathfrak{g}\mathfrak{f}_\infty)$ the set of isomorphism classes of 
totally real or totally imaginary $A_4$-quartic fields $K$ such that the unique cubic resolvent ring of 
the maximal order $\mathcal{O}_K$ is isomorphic to $R_g=\mathbb{Z}+\mathfrak{g}$. 
By the same argument as above, we obtain the following formula:
If $\mathfrak{f}$ is the first kind, then we have
\begin{equation}\label{eq:NUM-QRING-A4-L1L3K1}
\sum_{g|f} \sum_{K \in \mathcal{K}_k(\mathfrak{g}\mathfrak{f}_\infty)} a_K(\mathfrak{h})
=\sum_{g|f} |\Cl_{R_g,+}^{(2)}| (2-|X_+(\mathcal{O},R_g)|)-2^{\omega(f)}.
\end{equation}
If $\mathfrak{f}$ is the second kind, then we have
\begin{equation}\label{eq:NUM-QRING-A4-L1L3K2}
\sum_{g|f} \sum_{K \in \mathcal{K}_k(\mathfrak{g}\mathfrak{f}_\infty)} a_K(\mathfrak{h})
=\frac{2^{\omega(f)}}{3} \left(|\Cl_{k,+}^{(2)}|-1 \right).
\end{equation}

It remains to count the number of quartic rings $Q$ contained in the quartic algebra 
$K=\mathbb{Q}\oplus k$ whose cubic resolvent rings are isomorphic to $\mathcal{O}$. 
By the same calculation as before, we see that $a_K(\mathfrak{f})=3 \cdot 2^{\omega(f)}$ 
if $\mathfrak{f}$ is the first kind and $a_K(\mathfrak{f})=2^{\omega(f)}$ 
if $\mathfrak{f}$ is the second kind. 
Since $\Aut(K)$ is a cyclic group of order three, we have $|\Aut(Q)|=1$ or $3$. 
We denote by $[Q]$ the isomorphism class of $Q$. 
If we denote by $\psi$ the correspondence $Q\mapsto [Q]$, then 
the equation $|\psi^{-1}([Q])|\cdot |\Aut(Q)|=3$ always holds. 
Hence we have 
\[ \sum_{[Q]} \frac{1}{|\Aut(Q)|}
=\sum_{[Q]} \frac{1}{3}|\psi^{-1}([Q])|
=\frac{1}{3} \sum_{Q} 1
=\frac{1}{3}a_K(\mathfrak{f}).
\]
So we have
\begin{equation}\label{eq:A4-SPLIT-ALGEBRA}
\sum_{[Q]} \frac{1}{|\Aut(Q)|}
=\left\{
\begin{array}{ll}
2^{\omega(f)} , &\quad \text{if $\mathfrak{f}$ is the first kind},  \\
2^{\omega(f)}/3 , &\quad \text{if $\mathfrak{f}$ is the second kind}. 
\end{array}
\right.
\end{equation}
It is clear that $\Aut(Q)$ is trivial for any quartic ring $Q$ contained in some $A_4$-quartic fields. 
Now it follows from Theorem \ref{thm:BHARGAVA3}, Corollary \ref{cor:UNIQUE-RESOLVENT}, \eqref{eq:NUM-QRING-A4-K1}, 
\eqref{eq:NUM-QRING-A4-K2} and \eqref{eq:A4-SPLIT-ALGEBRA} that
\begin{align}\label{eq:L1-A4CLN}
\lefteqn{\sum_{x \in \Gamma\backslash L_1(\mathcal{O})} \mu(x)}\\
 &=\left\{
\begin{array}{ll}
\sum_{g|f} |\Cl_{R_g}^{(2)}| (2-|X(\mathcal{O},R_g)|) , &\quad \text{if $\mathfrak{f}$ is the first kind},  \\
 (2^{\omega(f)}/3)\, |\Cl_k^{(2)}|, &\quad \text{if $\mathfrak{f}$ is the second kind}.
\end{array}
\right.\nonumber
\end{align}
Similarly it follows from Theorem \ref{thm:BHARGAVA3}, Corollary \ref{cor:UNIQUE-RESOLVENT}, \eqref{eq:NUM-QRING-A4-L1L3K1}, 
\eqref{eq:NUM-QRING-A4-L1L3K2} and \eqref{eq:A4-SPLIT-ALGEBRA} that 
\begin{align}\label{eq:L1L3-A4CLN}
\lefteqn{\sum_{x \in \Gamma\backslash L_1(\mathcal{O})} \mu(x)+\sum_{x \in \Gamma\backslash L_3(\mathcal{O})} \mu(x)}\\
 &=\left\{
\begin{array}{ll}
\sum_{g|f} |\Cl_{R_g,+}^{(2)}| (2-|X_+(\mathcal{O},R_g)|) , &\quad \text{if $\mathfrak{f}$ is the first kind},  \\
 (2^{\omega(f)}/3)\, |\Cl_{k,+}^{(2)}|, &\quad \text{if $\mathfrak{f}$ is the second kind}.
\end{array}
\right.\nonumber
\end{align}
On the other hand, $|\Aut(\mathcal{O})|$ equals $1$ or $3$ according as $\mathfrak{f}$ is the first kind or the second kind. 
It follows from Corollary \ref{cor:LHAT-CLN2} and Lemma \ref{lem:LHAT-CLN3} that 
\begin{align}
\lefteqn{\sum_{y \in \Gamma\backslash \hat{L}_1(\mathcal{O})} \mu(y)} \label{eq:LHATL1-CLNA4}\\
 &=\left\{
\begin{array}{ll}
\sum_{g|f} |\Cl_{R_g,+}^{(2)}| \left(2-|X_+(\mathcal{O},R_g) \right) , &\quad \text{if $\mathfrak{f}$ is the first kind}, \\
(2^{\omega(f)}/3)\,|\Cl_{k,+}^{(2)}| , &\quad \text{if $\mathfrak{f}$ is the second kind},
\end{array}
\right. \nonumber \\
\lefteqn{\sum_{y \in \Gamma\backslash \hat{L}_1(\mathcal{O})} \mu(y)+\sum_{y \in \Gamma\backslash \hat{L}_3(\mathcal{O})} \mu(y)} \label{eq:LHATL1L3-CLNA4}\\
 &=\left\{
\begin{array}{ll}
4\sum_{g|f} |\Cl_{R_g}^{(2)}| \left(2-|X(\mathcal{O},R_g) \right) , &\quad \text{if $\mathfrak{f}$ is the first kind}, \\
4(2^{\omega(f)}/3)\,|\Cl_{k}^{(2)}| , &\quad \text{if $\mathfrak{f}$ is the second kind}.
\end{array}
\right.
\nonumber
\end{align}
Here we used the elementary fact $\sum_{g|f} 1=2^{\omega(f)}$. 
So we finally complete the proof of Theorem \ref{thm:MAIN-S4A4} for a Galois cubic field $k$ 
by the equations \eqref{eq:L1L3-A4CLN}, \eqref{eq:LHATL1-CLNA4}, \eqref{eq:L1-A4CLN} and \eqref{eq:LHATL1L3-CLNA4}. 

\section{Proof of Theorem \ref{thm:MAIN}}
We say that a rational intger $n$ is a \textit{fundamental discriminant} 
if $n$ equals the discriminant of a quadratic field. 
Let $k_1$ be a quadratic field. We fix $k_1$ and put $k=\mathbb{Q}\oplus k_1$. 
In this section, we consider quartic rings whose cubic resolvent rings are contained in 
the cubic algebra $k$. If $Q$ is such a quartic ring, then the quartic algebra $K=Q\otimes_\mathbb{Z} \mathbb{Q}$ 
is one of the followings: $K$ is a quartic field with $G=D_4$ or $C_4$; $K$ is a direct sum of two distinct quadratic fields; 
$K=\mathbb{Q}\oplus\mathbb{Q}\oplus k_1$.

We restrict ourselves to the case such that $\Disc(Q)=\Disc(k_1)$. 
This implies that $Q=\mathcal{O}_K$, $\Disc(K)=\Disc(k_1)$ 
and the cubic resolvent ring of $Q$ is isomorphic to $\mathcal{O}_k=\mathbb{Z}\oplus \mathcal{O}_{k_1}$. 
We suppose that $K$ is a quartic field with $G=C_4$. Then $k_1$ is the unique quadratic subfield of $K$. 
Let $\chi$ be the Dirichlet character of order $4$ corresponding to the cyclic quartic field $K$ 
and denote by $f_\chi>1$ the conductor of $\chi$. By the discriminant-conductor formula, 
$\Disc(K)=f_\chi^2 \Disc(k_1) \neq \Disc(k_1)$. This is a contradiction. 
We next suppose that $K$ is a quartic field with $G=D_4$. 
Then $K$ is a quadratic extension of a quadratic field $k_2$. 
Hence $\Disc(K)=N(\Disc(K/k_2))\Disc(k_2)^2$. 
This contradicts $\Disc(K)=\Disc(k_1)$. 
Therefore $K$ is a direct sum of two distinct quadratic fields or $K=\mathbb{Q}\oplus \mathbb{Q} \oplus k_1$. 
We denote by $\mathcal{Q}(k_1)$ the set of the isomorphism classes of such quartic algebras. 
If $\Disc(k_1)>0$, then we also denote by $\mathcal{Q}_+(k_1)$ the set of the isomorphism classes of such totally real quartic algebras. 
Since $Q=\mathcal{O}_K$, it suffices to count $|\mathcal{Q}(k_1)|$ and $|\mathcal{Q}_+(k_1)|$ for our purpose. 

If $K=\mathbb{Q}\oplus \mathbb{Q} \oplus k_1$, then $\mathcal{O}_K=\mathbb{Z} \oplus \mathbb{Z} \oplus \mathcal{O}_{k_1}$, hence 
$\Disc(K)=\Disc(k_1)$. 
If $K=k_2 \oplus k_3$ where $k_2$ and $k_3$ are two distinct quadratic fields, 
Then $\Disc(k_1)=\Disc(K)=\Disc(k_2)\Disc(k_3)$. 
We put $d_i=\Disc(k_i)$ for $i=1, 2, 3$. 
We denote by $t$ the number of prime divisors of $d_1$. 
We must count the number of expressions of writing $d_1$ as a product of two fundamental discriminants $d_2$ and $d_3$. 
We denote it by $m(d_1)$. If $d_1>0$, we denote by $m_+(d_1)$ the number of expressions of 
writing $d_1$ as a product of two positive fundamental discriminants. 

We use the following theorem of Gauss (cf. \cite[Theorem 3.70]{MOL}). 
\begin{thm}\label{thm:QUAD-CL-2RANK}
Let $k_1$, $d_1$ and $t$ be as above. Then $|\Cl_{k_1,+}^{(2)}|=2^{t-1}$. 
Further $|\Cl_{k_1}^{(2)}|=2^{t-2}$ if $d_1>0$ and $d_1$ has a prime divisor $p\equiv 3 \pmod{4}$,  
otherwise $|\Cl_{k_1}^{(2)}|=2^{t-1}$. 
\end{thm}

We first assume $t=1$. Then $d_1$ is one of the followings: $d_1=p$, $p$ is a prime number with $p \equiv 1 \pmod{4}$; 
$d_1=-p$, $p$ is a prime number with $p \equiv 3 \pmod{4}$; 
$d_1=8$; $d_1=-4$; $d_1=-8$.
These $d_1$'s can not be expressed as a product of two fundamental discriminants. 
So $m(d_1)=0$, hence the isomorphism class of $\mathbb{Q}\oplus \mathbb{Q} \oplus k_1$ is the only element of $\mathcal{Q}(k_1)$. 
By Theorem \ref{thm:QUAD-CL-2RANK}, we have $|\mathcal{Q}(k_1)|=1=|\Cl_{k_1,+}^{(2)}|$. 
We also have $|\mathcal{Q}_+(k_1)|=1=|\Cl_{k_1}^{(2)}|$ if $d_1>0$.

We next assume $t \geq 2$. We write $d_1=\pm 2^{e_0} p_1 \cdots p_r q_1 \cdots q_s$ where 
$p_i$'s and $q_j$'s are distinct prime numbers such that $p_i\equiv 1 \pmod{4}$ and 
$q_j\equiv 3 \pmod{4}$. We put $q_j^*=-q_j$. 

Case 1. $e_0=0$. Since $d_1\equiv 1 \pmod{4}$,  we have $d_1=p_1\cdots p_r q_1^* \cdots q_s^*$ and $t=r+s$. 
It is clear that $d_2=\prod_i p_1^{a_i} \prod_j (q_j^*)^{b_j}$ for some $a_i, b_j \in \{0,1\}$ with $d_2 \neq 1, d_1$. 
Hence $m(d_1)=2^t-2$. Since $k_2 \oplus k_3\cong k_3 \oplus k_2$, 
we have $|\mathcal{Q}(k_1)|=m(d_1)/2+1=2^{t-1}=|\Cl_{k_1,+}^{(2)}|$ taking account of $\mathbb{Q}\oplus \mathbb{Q} \oplus k_1$. 
We assume $d_1>0$, so $s$ is even. If $s=0$, then obviously $\mathcal{Q}_+(k_1)=\mathcal{Q}(k_1)$, hence 
$|\mathcal{Q}_+(k_1)|=2^{t-1}=|\Cl_{k_1}^{(2)}|$ by Theorem \ref{thm:QUAD-CL-2RANK}. 
If $s\geq 2$, then $d_2=\prod_i p_1^{a_i} \prod_j (q_j^*)^{b_j}$ for some $a_i, b_j \in \{0,1\}$ with $d_2 \neq 1, d_1$ 
and $\sum_j b_j \equiv 0 \pmod{2}$. 
Hence $m_+(d_1)=2^{t-1}-2$, $|\mathcal{Q}_+(k_1)|=m_+(d_1)/2+1=2^{t-2}=|\Cl_{k_1}^{(2)}|$ by Theorem \ref{thm:QUAD-CL-2RANK}. 

Case 2. $e_0=3$. Then we have $d_1=\pm 2^3 p_1\cdots p_r q_1^* \cdots q_s^*$ and $t=r+s+1$. 
We may assume that $d_2$ is odd and $d_3$ is even. 
It is clear that $d_2=\prod_i p_1^{a_i} \prod_j (q_j^*)^{b_j}$ for some $a_i, b_j \in \{0,1\}$ with $d_2 \neq 1$. 
Hence we have $|\mathcal{Q}(k_1)|=2^{t-1}=|\Cl_{k_1,+}^{(2)}|$. 
We assume $d_1>0$. If $s=0$, then obviously $\mathcal{Q}_+(k_1)=\mathcal{Q}(k_1)$, hence 
$|\mathcal{Q}_+(k_1)|=2^{t-1}=|\Cl_{k_1}^{(2)}|$ by Theorem \ref{thm:QUAD-CL-2RANK}. 
If $s\geq 1$, then $d_2=\prod_i p_1^{a_i} \prod_j (q_j^*)^{b_j}$ for some $a_i, b_j \in \{0,1\}$ with $d_2 \neq 1$ 
and $\sum_j b_j \equiv 0 \pmod{2}$. 
Hence $|\mathcal{Q}_+(k_1)|=2^{t-2}=|\Cl_{k_1}^{(2)}|$ by Theorem \ref{thm:QUAD-CL-2RANK}.

Case 3. $e_0=2$. Then we must have $d_1=(-1)^{s+1}2^2 p_1\cdots p_r q_1 \cdots q_s$ and $t=r+s+1$. 
We may assume that $d_2$ is odd and $d_3$ is even. 
It is clear that $d_2=\prod_i p_1^{a_i} \prod_j (q_j^*)^{b_j}$ for some $a_i, b_j \in \{0,1\}$ with $d_2 \neq 1$. 
Hence we have $|\mathcal{Q}(k_1)|=2^{t-1}=|\Cl_{k_1,+}^{(2)}|$. 
We assume $d_1>0$, so $s$ is odd. Then $d_2=\prod_i p_1^{a_i} \prod_j (q_j^*)^{b_j}$ for some $a_i, b_j \in \{0,1\}$ with $d_2 \neq 1$ 
and $\sum_j b_j \equiv 0 \pmod{2}$. 
Hence $|\mathcal{Q}_+(k_1)|=2^{t-2}=|\Cl_{k_1}^{(2)}|$ by Theorem \ref{thm:QUAD-CL-2RANK}. 

For any $K \in \mathcal{Q}(k_1)$, it is obvious that $\Aut(\mathcal{O}_K)$ is isomorphic to $(\mathbb{Z}/2\mathbb{Z})^2$. 
So we obtain the following proposition by Theorem \ref{thm:BHARGAVA3}, Corollary \ref{cor:UNIQUE-RESOLVENT}:
\begin{prop}\label{prop:QPLUSQUAD-RESOLVENT}
Let $k_1$ be a quadratic field. We put $k=\mathbb{Q}\oplus k_1$ and $\mathcal{O}_k=\mathbb{Z}\oplus\mathcal{O}_{k_1}$. 
If $\Disc(k_1)>0$, then we have 
\begin{align*}
\sum_{x \in \Gamma\backslash L_1(\mathcal{O}_k)} \mu(x)
 &=\frac{1}{4}|\Cl_{k_1}^{(2)}|, \\
\sum_{x \in \Gamma\backslash L_1(\mathcal{O}_k)} \mu(x)
+\sum_{x \in \Gamma\backslash L_3(\mathcal{O}_k)} \mu(x)
 &=\frac{1}{4}|\Cl_{k_1,+}^{(2)}|.
\end{align*}
If $\Disc(k_1)<0$, then we have 
\[ \sum_{x \in \Gamma\backslash L_2(\mathcal{O}_k)} \mu(x)
 =\frac{1}{4}|\Cl_{k_1}^{(2)}|.
\]
\end{prop}

Let $k_1$, $k$ be as above. We recall that $U^+(\mathcal{O}_k)$ and $U_+(\mathcal{O}_k)$ 
denote the group of units of $\mathcal{O}_k$ having positive norm 
and the group of totally positive units of $\mathcal{O}_k$, respectively. 
By \eqref{eq:UPLUS-INDEX}, we have
\[(U^+(\mathcal{O}_0):U^+(\mathcal{O}_0)^2)=\left\{
\begin{array}{ll}
4 , &\quad  \Disc(k)>0, \\
2 , &\quad \Disc(k)<0,
\end{array}
\right.
\]
Since $U_+(\mathcal{O}_k)=1 \times U_+(\mathcal{O}_{k_1})$, 
we have $(U_+(\mathcal{O}_k):U_+(\mathcal{O}_k)^2)=2$. 
We recall that $U_2(\mathcal{O}_k)$ denotes the group of units in $\mathcal{O}_k$ having order dividing 2. 
Then we have $U_2^+(\mathcal{O}_k)=U_2(\mathcal{O}_k) \cap U^+(\mathcal{O}_k)$. 
Hence we have $U_2(\mathcal{O}_k)=\{(1,1),(1,-1),(-1,1),(-1,-1)\}$ and 
$U_2^+(\mathcal{O}_k)=\{(1,1),(1,-1)\}$. 
Further we have $|\Aut(\mathcal{O}_k)|=2$, $\Cl_k^{(2)}=\Cl_{k_1}^{(2)}$ and $\Cl_{k,+}^{(2)}=\Cl_{k_1,+}^{(2)}$. 
We obtain the following formulae by Proposition \ref{prop:LHAT-CLN1}: 
If $\Disc(k_1)>0$, then we have 
\begin{equation}\label{eq:LHAT-CLN-QPLUSREALQUAD}
\begin{split}
\sum_{y \in \Gamma\backslash \hat{L}_1(\mathcal{O}_k)} \mu(y)
&=\frac{1}{4}|\Cl_{k_1,+}^{(2)}|, \\
\sum_{y \in \Gamma\backslash \hat{L}_1(\mathcal{O}_k)} \mu(y)
+\sum_{y \in \Gamma\backslash \hat{L}_3(\mathcal{O}_k)} \mu(y)
 &=|\Cl_{k_1}^{(2)}|
\end{split}
\end{equation}
and if $\Disc(k)<0$, then we have 
\begin{equation}\label{eq:LHAT-CLN-QPLUSIMAGQUAD}
\sum_{y \in \Gamma\backslash \hat{L}_2(\mathcal{O}_k)} \mu(y)
=\frac{1}{2}\,|\Cl_{k_1}^{(2)}|.
\end{equation}
By Proposition \ref{prop:SPLIT-RESOLVENT2}, \eqref{eq:LHAT-CLN-QPLUSREALQUAD} and \eqref{eq:LHAT-CLN-QPLUSIMAGQUAD}, 
we finally obtain the following proposition:
\begin{prop}\label{prop:SPLIT-RESOLVENT2}
Let $k_1$ be a quadratic field. We put $k=\mathbb{Q}\oplus k_1$ and $\mathcal{O}_k=\mathbb{Z}\oplus\mathcal{O}_{k_1}$. 
If $\Disc(k_1)>0$, then we have 
\begin{align*}
\sum_{y \in \Gamma\backslash \hat{L}_1(\mathcal{O}_k)} \mu(y)
+\sum_{y \in \Gamma\backslash \hat{L}_3(\mathcal{O}_k)} \mu(y)
 &=4\sum_{x \in \Gamma\backslash L_1(\mathcal{O}_k)} \mu(x), \\
\sum_{y \in \Gamma\backslash \hat{L}_1(\mathcal{O}_k)} \mu(y)
&=\sum_{x \in \Gamma\backslash L_1(\mathcal{O}_k)} \mu(x)
+\sum_{x \in \Gamma\backslash L_3(\mathcal{O}_k)} \mu(x).
\end{align*}
If $\Disc(k_1)<0$, then we have 
\[
\sum_{y \in \Gamma\backslash \hat{L}_2(\mathcal{O}_k)} \mu(y)
=2 \sum_{x \in \Gamma\backslash L_2(\mathcal{O}_k)} \mu(x).
\]
\end{prop}

Let $n\in \mathbb{Z}$ be a fundamental discriminant. 
Theorem \ref{thm:MAIN} now follows from Proposition \ref{prop:SPLIT-RESOLVENT2} and Theorem \ref{thm:MAIN-S4A4} 
applied to the non-Galois cubic fields $k$ having discriminant $n$ 
and $\mathcal{O}=\mathcal{O}_k$.

\end{document}